\documentclass[11pt,leqno]{article}

\usepackage{amsmath,amsfonts,amscd,amssymb,theorem}

\long\def\comment#1\endcomment{}

\comment
\endcomment

\makeatletter
\begingroup
\gdef\th@dotted{\normalfont\itshape
  \def\@begintheorem##1##2{%
        \item[\hskip\labelsep \theorem@headerfont ##1\ ##2.]}%
\def\@opargbegintheorem##1##2##3{%
   \item[\hskip\labelsep \theorem@headerfont ##1\ ##2\ (##3).]}}
\endgroup
\makeatother

\theoremstyle{dotted}

\newtheorem{theorem}{Theorem}[section]
\newtheorem{lemma}[theorem]{Lemma}

\newtheorem{prop}[theorem]{Proposition}
\newtheorem{corr}[theorem]{Corollary}

\makeatletter
\begingroup
\gdef\th@upshape{\normalfont
  \def\@begintheorem##1##2{%
        \item[\hskip\labelsep \theorem@headerfont ##1\ ##2.]}%
\def\@opargbegintheorem##1##2##3{%
   \item[\hskip\labelsep \theorem@headerfont ##1\ ##2\ (##3).]}}
\endgroup
\makeatother

\theoremstyle{upshape}

\newtheorem{defn}[theorem]{Definition}
\newtheorem{remark}[theorem]{Remark}
\newtheorem{exa}[theorem]{Example}

\makeatletter
\renewcommand{\subsection}{\@startsection{subsection}{2}{0pt}{-3ex
plus -1ex minus -0.2ex}{-2mm plus -0pt minus
-2pt}{\normalfont\bfseries}} 
\renewcommand{\subsubsection}{\@startsection{subsubsection}{3}{0pt}{-3ex
plus -1ex minus -0.2ex}{-2mm plus -0pt minus
-2pt}{\normalfont\bfseries}} 
\makeatother

\makeatletter
\@addtoreset{equation}{section}
\makeatother

\newcommand{\cntrct}                
{\hspace{2pt}\raisebox{1pt}{\text{$\lrcorner$}}\hspace{2pt}}

\newcommand{\proof}[1][Proof.]{\smallskip\noindent{\em #1}}
\def\endproof{\hfill\ensuremath{\square}\par\medskip}

\renewcommand{\labelenumi}{{\normalfont(\roman{enumi})}}

\def\eqref#1{\thetag{\ref{#1}}}

\let\latexref=\ref
\def\ref#1{{\normalfont{\latexref{#1}}}}

\newcommand{\wt}{\widetilde}
\newcommand{\wh}{\widehat}

\newcommand{\idot}{{\:\raisebox{1pt}{\text{\circle*{1.5}}}}}
\newcommand{\hdot}{{\:\raisebox{3pt}{\text{\circle*{1.5}}}}}
\newcommand{\Sp}{\operatorname{\sf Sp}}

\newcommand{\eps}{\varepsilon}
\renewcommand{\phi}{\varphi}

\newcommand{\vH}{\check{H}}

\def\dlim_#1{{\displaystyle\lim_{#1}}^\hdot}

\newcommand{\Hom}{\operatorname{Hom}}

\newcommand{\RHom}{\operatorname{RHom}}

\newcommand{\id}{\operatorname{\sf id}}
\newcommand{\Id}{\operatorname{\sf Id}}
\newcommand{\gr}{\operatorname{\sf gr}}

\newcommand{\A}{{\cal A}}
\newcommand{\D}{{\cal D}}
\newcommand{\DM}{{\cal D}{\cal M}}
\newcommand{\DN}{{\cal D}{\cal N}}
\newcommand{\wDM}{\widehat{\DM}}
\newcommand{\DS}{{\cal D}{\cal S}}

\newcommand{\C}{{\cal C}}
\newcommand{\B}{{\cal B}}

\newcommand{\Q}{{\cal Q}}

\newcommand{\Sets}{\operatorname{Sets}}
\newcommand{\Supp}{\operatorname{Supp}}

\newcommand{\Iso}{{\operatorname{Iso}}}
\newcommand{\Aut}{{\operatorname{Aut}}}

\newcommand{\Add}{\operatorname{\sf Add}}

\newcommand{\amod}{{\text{\rm -mod}}}

\newcommand{\ppt}{{\sf pt}}

\newcommand{\lotimes}{\overset{\sf\scriptscriptstyle L}{\otimes}}

\newcommand{\bPhi}{\overline{\Phi}}

\newcommand{\wPsi}{\widetilde{\Psi}}
\newcommand{\wPhi}{\widetilde{\Phi}}

\newcommand{\M}{\operatorname{\mathcal M}}
\newcommand{\N}{\operatorname{\mathcal N}}
\newcommand{\Nn}{\operatorname{\sf N}}

\newcommand{\Infl}{\operatorname{\sf Infl}}
\newcommand{\Ind}{\operatorname{\sf Ind}}

\newcommand{\copr}{\sqcup}

\newcommand{\Z}{{\mathbb Z}}
\newcommand{\QQ}{{\mathbb Q}}
\newcommand{\RR}{{\mathbb R}}

\newcommand{\W}{\mathbb{W}}

\newcommand{\hh}{\mathcal{H}\mathcal{H}}

\newcommand{\Fun}{\operatorname{Fun}}

\newcommand{\E}{\mathcal{E}}

\newcommand{\wM}{\operatorname{\wh{\mathcal{M}}}}
\newcommand{\wA}{\operatorname{\wh{\mathcal{A}}}}
\newcommand{\wAG}{\operatorname{\wh{\mathcal{A}^G}}}
\newcommand{\wBG}{\operatorname{\wh{\mathcal{B}^G}}}
\newcommand{\wAZ}{\operatorname{\wh{\mathcal{A}^{\Z}}}}
\newcommand{\wGamma}{\operatorname{\widehat{\Gamma}}}

\newcommand{\Inj}{\operatorname{{\sf Inj}}}
\newcommand{\Av}{\operatorname{{\sf Av}}}

\newcommand{\vC}{\check{C}}

\newcommand{\DA}{\D^\alpha}

\newcommand{\hGamma}{\overline{\Gamma}}

\newcommand{\Tel}{\operatorname{Tel}}
\newcommand{\tel}{\operatorname{tel}}

\newcommand{\R}{\mathcal{R}}

\title{Mackey profunctors}

\author{D. Kaledin}

\begin{document}

\maketitle

\tableofcontents

\section*{Introduction.}

The so-called ``Mackey functors'' associated to a finite group $G$
have long been a standard tool both in group theory and in algebraic
topology, where they appear as natural coefficients for
$G$-equivariant cohomology theories.  The reader can find modern
introductions to Mackey functor in the topological context e.g. in
\cite{may1}, \cite{may2}, \cite{tD}, or a more algebraic treatment
in \cite{the}. A derived version of the theory has been suggested in
\cite{mackey}.

In topological applications, the group $G$ does not have to be
finite --- one can generalize the notion of a Mackey functor to
allow $G$ to be an arbitrary compact Lie group equipped with its
standard topology.

The goal of the present paper is to suggest another
generalization. In our theory, the group $G$ is infinite but
discrete --- the only requirement is that it is finitely generated.

We should note that formally, the definition of a $G$-Mackey functor
depends not on the group $G$ but on the category of finite
$G$-sets. Thus already in the original definition, one can allow $G$
to be infinite --- it is only $G$-sets we consider that have to be
finite. However, the finer points of the theory break down, and the
resulting category of $G$-Mackey functors seems not to be the right
thing to consider. What we suggest in the present paper is an
alternative theory of ``Mackey profunctors'' that seems to behave
better and preserve most of the nice properties of Mackey functors
for finite groups. The basic idea is to relax the finiteness
requirement on $G$-sets, and allow $G$-sets $S$ such that the
stabilizer $H \subset G$ of any element $s \in S$ is cofinite, while
conversely, the fixed point subset $S^H \subset S$ is finite for any
cofinite subgroup $G$ (this is our Definition~\ref{adm.G.set}). We
call such $G$-sets {\em admissible}, motivated by a similar notion
of an admissible representation of a $p$-adic Lie group of
\cite{ber}, and we show that the resulting theory is a ``profinite
completion'' of the usual theory. Roughly speaking, giving a
$G$-Mackey profunctor is equivalent to giving a system of $W$-Mackey
profunctors for all finite quotients $W$ of the group $G$, related
by some natural compatibility maps. This allows to transfer most of
the finite theory to our profinite case without much effort.

We note that both the naive generalization of Mackey functors and
our theory of Mackey profunctors only depend on the profinite
completion $\wh{G}$ of the group $G$. This is still a compact group
with respect to its natural profinite topology. However, the
topology is very different from the standard topology on Lie groups,
and the theory we develop also looks differently --- in effect, it
is completely orthogonal to the theory of Mackey functors for
compact Lie groups. At present, we do not know whether there is a
general ``adelic'' context that unifies the two.

\medskip

Formally, the theory developed in the paper is self-contained, and
the whole paper could be considered simply an extended exercise in
homological algebra. However, given the length of the exercise, we
should say at least a couple of words about motivations.

\medskip

The main immediate application of the theory that we have in mind
concerns one specific infinite group, namely, the infinite cyclic
group $\Z$. The Mackey profunctors that we are after occur in the
theory of Topological Cyclic Homology of \cite{BHM} and its recent
algebraic version, Hochschild-Witt Homology of \cite{ka-wi},
\cite{ka-hw}.

The starting point of the story is the well-known relation between
Mackey functors and Witt vectors. Its simplest instance is the
isomorphism $W_n(\Z) \cong A(\Z/p^n\Z)$ between the ring $W_n(\Z)$
of $n$-truncated $p$-typical Witt vectors of the ring $\Z$ and the
Burnside ring $A(\Z/p^n\Z)$ of the cyclic group $\Z/p^n\Z$. To get a
similar interpretation of the universal Witt vectors ring $\W(\Z)$
one needs to define correctly Burnside rings of profinite groups,
and this has been successfully accomplished in \cite{dr-si}, in
essentially the same way as in our paper (our admissible $G$-sets
appear in \cite[2.2]{dr-si} as {\em almost finite} $G$-sets).

The theory of Topological Cyclic Homology contains both Witt vectors
and Mackey functors. Witt vectors are in fact an integral part of
the story, to the point that L. Hesselholt \cite{hewi} used
Topological Cyclic Homology to give a generalization of Witt vectors
to non-commutative rings. Mackey functors are even more central,
since the main technical tool of \cite{BHM} is equivariant stable
homotopy theory, with the group in question being the unit circle
$S^1$. In practice, what one considers is not the full unit circle
but all its finite subgroups --- that is, the subgroups $\Z/n\Z
\subset S^1$ formed by roots of unity of order $n$, $n \geq 1$. One
observes that for any ring spectrum $A$, the Topological Hochschild
Homology spectrum $THH(A)$ of B\"okstedt \cite{bo} is naturally an
$S^1$-equivariant spectrum, and observes that fixed point spectra
$THH(A)^{\Z/n\Z}$ form an inverse system. Its limit is denoted
$TR(A)$. Hesselholt's non-commutative Witt vectors $\W(A)$ of an
associative ring $A$ are by definition elements in $\pi_0(TR(A))$.

In the existing theory, $TR(A)$ is only a spectrum with no obvious
equivariant structure. However, one can package the same groups
$\Z/n\Z$ differently --- they are also finite quotients of the
infinite cyclic group $\Z$. Moreover, the fixed point spectrum
$THH(A)^{\Z/n\Z}$ is naturally $\Z/n\Z$-equivariant for any $n \geq
1$, and once one has a good equivariant homotopy theory for infinite
groups, $TR(A)$ should be $\Z$-equivariant. Considering it as such,
we keep track of all the information that went into construction,
and in particular, we can recover back all the fixed points spectra
$THH(A)^{\Z/n\Z}$.

The big advantage of this point of view is that it should allow one to
develop the theory with coefficients. Namely, the usual Hochschild
homology $HH(A)$ of an associative algebra $A$ can be extended to a
two-variable homology theory $HH(A,M)$ of an algebra $A$ with
coefficients in a bimodule $M$, and this additional flexibility
helps very much in practical computations. Similarly, one can define
$THH(A,M)$, but this will not be a $S^1$-equivariant spectrum
anymore, so the constructions of \cite{BHM} would not
apply. However, it seems plausible that one can define $TR(A,M)$
directly as a $\Z$-equivariant spectrum, and recover all the data of
\cite{BHM} starting from the other end.

One situation where this has actually been done is the
Hochschild-Witt homology theory of \cite{ka-hw}, an algebraic
counterpart of \cite{BHM}. Here $A$ is an associative unital algebra
over a finite field $k$ of characteristic $p$. It is known
(\cite{HM}) that in this situation, $TR(A)$ is an Eilenberg-Mac Lane
spectrum, thus effectively a complex of abelian groups, and the
Hochschild-Witt complex $WCH_\idot(A)$ of \cite{ka-hw} is expected
to be quasiisomorphic to the $p$-typical part $TR(A;p)$ of $TR(A)$
(this is proved in degree $0$, where both give Hesselholts's
non-commutative Witt vectors, and also in the case when $A$ is
commutative and smooth, where both theories give the de Rham-Witt
differential forms of Deligne and Illusie \cite{ill}). By its very
construction, $WCH_\idot(A)$ comes equipped with an extension to a
two-variable theory $WCH_\idot(A,M)$, and it also explicitly has a
structure of a Mackey profunctor with respect to an infinite group.

Since the theory in \cite{ka-hw} is $p$-typical, the group in
question is not $\Z$ but its $p$-adic completion $\Z_p$. In this
case, the theory of Mackey profunctors is much simpler, so that one
can explicitly describe them by a small amount of data (and this is
how it is done in \cite{ka-hw}). Defining and studying the universal
counterpart $\W CH_\idot(A)$ of the Hochschild-Witt complex that
corresponds to the entire $TR(A)$ would require the full theory of
Mackey profunctors as developed in this paper.

\medskip

Given that we are mostly concerned with this set of applications, we
could have restricted ourselves to the case $G = \Z$. The reason we
chose to do things in full generality is twofold.

\medskip

Firstly, the general case is not essentially different from the
special case we really need. The trivial case for the theory would
be $\Z_p$, the group of $p$-adic integers; all its finite quotients
are just cyclic groups, and the lattice of cofinite subgroups is
simply the lattice of positive integers. One can show that in such a
simple case, our Mackey profunctors reduce to Mackey functors. For
$\Z$, the finite quotients are still cyclic groups, but the subgroup
lattice is much bigger, and in fact is not finitely generated. As
the result of this, all the finer points in the theory are as
non-trivial for $\Z$ as for any finitely generated group $G$, and
the fact that $\Z$ is commutative only obscures things (although
after the general theory has been developed, one can prove further
results valid only for $\Z$).

\medskip

Secondly, we do have further applications in mind, although at this
point they are rather speculative.

\medskip

One potential application concerns the study of Hochschild
cohomology. For algebras over $\QQ$, Hochschild homology and
Hochschild cohomology come as a package, and both carry additional
structures: the Connes-Tsygan differential $B$ on $HH_\idot(A)$, and
the Gerstenhaber bracket on $HH^\hdot(A)$ (see e.g.\ \cite{TT} for
an overview). The differential $B$ is best described and constructed
in terms of A. Connes' cyclic category $\Lambda$ whose objects
correspond to configurations of points on a circle. For
$HH^\hdot(A)$, the Gerstenhaber bracket comes as a part of a
structure of $E_2$-algebra, and this is controlled by configuration
spaces of points on a $2$-sphere $S^2$. In positive characteristic,
and over $\Z_p$, we have all this but also much more. On
$HH_\idot(A)$, we have higher ``motivic'' structures in the sense of
\cite{ka-icm} such as the filtered Dieudonn\'e module
structure. These can be again described by a version of the category
$\Lambda$, except that now one has to also consider $n$-fold
coverings of the circle by itself. The resulting category is
essentially equivalent to the category $\Z$-Mackey profunctors
equipped with some additional compatibility with the cyclic category
$\Lambda$.  What happens for $HH^\hdot(A)$, however, is completely
unknown at the moment (except for a trivial observation that
$HH^1(A)$ with the Gerstenhaber bracket is a restricted Lie
algebra). It is natural to expect that the relevant combinatorics is
then encoded by coverings of the punctured $2$-sphere. But now the
group in question is a free group on several generators, not $\Z$,
and one would need a theory of Mackey profunctors at least for free
groups --- and possibly also for the braid groups and the mapping
class groups.

Another set of applications is arithmetic, with $G$ being the Galois
group of a number field or, more generally, an \'etale fundamental
group of a scheme. Here we have Voevodsky's theory of Artin motives
and ``sheaves with transfer'' \cite{V}. This is structurally similar
to Mackey functors, but there is one difference: instead of allowing
all correspondences between finite $G$-sets $S_1$, $S_2$ as maps,
one only considers correspondences $S$ such that the map $S \to S_1
\times S_2$ in injective. It is reasonable to expect that the full
category of $G$-Mackey functors also plays a role in the story, and
then it would be easier and more natural to work with the whole
Galois group and Mackey profunctors. However, it seems that so far,
this has not been studied at all.

Finally, let us mention that while in this paper, we work entirely
within the framework of homological algebra, a similar treatment of
equivariant stable homotopy theory is certainly also possible, and
this would result in a good equivariant stable homotopy category for
a profinite group. One place where such a treatment was attempted is
recent preprint \cite{clark} by C. Barwick. What are the topological
motivations for such a theory, or even whether there are any, I am
not qualified to say.

\bigskip

Let us now give a brief overview of the paper. As we have explained,
we start with a completely arbitrary group $G$. At some point in the
story, we need to require it to be finitely generated, but this is
the only restriction. We develop profinite counterparts both for the
classic theory of Mackey functors and for its derived version
constructed in \cite{mackey}. Section~\ref{intro.sec} contains
preliminaries and notation. Section~\ref{mack.class.sec} is a brief
recapitulation of the standard theory of Mackey
functors. Section~\ref{mackey.sec} introduces $G$-Mackey profunctors
and contains the beginning of the theory, up to the point when we
need to go to the derived setting. In particular, we introduce the
notion of a ``normal system'', an axiomatization of a compatible
family of $W$-Mackey functors for all the finite quotients $W$ of
the group $G$. To go to the derived setting, we use some technology
developed in \cite{mackey}. In Section~\ref{S.sec}, we recall this
technology and we develop it further (in particular,
Subsection~\ref{bc.subs} is new). Section~\ref{add.sec} is also new
--- while its main result is a version of a result proved in
\cite{mackey}, we give an alternative proof that is easier and
cleaner. In Section~\ref{der.sec}, we start developing the theory of
derived Mackey profunctors. We use both the technology of
\cite{mackey} and new tools created for this paper. These new tools
also allow to clean up and strengthen some of the results of
\cite{mackey}; we take the opportunity to do so in
Section~\ref{repr.sec}. Then in Section~\ref{ns.sec}, we continue
with derived Mackey profunctors, and we prove the results that
depend on a derived version of the notion of a normal
system. Finally, in Section~\ref{mack.cycl.subs}, we show how our
abstract machinery works in the particular case $G=\Z$, the infinite
cyclic group.

\subsection*{Acknowledgements.} I am very grateful to L. Hesselholt
who, among other things, introduced me to Mackey functors, and to
all the people with whom I had an opportunity to discuss this
subject. I am grateful to the referees for very careful reports and
many useful suggestions, and in particular, for bringing the paper
\cite{dr-si} to my attention. A large part of this work was done in
April 2014, during my visit to J. Dieudonn\'e Laboratory of the
University of Nice-Sophia Antipolis, with its wonderful people and
great working atmosphere. I am very grateful to E. Balzin and
C. Simpson for setting up the visit, and for many useful
discussions. In fact, the only distractions during my visit were
continuous news about my home country behaving as a fascist state it
has become, reaching new lows each week. I am grateful to the people
of Ukraine for their courage and humanity, since this is the only
thing that stands between us and complete darkness, and I wish them
total victory in the struggle against their barbarian neighbor.

\section{Preliminaries.}\label{intro.sec}

We start with recalling some facts about combinatorics of simplicial
sets and homology of small categories that we will need in the
paper. All the material is quite standard; we mostly include it to
fix notation.

\subsection{Homology of small categories.}

For any two objects $c,c' \in \C$ of a category $\C$, we will denote
by $\C(c,c')$ the set of maps in $\C$ from $c$ to $c'$. For any
object $c \in \C$ in a category $\C$, we will denote by $\C/c$ the
category of objects $\wt{c} \in \C$ equipped with a map $\wt{c} \to
c$. Any map $f:c \to c'$ in $\C$ induces a functor
\begin{equation}\label{f.l}
f_!:\C/c \to \C/c'
\end{equation}
sending a map $\wt{c} \to c$ to its composition with $f$, and if
$\C$ has fibered products, then this functor has a right-adjoint
\begin{equation}\label{f.r}
f^*:\C/c' \to \C/c
\end{equation}
sending a map $\wt{c} \to c'$ to the natural projection $\wt{c}
\times_{c'} c \to c$.

For any category $\C$, an {\em inverse system} $\{c_i\}$ of objects
in $\C$ is a collection of objects $c_i \in \C$, $i \geq 1$,
equipped with transition maps $c_{i+1} \to c_i$.

For any small category $\C$ and ring $R$, we denote by $\Fun(\C,R)$
the category of functors from $\C$ to the category of
$R$-modules. For any $R$-module $M$, we will denote by $M_{\C} \in
\Fun(\C,R)$ the constant functor with value $M$, and we will
sometimes shorten it to $M$ if $\C$ is clear from the context. For
any object $c \in \C$, we will denote by $M_c \in \Fun(\C,R)$ the
functor given by
\begin{equation}\label{m.c.def}
M_c(c') = M[\C(c,c')],
\end{equation}
the sum of copies of $M$ numbered by elements in the set of maps
$\C(c,c')$.

The category $\Fun(\C,R)$ is abelian. We denote its derived category
by $\D(\C,R)$. If $\C$ is the point category $\ppt$, so that
$\Fun(\ppt,R)$ is the category of $R$-modules, we shorten the
notation by setting $\D(R) = \D(\ppt,R)$, the derived category of
$R$-modules. Objects $M_c$ of \eqref{m.c.def} generate the category
$\D(\C,R)$ in the sense that every object $E \in \D(\C,R)$ can be
represented by an $h$-projective complex whose terms are sums of
objects of the form \eqref{m.c.def}.

For any functor $\gamma:\C \to \C'$ between small categories,
composition with $\gamma$ gives a pullback functor
$\gamma^*:\Fun(\C',R) \to \Fun(\C,R)$. It has a left and a
right-adjoint functor known as the left and the right Kan extension;
we denote them by $\gamma_!,\gamma_*:\Fun(\C,R) \to
\Fun(\C',R)$. The derived functors
$L^\hdot\gamma_!,R^\hdot\gamma_*:\D(\C,R) \to \D(\C',R)$ are left
resp. right-adjoint to the pullback functor $\gamma^*:\D(\C',R) \to
\D(\C,R)$. If $\C'$ is the point category $\ppt$, and $\tau:\C \to
\ppt$ is the tautological projection, then the Kan extension
functors $\tau_!$, $\tau_*$ are just the colimit and limit over the
small category $\C$. We denote
$$
C_\idot(\C,E) = L^\hdot\tau_!E \in \D(R)
$$
for any $E \in \Fun(\C,R)$. The homology groups $H_\idot(\C,E)$ of
this complex are by definition the homology groups of the category
$\C$ with coefficients in $E$.

In some situations, computing the derived Kan extensions is
easy. For example, for any $c \in \C$, $\gamma:\C \to \C'$ and
$R$-module $M$, we have
\begin{equation}\label{mc.sh}
L^\hdot \gamma_!M_c \cong M_{\gamma(c)},
\end{equation}
where $M_c$, $M_{\gamma(c)}$ are the functors
\eqref{m.c.def}. Another case is a functor $\gamma:\C \to \C'$ that
admits a right-adjoint $\delta:\C' \to \C$. Then we have a natural
isomorphism
$$
\gamma_! \cong \delta^*,
$$
and in particular, $\gamma_!$ is an exact functor, so that we also
have an isomorphism
\begin{equation}\label{adj.pb}
L^\hdot\gamma_! \cong \delta^*
\end{equation}
on the level of derived categories.

Given a small category $\C$, two algebras $R_1$, $R_2$ flat over a
commutative ring $k$, and two objects $E_1 \in \Fun(\C,R_1)$, $E_2
\in \Fun(\C,R_2)$, we denote by
$$
E_1 \otimes_k E_2 \in \Fun(\C,R_1 \otimes_k R_2)
$$
their pointwise tensor product, and we use the same notation for the
derived pointwise tensor product of objects in the derived
categories $\D(\C,R_1)$, $\D(\C,R_2)$. If $k$ is clear from context,
we drop it from notation. If we have two small categories $\C_1$,
$\C_2$, and objects $E_1 \in \D(\C,R_1)$, $E_2 \in \D(\C,R_2)$, then
their {\em box product} is given by
$$
E_1 \boxtimes_k E_2 = \pi_1^*E_1 \otimes_k \pi_2^*E_2 \in \D(\C_1
\times \C_2,R_1 \otimes_k R_2),
$$
where $\pi_1$, $\pi_2$ are projections from $\C_1 \times \C_2$ to
$\C_1$ resp. $\C_2$. Again, if $k$ is clear from context, we drop it
from notation.

Assume given two small categories $\C$, $\C_1$, a ring $R$, and an
object $T \in \Fun(\C_1,\Z)$. Define a functor $j^T_\C:\D(\C,R) \to
\D(\C \times \C_1,R)$ by setting
\begin{equation}\label{j.c.t}
j^T_\C(E) = E \boxtimes T.
\end{equation}
If we have another small category $\C'$ and a functor $\gamma:\C'
\to \C$, we have an obvious isomorphism
\begin{equation}\label{j.c.t.g}
(\gamma \times \id)^* \circ j^T_\C \cong j^T_{\C'} \circ \gamma^*.
\end{equation}
We will need the following easy result.

\begin{lemma}\label{j.c.le}
Assume that for any object $c \in \C_1$, $T(c)$ is a finitely
generated flat $\Z$-module. Then for any small category $\C$, the
functor $j^T_\C$ of \eqref{j.c.t} has a left-adjoint $l^T_\C:\D(\C
\times \C_1,R) \to \D(\C,R)$, and for any functor $\gamma:\C' \to
\C$ between small categories, the base change map
\begin{equation}\label{j.c.bc}
l^T_{\C'} \circ (\gamma \times \id)^* \to \gamma^* \circ l^T_\C
\end{equation}
induced by the isomorphism \eqref{j.c.t.g} is itself an isomorphism.
\end{lemma}

\proof{} Since the category $\D(\C \times \C_1,R)$ is generated by
objects $M_{c \times c_1}$ of \eqref{m.c.def}, $c \times c_1 \in \C
\times \C_1$, $M$ an $R$-module, it suffices to construct $l^T_\C$
for such objects. By adjunction, it is given by
\begin{equation}\label{l.t.eq}
l^T_\C(M_{c \times c_1}) = M_c \otimes T(c_1).
\end{equation}
To prove that \eqref{j.c.bc} is an isomorphism, it again suffices to
check it after evaluating at an object $M_{c \times c_1}$, and
moreover, it suffices to consider the case when $\C' = \ppt$ is a
point category. In this case, the claim immediately follows from
\eqref{m.c.def} and \eqref{l.t.eq}.
\endproof

\subsection{Fibrations and cofibrations.}

Throughout the paper, we will heavily use the machinery of
\cite{SGA}. To fix the terminology, here are the basic ingredients.
\begin{itemize}
\item A morphism $f:c_1 \to c_2$ in a category $\C'$ is {\em
  cartesian} with respect to a functor $\gamma:\C' \to \C$ if any
  $f':c_1' \to c_2$ with $\pi(f) = \pi(f')$ factors uniquely as $f'
  = p \circ f$ with invertible $\pi(p)$.
\item A functor $\pi:\C' \to \C$ is a {\em prefibration} if for any
  $c \in \C'$ and any morphism $f:b \to \pi(c)$ in $\C$, there
  exists a cartesian $f':b' \to c$ with $\pi(f') = f$.
\item A prefibration is a {\em fibration} if the composition of any
  two cartesian maps is cartesian.
\item For any map $f:b' \to b$ in $\C$, the associated {\em
  transition functor} $f^*:\pi^{-1}(b) \to \pi^{-1}(b')$ between
  fibers of a fibration $\pi:\C' \to \C$ sends $c \in \pi^{-1}(b)$
  to the source $c'$ of the cartesian map $f':c' \to c$ with
  $\pi(f') = f$ (one checks that this is functorial, and $f^*$ is
  well-defined up to a canonical isomorphism).
\item For any two fibrations $\gamma_1:\C_1 \to \C$, $\gamma_2:\C_2
  \to \C$, a functor $F:\C_1 \to \C_2$ such that $\gamma_1 \cong
  \gamma_2 \circ F$ is {\em cartesian} if it sends cartesian maps to
  cartesian maps.
\end{itemize}
The notions of a cocartesian map, a cofibration $\pi:\C' \to \C$, a
transition functor $f_*:\pi^{-1}(b') \to \pi^{-1}(b)$, and a
cocartesian functors are obtained by passing to opposite categories.

For any cartesian square
\begin{equation}\label{cof.c}
\begin{CD}
\C_1' @>{\gamma'_1}>> \C'\\
@V{\nu_1}VV @VV{\nu}V\\
\C_1 @>{\gamma_1}>> \C
\end{CD}
\end{equation}
of small categories, if $\gamma$ is a fibration or a cofibration,
then $\gamma'$ is also a fibration resp. a cofibration. Moreover,
for any cartesian functor $F:\C_1 \to \C_2$ between fibrations
$\gamma_1:\C_1 \to \C$, $\gamma_2:\C_2 \to \C$, and an arbitrary
functor $\nu:\C' \to \C$, we can form the commutative diagram
\begin{equation}\label{cof.F}
\begin{CD}
\C_1' @>{F'}>> \C_2' @>{\gamma_2'}>> \C'\\
@V{\nu_1}VV @V{\nu_2}VV @VV{\nu}V\\
\C_1 @>{F}>> \C_2 @>{\gamma_2}>> \C
\end{CD}
\end{equation}
with cartesian squares. Then $\gamma_2'$, $\gamma_1' = F' \circ
\gamma_2'$ are fibrations, the functor $F'$ is cartesian, and for
any ring $R$, the isomorphism $F^{'*} \circ \nu_2^* \cong \nu_1^*
\circ F^*$ induces by adjunction a base change map
\begin{equation}\label{b.c.cart}
\nu_2^* \circ R^\hdot F_* \to R^\hdot F'_* \circ \nu_1^*.
\end{equation}
Dually, if $\gamma_1$, $\gamma_2$ are cofibrations, and $F$ is
cocartesian, we have the base change map
\begin{equation}\label{bc.iso}
L^\hdot F'_! \circ \nu_1^* \to \nu_2^* \circ L^\hdot F_!.
\end{equation}
We prove the following easy result for which we could not find a
convenient reference.

\begin{lemma}\label{bc.sga}
If $F$ in \eqref{cof.F} is a cartesian functor between fibrations,
then the base change map \eqref{b.c.cart} is an isomorphism, and if
$F$ is a cocartesian functor between cofibrations, then
\eqref{bc.iso} is an isomorphism.
\end{lemma}

\proof{} We will prove the claim for cofibrations; the proof for
fibrations is dual. First of all, it clearly suffices to prove the
claim when $\C' = \ppt$, and $\nu$ is the embedding onto an object
$c' \in \C$. Next, the category $\Fun(\C_1,R)$ is generated by
functors \eqref{m.c.def}, so it suffices to prove that for any $c
\in \C_1$ and any $R$-module $M$, the natural map
$$
\nu_2^*L^\hdot F_!M_c \to L^\hdot F_!'\nu_1^*M_c
$$
is an isomorphism. By \eqref{mc.sh}, we have $L^\hdot F_!M_c \cong
M_{F(c)}$, and since $\gamma_1$, $\gamma_2$ are cofibrations,
\eqref{m.c.def} gives canonical identifications
$$
\nu_1^*M_c \cong \bigoplus_{f \in \C(\gamma_1(c),c')}M_{f_*(c)},\qquad
\nu_2^*M_{F(c)} \cong \bigoplus_{f \in
  \C(\gamma_1(c),c')}M_{f_*(F(c))}.
$$
Since $F$ is cocartesian, we have $f_*(F(c)) \cong F(f_*(c))$, so
that to finish the proof, it remains to notice that $L^\hdot
F'_!M_{f_*(c)} \cong M_{F(f_*(c))}$ by \eqref{mc.sh}.
\endproof

In particular, by taking $\C_2 = \C$, we can apply
Lemma~\ref{bc.sga} to any cartesian diagram \eqref{cof.c}. This
shows that for any cofibration $\gamma:\C' \to \C$ and any $E \in
\D(\C',R)$, the value of the Kan extension $L^\hdot\gamma_!E$ at
some object $c \in \C$ can be expressed as
$$
L^\hdot\gamma_!E(c) \cong H_\idot(\gamma^{-1}(c),E),
$$
where $\gamma^{-1}(c) \subset \C'$ is the fiber of the cofibration
$\gamma$ over $c$.

If $\gamma:\C' \to \C$ is not a cofibration, then it can still be
factorized as
\begin{equation}\label{facto.co}
\begin{CD}
\C' @>{e}>> \gamma\backslash\C @>{t}>> \C,
\end{CD}
\end{equation}
where $\gamma\backslash\C$ is the comma-category of the functor $\gamma$
(that is, the category of triples $\langle c',c,f \rangle$, $c' \in
\C'$, $c \in \C$, $f:\gamma(c') \to c$), $t$ is the natural
projection sending $\langle c',c,f \rangle$ to $c$, and $e$ is the
natural embedding sending $c'$ to $\langle c',\gamma(c),\id
\rangle$. Then $t$ is a cofibration, while $e$ has a right-adjoint
$s:\gamma\backslash\C \to \C'$ sending $\langle c',c,f \rangle$ to
$c'$. Therefore by \eqref{adj.pb}, we have
$$
L^\hdot\gamma_! \cong L^\hdot t_! \circ L^\hdot e_! \cong L^\hdot
t_! \circ s^*.
$$
In particular, if $M$ is an $R$-module, and $M_{\C'} \in
\Fun(\C',R)$ is the constant functor with value $M$, then
$s^*M_{\C'} \cong M_{\gamma\backslash\C} \in \Fun(\gamma\backslash\C,R)$, and we
have
\begin{equation}\label{facto.t}
L^\hdot\gamma_!M \cong L^\hdot t_!M.
\end{equation}
For example, let $\C' = \ppt$, the point category, embedded by
$\gamma$ to an object $c \in \C$. Then \eqref{facto.co} reads as
$$
\begin{CD}
\ppt @>>> c\backslash \C @>{t}>> \C,
\end{CD}
$$
where $c \backslash \C$ is the category of objects $c' \in \C$
equipped with a map $c \to c'$, and $t$ is the tautological
projection forgetting the map. In this case, the cofibration $t$ is
discrete, so that $t_!$ is an exact functor, and we have
\begin{equation}\label{m.c.eq}
L^\hdot\gamma_!M \cong \gamma_! M \cong t_!M \cong M_c,
\end{equation}
where $M_c$ is the  functor \eqref{m.c.def}.

\subsection{Pointed sets.}

We will denote by $\Sets$ the category of sets, with $\emptyset$
being the empty set and $\ppt$ being the one-element set. A {\em
  pointed set} is a set $S$ equipped with a distinguished element $o
\in S$. We will denote the category of pointed sets by
$\Sets_+$. For any set $S$ equipped with a subset $A \subset S$, we
denote by $S/A \in \Sets_+$ the pointed set obtained by the pushout
square
$$
\begin{CD}
A @>>> S\\
@VVV @VVV\\
\ppt @>{i}>> S/A,
\end{CD}
$$
with the distiguished element $i(\ppt) \in S/A$. For example, for
any set $S$, the quotient $S_+ = S/\emptyset$ is given by
\begin{equation}\label{s.pl}
S_+ =  S \copr \{o\},
\end{equation}
that is, the union of $S$ with a new element $o$ which we take as
distinguished. The functor $\Sets \to \Sets_+$ sending $S$ to $S_+$
is left-adjoint to the forgetful functor $\Sets_+ \to \Sets$. The
category $\Sets_+$ has finite coproducts given by
\begin{equation}\label{vee.eq}
S \vee S' = (S \copr S')/\{o,o'\},
\end{equation}
and we have $(S \copr S')_+ = S_+ \vee S'_+$ for any two sets $S,S'
\in \Sets$. The {\em smash product} $S \wedge S'$ of two pointed
sets $S,S' \in \Sets_+$ is given by
\begin{equation}\label{sm.eq}
S \wedge S' = (S \times S')/((\{o\} \times  S')\cup(S \times
\{o'\})),
\end{equation}
and we have $(S \times S')_+ = S_+ \wedge S_+'$. The one-point set
$\emptyset_+$ is both the initial and the terminal object of
$\Sets_+$, and for any $S \in \Sets_+$, we have
\begin{equation}\label{o.0}
S \wedge \emptyset_+ = \emptyset_+.
\end{equation}

\subsection{Simplicial objects.}\label{simpl.subs}

As usual, we denote by $\Delta$ the category of non-empty finite
totally ordered sets, -- or in other words, finite ordinals, -- with
$[n] \in \Delta$ being the set with $n$ elements, $n \geq 1$. A {\em
  simplicial object} in a category $\C$ is a functor $\Delta^o \to
\C$ from the opposite category $\Delta^o$. More generally, an {\em
  $n$-simplicial object} is a functor $(\Delta^o)^n \to \C$ from the
$n$-fold self-product $(\Delta^o)^n = \Delta^o \times \dots \times
\Delta^o$; we extend it to $n=0$ by setting $(\Delta^o)^0 = \ppt$, so
that a $0$-simplicial object is just an object in $\C$. We denote
the category of $n$-simplicial objects in $\C$ by
$(\Delta^o)^n\C$. For any object $c \in \C$ and any $n$-simplicial
object $c' \in (\Delta^o)^n\C$, we denote by $\C(c,c') \in
(\Delta^o)^n\Sets$ the $n$-simplicial set given by
\begin{equation}\label{C.cc}
\C(c,c')([m]) = \C(c,c'([n])), \qquad [m] \in (\Delta^o)^n.
\end{equation}
For any abelian category $\E$, we have the Dold-Kan equivalence
\begin{equation}\label{dk.eq}
\Delta^o\E \cong C_{\geq 0}(\E)
\end{equation}
between $\Delta^o\E$ and the category $C_{\geq 0}(\E)$ of chain
complexes in $\E$ concentrated in non-negative homological
degrees. The equivalence sends a simplicial object $E \in
\Delta^o\E$ to its normalized chain complex $C_\idot(E)$. For any
ring $R$ and any $E \in \Fun(\Delta^o,R)$, the complex $C_\idot(E)$
it quasiisomorphic to the homology complex $C_\idot(\Delta^o,E)$.

For any simplicial set $S \in \Delta^o\Sets$ and any abelian group
$M$, we have the standard normalized chain complex
\begin{equation}\label{c.simpl}
C_\idot(X,M)
\end{equation}
that computes the homology of $X$ with coefficients in $M$. If we
consider the functor $M[X] \in \Fun(\Delta^o,\Z)$ given by
\begin{equation}\label{m.x}
M[X]([n]) = M[X([n])],
\end{equation}
then the complex \eqref{c.simpl} corresponds to $M[X]$ under the
Dold-Kan equivalence \eqref{dk.eq}. In particular, it computes the
homology of the category $\Delta^o$ with coefficients in $M[X]$, so
that we have a natural quasiisomorphism
$$
C_\idot(X,M) \cong C_\idot(\Delta^o,M[X]).
$$
For an $n$-simplicial set $X$, there are two ways to define
homology. Firstly, we can consider the restriction $\delta^*X \in
\Delta^o\Sets$ with respect to the diagonal embedding
$\delta:\Delta^o \to (\Delta^o)^n$, and set
$$
C_\idot(X,M) = C_\idot(\delta^*X,M).
$$
Secondly, we can take
$$
C_\idot(X,M) = C_\idot((\Delta^o)^n,M[X]),
$$
where $M[X]$ has the same meaning as in \eqref{m.x}. The two
definitions are the same up to a quasiisomorphim: for any ring $R$
and any $E \in \Fun((\Delta^o)^n,R)$, the adjunction map gives a
canonical quasiisomorphism
\begin{equation}\label{shuffle.eq}
C_\idot(\Delta^o,\delta^*E) \to C_\idot((\Delta^o)^n,E).
\end{equation}
In particular, for any two simplicial sets $X_1,X_2 \in
\Delta^o\Sets$, their {\em box product} $X_1 \boxtimes X_2$ is a
$2$-simplicial set given by
$$
X_1 \boxtimes X_2 ([n_1] \times [n_2]) = X_1([n_1]) \times
X_2([n_2]),
$$
the restriction $\delta^*(X_1 \boxtimes X_2)$ is the pointwise
product $X_1 \times X_2$, and \eqref{shuffle.eq} gives the K\"unneth
quasiisomorphism
$$
C_\idot(X_1 \times X_2,\Z) \cong C_\idot(X_1 \boxtimes X_2,\Z) =
C_1(X_1,\Z) \otimes C_2(X_2,\Z).
$$
For pointed simplicial set $X \in \Delta^o\Sets_+$, the {\em reduced
  chain complex} is given by
$$
\overline{C}_\idot(X,\Z) = C_\idot(Z,\Z)/\Z \cdot o,
$$
there $o \in X([1])$ is the distiguished element. For reduced chain
complexes, the K\"unneth formula reads as
$$
\overline{C}_\idot(X_1 \wedge X_2,\Z) \cong
\overline{C}_\idot(X_1,\Z) \otimes \overline{C}_\idot(X_2,\Z),
$$
where $X_1 \wedge X_2$ is the pointwise smash product of pointed
simplicial sets $X_1,X_2 \in \Delta^o\Sets_+$.

\subsection{Contractible simplicial sets.}\label{contr.subs}

Denote by $\Delta_+ \subset \Delta$ the subcategory of non-empty
finite ordinals and maps between them that send the terminal element
to the terminal element. The embedding $\Delta_+ \subset \Delta$ has
a left-adjoint functor $\kappa:\Delta \to \Delta_+$ that adds a new
terminal element to an ordinal. We will denote by $[n]_+ \in
\Delta_+$, $n \geq 0$ the ordinal with $n+1$ elements, so that
$\kappa([n]) = [n]_+$.

We say that a $0$-simplicial set $S \in \Sets$ is {\em contractible}
if $S = \ppt$, and we introduce the following inductive definition.

\begin{defn}\label{contr.def}
  For any $n \geq 1$, an $n$-simplicial set $X \in (\Delta^o)^n\Sets$
  is {\em contractible} if we have
$$
X \cong (\kappa \times \id^{n-1})^*\wt{X}
$$
for some functor $\wt{X}:\Delta_+^o \times (\Delta^o)^{n-1} \to
\Sets$ called an {\em extension} of $X$, and the restriction of
$\wt{X}$ to $[0]_+ \times (\Delta^o)^{n-1} \subset \Delta^o_+ \times
(\Delta^o)^{n-1}$ is a contractible $(n-1)$-sim\-pli\-ci\-al set.
\end{defn}

We say that a pointed $n$-simplicial set is contractible if it
becomes contractible after forgetting the distiguished point, and we
note that the product of contractible $n$-simplicial sets and the
smash-product of contractible pointed $n$-simplicial sets is
obviously contractible. Moreover, by \eqref{o.0}, the smash product
$S \wedge X$ of a contractible pointed $n$-simplicial $X$ and a
pointed set $S$ is automatically contractible.

\begin{lemma}\label{contr.le}
For any contractible $n$-simplicial set $X$, we have a natural
isomorphism
$$
C_\idot(X,\Z) \cong \Z.
$$
\end{lemma}

\proof{} Applying \eqref{shuffle.eq} inductively, we see that it
suffices to prove that for any $E \in \Fun(\Delta^o_+,\Z)$, the
natural map
$$
C_\idot(\Delta^o,\kappa^*E) \to E([0]_+)
$$
is a quasiisomorphism. Since $\kappa:\Delta^o \to \Delta^o_+$ is
right-adjoint to the embedding $\Delta^o_+ \subset \Delta^o$, we
have
$$
C_\idot(\Delta^o,\kappa^*E) \cong C_\idot(\Delta^o_+,E),
$$
and since $[0]_+ \in \Delta^o_+$ is the terminal object, the claim
is obvious.
\endproof

We will need the following standard examples of contractible
$n$-simplicial sets.

\begin{exa}
  For any $n \geq 0$, the {\em elementary $n$-simplex} $I_n \in
  \Delta^o\Sets$ given by
$$
I_n([m]) = \Delta([m],[n+1]), \qquad [m] \in \Delta^o
$$
is contractible ($\wt{I_n} \in \Delta^o_+\Sets$ is given by
$I_n([m]_+) = \Delta_+([m]_+,[n]_+)$).
\end{exa}

\begin{exa}\label{E.exa}
  Let $\rho:\Delta \to \Sets$ be the forgetful functor taking an
  ordinal to itself considered simply as a set, and extend it to the
  forgetful functor $\rho_+:\Delta_+ \to \Sets_+$ by taking the
  terminal element as the distiguished point. For any set $S$, let
  $ES \in \Delta^o\Sets$ be the simplicial set given by
$$
ES([n]) = S^n = \Sets(\rho([n]),S).
$$
Then as soon as $S$ is not empty, $ES$ is contractible -- to obtain an
extension $\wt{ES} \in \Delta^o_+\Sets$, let
$$
\wt{ES}([n]_+) = \Sets_+(\rho_+([n]_+),S),
$$
where we choose an arbitrary element in $S$ and treat it as a
pointed set.
\end{exa}

\begin{exa}\label{cone.exa}
  For any simplicial set $X \in \Delta^o\Sets$, let $C(X)$ be the
  $2$-simplicial pointed set obtained by pushout square
$$
\begin{CD}
X \boxtimes \{s,t\} @>>> \{s,t\}\\
@VVV @VVV\\
X \boxtimes I @>>> C(X),
\end{CD}
$$
where $I=I_1 \in \Delta^o\Sets$ is the simplicial interval, $\{s,t\}
\subset I$ is the constant simplicial set formed by its two ends,
and we choose $s \in \{s,t\} \subset C(X)$ as the
distiguished point.  Then if $X$ is contractible, $C(X)$ is also
contractible. Indeed, an extention $\wt{C(X)}:\Delta^o_+ \times
\Delta^o \to \Sets$ is given by the same square with $X$ replaced by
its extension $\wt{X}$, and since $\wt{X}([0]_+)$ must be the
one-element set $\ppt$, the restriction of $\wt{C(X)}$ to $[0]_+
\times \Delta^o$ is the interval $I$.
\end{exa}

\begin{remark}
Example~\ref{cone.exa} is the reason we need to bother with
$n$-simplicial sets -- I do not know whether one can define a
contractible $C(X) \in  \Delta^o\Sets_+$ with the same properties.
\end{remark}

\subsection{Triangulated categories.}\label{tria.subs}

Finally, we need some results on triangulated categories. It is
worthwhile to mention that as in \cite{mackey}, our notion of a
triangulated category is the original one of Verdier, with derived
categories $\D(\C,R)$ giving examples. We also assume known the
notion of a {\em $t$-structure} on a triangulated category $\D$
introduced in \cite{BBD}, and the fact that such a structure is
completely defined by the subcategories $\D^{\leq i} \subset \D$,
$i$ an integer. The {\em standard $t$-structure} on $\D = \D(\C,R)$
is obtained by taking as $\D^{\leq i}(\C,R)$ the full subcategory
spanned complexes concentrated in cohomological degrees $\leq
i$. For any $t$-structure, we denote by $\D^- \subset \D$ the union
of all the full subcategories $\D^{\leq i} \subset \D$.

For an inverse system $E_i$ of objects in a triangulated category
$\D$ with infinite products, the {\em telescope} $\Tel(E_i)$ is
defined by the distiguished triangle
\begin{equation}\label{Tel}
\begin{CD}
\Tel(E_i) @>>> \prod_iE_i @>{\id - t}>> \prod_iE_i @>>>,
\end{CD}
\end{equation}
where $t$ is the product of the transition maps $E_{i+1} \to
E_i$. If we throw away a finite number of terms in the inverse
system, the telescope does not change, and if the inverse sysem is a
constant one with value $E$, then its telescope is also isomorphic
to $E$. If $\D = \D(\C,R)$, then the inverse system $E_i$ comes from
an inverse system of complexes in $\Fun(\C,R)$, and we have
\begin{equation}\label{tel.eq}
\Tel(E_i) \cong \dlim_{\overset{i}{\gets}}E_i,
\end{equation}
where $\lim^\hdot$ in the right-hand side is the derived functor of
the inverse limit.

Analogously, for a direct system $E_i$ of objects in $\C$, its
telescope $\tel(E_i)$ is defined by the distinguished triangle
\begin{equation}\label{tel}
\begin{CD}
\bigoplus_iE_i @>{\id - t}>>\bigoplus_i E_i @>>> \tel(E_i) @>>>,
\end{CD}
\end{equation}
where $t$ is the sum of the transition maps $E_i \to E_{i+1}$. It
enjoys the properties similar to $\Tel$. Moreover, if $\D = \D(R)$
is the derived category of modules over a ring $R$ -- or more
generally, the derived category of an abelian category satisfying
$AB5$ -- then the homology objects $\hh_\idot(\tel(E_i))$ are given by
$$
\hh_\idot(\tel(E_i)) = \lim_{\overset{i}{\to}}\hh_\idot(E_i),
$$
where $\hh_\idot(E_i)$ are the homology objects of $E_i$.

We will also need one more technical result. Assume given two
triangulated categories $\D_1$, $\D_2$, and a pair of triangulated
functors $\phi:\D_2 \to \D_1$, $\rho:\D_1 \to \D_2$ such that $\rho$
is right-adjoint to $\phi$. Moreover, assume given full triangulated
subcategories $\D_1' \subset \D_1$, $\D_2' \subset \D_2$ such that
the quotients $\D_1/\D_1'$, $\D_2/\D_2'$ are well-defined, and
denote by $q_1:\D_1 \to \D_1/\D_1'$, $q_2:\D_1 \to \D_1/\D_1'$ the
quotient functors. Finally, assume that $\D_1' \subset \D_1$ is
left-admissible -- that is, $q_1$ has a right-adjoint functor
$r:\D_1/\D_1' \to \D_1$.

\begin{lemma}\label{tria.le}
In the assumptions above, assume that $\phi:\D_2 \to \D_1$ sends
$\D_2' \subset \D_2$ into $\D_1' \subset \D_1$, thus induces a
functor
$$
\overline{\phi}:\D_2/\D_2' \to \D_1/\D_1'.
$$
Then the functor
$$
\overline{\rho} = q_2 \circ \rho \circ r:\D_1/\D_1' \to \D_2/\D_2'
$$
is right-adjoint to the functor $\overline{\phi}$.
\end{lemma}

\proof{} We note that if the projection $q_2$ has a left-adjoint
functor $l:\D_2/\D_2' \to \D_2$, the statement is obvious:
$\overline{\rho}$ is right-adjoint to $q_1 \circ \overline{\phi}
\circ l \cong \phi \circ q_2 \circ l$, and since $q_2$ is the
localization functor, $l$ must be fully faithful, so that $q_2 \circ
l \cong \Id$.

In the general case, note that for an arbitrary object $E \in \D_2$,
we have a natural map
\begin{equation}\label{eq.1}
\begin{CD}
  q_2(E) @>{q_2(a(\phi))}>> (q_2 \circ \rho \circ \phi)(E) @>{(q_2
    \circ \rho)(a(q_1))}>> (q_2 \circ \rho \circ r \circ q_1 \circ
  \phi)(E)
\end{CD}
\end{equation}
where $a(\phi):\Id \to \rho \circ \phi$, $a(q_2):\Id \to r \circ
q_1$ are the adjunction maps, and by definition, we have
$$
(q_2 \circ \rho \circ r \circ q_1 \circ \phi)(E) \cong (q_2 \circ
\rho \circ r \circ \overline{\phi} \circ q_2)(E) = (\overline{\rho}
\circ \overline{\phi})(q_2(E)).
$$
The map \eqref{eq.1} is functorial in $E$, and since $q_2$ is a
localization functor, it induces a map
$$
a(\overline{\phi}):\Id \to \overline{\rho} \circ \overline{\phi}.
$$
Analogously, for any $E \in \D_1$, we have a diagram
$$
\begin{CD}
  (\phi \circ \rho \circ r \circ q_1)(E) @>{a(\rho)}>> (r \circ
  q_1)(E) @<{a(q_1)}<< E,
\end{CD}
$$
where $a(\rho):\rho \circ \phi \to \Id$ is the adjunction map, and
by the definition of Verdier localization, it induces a map
$$
a(\overline{\rho}):\overline{\phi} \circ \overline{\rho} \to \Id.
$$
To prove that the maps $a(\overline{\phi})$ and $a(\overline{\rho})$
define an adjunction between $\overline{\phi}$ and
$\overline{\rho}$, it remains to check that the composition morphisms
$$
\begin{CD}
\overline{\rho} @>{a(\overline{\rho})}>> \overline{\rho} \circ
\overline{\phi} \circ \overline{\rho}
@>{\overline{\rho}(a(\overline{\phi}))}>> \overline{\rho},\qquad
\overline{\phi} @>{\overline{\phi}(a(\overline{\rho}))}>>
\overline{\phi} \circ \overline{\rho} \circ \overline{\phi}
@>{a(\overline{\phi})}>> \overline{\phi},
\end{CD}
$$
are equal to the identity maps. This is a straightforward exercise
that we leave to the reader.
\endproof

\section{Recollection on Mackey functors.}\label{mack.class.sec}

We will also need some facts from the classical theory of Mackey
functors; this is reasonably well-known, but we include a brief
sketch for the convenience of the reader. We follow the exposition
in \cite[Section 2]{mackey}.

\subsection{Definitions.}\label{ma.def.subs}

For any group $G$, let $\B^G$ be the category with the same objects
as $\Gamma_G$, the category of finite $G$-sets, and with $\Hom$-sets
given by
\begin{equation}\label{B.G}
\B^G(S_1,S_2) = \Z[\Iso(\Gamma_G/S_1 \times S_2)]/([S \copr S'] -
    [S] - [S']).
\end{equation}
That is, $\B^G(S_1,S_2)$ is the free abelian group spanned by
isomorphism classes of diagrams
\begin{equation}\label{domik}
\begin{CD}
S_1 @<<< S @>>> S_2
\end{CD}
\end{equation}
in $\Gamma_G$, modulo the relations $[S \copr S'] = [S] + [S']$,
where $\copr$ stands for disjoint union. Composition of diagrams is
defined by taking fibered products. One checks easily that
$\B^G$ is an additive category, with direct sum given by disjoint
union of $G$-sets. 

\begin{defn}
A {\em $G$-Mackey functor} with values in some ring $R$ is an
additive functor from $\B^G$ to the category of $R$-modules. The
category of $R$-valued $G$-Mackey functors is denoted by $\M(G,R)$.
\end{defn}

\begin{exa}\label{burn.exa}
For any $S \in \Gamma_G$, let $\A(S) = \B^G(\ppt,S)$, where $\ppt$
is the one-point set with the trivial $G$-action. Then $\A$ is a
$\Z$-valued Mackey functor known as {\em Burnside Mackey
  functor}. Explicitly, we have
$$
\A(S) = \Z[\Iso(\Gamma_G/S)]/([S \copr S'] - [S] - [S']).
$$
\end{exa}

We further note that if $S = \ppt$, then $\A(\ppt) =
\B^G(\ppt,\ppt)$ is a ring, known as the {\em Burnside ring} of the
group $G$ and denoted $\A^G$. As a group, $\A^G$ is freely generated
by the set of isomorphisms classes of $G$-orbits, that is, $G$-sets
of the form $G/H$, $H \subset G$ a cofinite subgroup. The product is
induced by the cartesian product of orbits. In particular, it is
obviously commutative. Moreover, for a $G$-orbit $S = G/H$, we have
a natural equivalence
\begin{equation}\label{orb.fib}
\Gamma_G/S \cong \Gamma_H,
\end{equation}
so that $\A(G/H) = \A^H$, the Burnside ring of the group $H$. For
any ring $R$, the Burnside ring $\A^G$ maps naturally to the center
of the category $\M(G,R)$ -- that is, for any $M \in \M(G,R)$ and
any $a \in \A^G$, we have a natural map $a:M \to M$. Explicitly, if
$a = [S_0] \in \A^G$ is the class of a $G$-set $S_0 \in \Gamma_G$,
then for any $S \in \Gamma_G$, the map $a(S):M(S) \to M(S)$
corresponds to the diagram
\begin{equation}\label{burn.act.eq}
\begin{CD}
S @<{f}<< S \times S_0 @>{f}>> S,
\end{CD}
\end{equation}
where $f:S \times S_0 \to S$ is the natural projection.

A useful alternative description of Mackey functors is due to
Lindner \cite{lind}. Let $Q\Gamma_G$ be the category with the same
objects as $\Gamma_G$, and with morphisms given by isomorphism
classes of diagrams \eqref{domik}. Then by definition, both
$\Gamma_G$ and the opposite category $\Gamma_G^o$ are subcategories
in $Q\Gamma_G$: objects are the same, and a diagram \eqref{domik}
defines a map in $\Gamma_G \subset Q\Gamma_G$ resp. $\Gamma_G^o
\subset Q\Gamma_G$ if and only if the map $S \to S_1$ resp. $S \to
S_2$ is an isomorphism.

\begin{defn}\label{add.def}
  A functor $M$ from $\Gamma_G^o$ to an additive category is {\em
    additive} if the natural map
\begin{equation}\label{add.def.eq}
M(S_1 \copr S_2) \to M(S_1) \oplus M(S_2)
\end{equation}
is an isomorphism for any $S_1,S_2 \in \Gamma_G$. A functor $M$ from
$Q\Gamma_G$ is additive if so is its restriction to $\Gamma_G^o
\subset Q\Gamma_G$.
\end{defn}

Then for any ring $R$, the full subcategory
$$
\Fun_{add}(Q\Gamma_G,R) \subset \Fun(Q\Gamma_G,R)
$$
spanned by additive functors is naturally identified with $\M(G,R)$,
and the embedding has a left-adjoint additivization functor
\begin{equation}\label{Add}
\Add:\Fun(Q\Gamma_G,R) \to \Fun_{add}(Q\Gamma_G,R) \cong \M(G,R).
\end{equation}
The maps \eqref{add.def.eq} for all $S_1,S_2 \in \Gamma_G$ can be
bundled together into a single map of functors
$$
s^* \to p_1^* \oplus p_2^*,
$$
where we denote by $p_1,p_2:Q\Gamma_G \times Q\Gamma_G \to
Q\Gamma_G$ the natural projections, and we denote by $s:Q\Gamma_G
\times Q\Gamma_G \to Q\Gamma_G$ the disjoint union functor. Then $M
\in \Fun(Q\Gamma_G,R)$ is additive if and only if the map $s^*M \to
p_1^*M \oplus p_2^*M$ is an isomorphism. We also note that $s$ is
right and left-adjoint to the diagonal embedding $\delta:Q\Gamma_G
\to Q\Gamma_G \times Q\Gamma_G$, and the projections $p_1$, $p_2$
are right and left-adjoint to the embeddings $i_1,i_2:Q\Gamma_G \to
Q\Gamma_G \times Q\Gamma_G$ sending $S$ to $S \times \emptyset$
resp. $\emptyset \times S$. Thus $s^* \cong \delta_!$, $p_1^* \cong
i_{1!}$, $p_2^* \cong i_{2!}$, and saying that \eqref{add.def.eq} is
an isomorphism is equivalent to saying that the natural map
\begin{equation}\label{add.def.eq.2}
\delta_!M \to i_{1!}M \oplus i_{2!}M
\end{equation}
is an isomorphism

It is clear from Lindner's description that an $R$-valued Mackey
functor $M \in \M(G,R)$ is completely defined by the following data:
\begin{enumerate}
\item an $R$-module $M(S)$ for any finite $G$-set $S \in \Gamma_G$,
  and
\item two maps
\begin{equation}\label{f.st.st}
f_*:M(S) \to M(S'), \qquad f^*:M(S') \to M(S)
\end{equation}
for any map $f:S \to S'$ in $\Gamma_G$.
\end{enumerate}
By additivity, it suffices to specify $M(S)$ and $f_*$, $f^*$ for
$G$-orbits $S = G/H$. Traditionally, $M(G/H)$ is denoted by $M^H$.
The maps $f_*$ and $f^*$ should satisfy some compatibility
conditions encoded in the structure of the category
$Q\Gamma_G$. Explicitly, for any two maps $f':G/H' \to G/H$,
$f'':G/H'' \to G/H$ induced by embeddings $H',H'' \subset H$, the
fibered product $(G/H') \times_{G/H} (G/H'')$ decomposes into a
disjoint union
\begin{equation}\label{orb.spl}
(G/H') \times_{G/H} (G/H'') = \coprod_{s \in S}G/H_s
\end{equation}
of $G$-orbits indexed by the finite set $S = H' \setminus\! H / H''$,
and we must have
\begin{equation}\label{double}
f^{'*} \circ f''_* = \sum_{s \in S}\wt{f}''_{s*} \circ
\wt{f}^{'*}_s,
\end{equation}
where $f'_s:G/H_s \to G/H''$, $f''_s:G/H_s \to G/H'$ are projections of
the component $G/H_s$ of the decomposition \eqref{orb.spl}. This is
known as the {\em double coset formula}.

\subsection{Fixed points.}\label{ma.fp.subs}

Assume now given a subgroup $H \subset G$, let $N_H \subset G$ be
its normalizer, and let $W = N_H/H$. Then for any finite $G$-set
$S \in \Gamma_G$, the subset $S^H \subset S$ of $H$-fixed points is
naturally a finite $W$-set. Sending $S \in \Gamma_G$ to $S^H$
gives a functor
\begin{equation}\label{phi.H}
\phi^H:\Gamma_G \to \Gamma_W.
\end{equation}
This functor is left-exact, thus extends to a functor
$Q(\phi^H):Q\Gamma_G \to Q\Gamma_W$. Since $Q(\phi^H)$ obviously
commutes with the functors $\delta$, $i_1$, $i_2$ of
\eqref{add.def.eq.2}, the corresponding left Kan extension functor
$$
Q(\phi^H)_!:\Fun(Q\Gamma_G,R) \to \Fun(Q\Gamma_W,R)
$$
preserves additivity in the sense of Definition~\ref{add.def}. By
definition, the {\em geometric fixed points functor} 
$$
\Phi^H:\M(G,R) = \Fun_{add}(Q\Gamma_G,R) \to \M(W,R) =
\Fun_{add}(Q\Gamma_W,R)
$$
is induced by the functor $Q(\phi^H)_!$.

On the other hand, one can simply restrict a $G$-action on a set $S$
to an action of the subgroup $H \subset G$; this gives a left-exact
functor $\rho^H:\Gamma_G \to \Gamma_H$ and its extension
$Q(\rho^H):Q\Gamma_G \to \Gamma_H$. The corresponding left Kan
extension functor $Q(\rho^H)_!$ also preserves additivity and
induces the {\em categorical fixed points functor}
$$
\Psi^H = Q(\rho^H)_!:\M(G,R) \to \M(H,R).
$$
We note that for two subgroups $H' \subset H \subset G$, we have an
obvious isomorphism
\begin{equation}\label{Phi.Psi.iso}
\phi^{H'} \circ \rho^H \cong \rho^{W_H} \circ \phi^{H'},
\end{equation}
where we denote $W' = N_{H'}/H'$ and $W_H = (N_{H'} \cap H)/H'
\subset W'$, and it induces an isomorphism
\begin{equation}\label{phi.psi.iso}
\Phi^{H'} \circ \Psi^H \cong \Psi^{W^H_{H'}} \circ \Phi^{H'}.
\end{equation}
We also note that for any subgroup $H \subset G$, the centralizer
$Z_H \subset G$ of the group $H$ acts on the functor $\rho^H$, thus
on $\Psi^H$, so that it can be promoted to a functor
\begin{equation}\label{psi.wt}
\wt{\Psi}^H:\M(G,R) \to \M(H,R[Z_H]),
\end{equation}
where $R[Z_H]$ is the group algebra of the centralizer. If $H
\subset G$ is cofinite, then the functor $\rho^H$ has a left-adjoint
functor
\begin{equation}\label{gamma.h}
\gamma_H:\Gamma_H \to \Gamma_G
\end{equation}
given by the identification \eqref{orb.fib} composed with the
forgetful functor (explicitly, we have $\gamma_H(S) = (G \times
S)/H$). The functor $\gamma_H$ is also left-exact, and the
corresponding functor $Q(\gamma_H):Q\Gamma_H \to Q\Gamma_G$ is adjoint
to $Q(\rho^H)$ both on the left and on the right. Thus we have
$$
\Psi^H = Q(\rho^H)_! \cong Q(\gamma_H)^*,
$$
and this allows to compute $\Psi^H$ rather explicitly (in
particular, it is an exact functor).

Geometric fixed points are usually much more difficult to compute
explicitly, and we will only give one result in this
direction. Denote by
\begin{equation}\label{infl.eq}
\Infl^H = Q(\phi^H)^*:\M(W,R) \to \M(G,R)
\end{equation}
the {\em inflation functor} right-adjoint to $\Phi^H$. Assume that the
subgroup $H \subset G$ is normal, so that $W = G/N$, and Let
\begin{equation}\label{lambda.N}
\lambda_N:\Gamma_{G/N} \to \Gamma_G
\end{equation}
be the fully faithful embedding sending a $G/N$-set $S$ to the same
set on which $G$ acts via the map $G \to G/N$. For every $G$-orbit
$S$, we have one of the two alternatives: either $S^N$ is empty, or
$S^N = S$, and $S$ lies in the image of the embedding
$\lambda_N$. Let
$$
Q(\lambda_N):Q\Gamma_{G/H} \to Q\Gamma_G
$$
be the functor induced by the embedding $\lambda_N$, and let
$$
\M_N(G,R) \subset \M(G,R)
$$
be the full subcategory spanned by $M \in \M(G,R)$ such that
$M(S)=0$ for any $S \in \Gamma_G$ with empty $S^N$.

\begin{lemma}\label{infl.le}
The inflation functor $\Infl^N$ is fully faithful and identifies the
category $\M(W,R)$ with the full subcategory $\M_N(G,R) \subset
\M(G,R)$, with inverse equivalence given by $Q(\lambda_N)^*$. For any
$M \in \M(G,R)$, the adjunction map
\begin{equation}\label{infl.adj}
M \to \Infl^N\Phi^NM
\end{equation}
is surjective, and for any $S \in \Gamma_{G/N}$, we have a short
exact sequence
\begin{equation}\label{phi.exp}
\begin{CD}
\displaystyle\bigoplus_{f:S' \to S}M(S') @>{\sum f_*}>> M(S) @>>>
\Phi^N(M)(S) @>>> 0,
\end{CD}
\end{equation}
where $f_*$ is as in \eqref{f.st.st}, and the sum is over all maps
$f:S' \to S$ in $\Gamma_G$ such that $S'$ has no elements fixed
under $N \subset G$.
\end{lemma}

\begin{remark}\label{f.uni}
Note that by additivity of $M \in \Fun(Q\Gamma_G,R)$, the image of
the map $\sum f_*$ in \eqref{phi.exp} is not only the sum of the
images of individual maps $f_*$, it is actually the union of these
images.
\end{remark}

\proof{} The functor $\lambda_N:\Gamma_{G/N} \to \Gamma_G$ is a
fully faithful embedding, and $\rho^N$ is by definition
right-adjoint to this embedding. Therefore for any $M \in \M(W,R)$,
we have an adjunction map
\begin{equation}\label{adj.1}
\overline{M} \to \lambda_N^*\phi^{N*}\overline{M},
\end{equation}
where $\overline{M} \in \Fun(\Gamma_W,R)$ is the restriction of the
Mackey functor $M$ to $\Gamma_W \subset Q\Gamma_W$, and analogously,
for any $M \in \M(G,R)$, we have an adjunction map
\begin{equation}\label{adj.2}
\phi^{N*}\lambda_N^*\overline{M} \to \overline{M},
\end{equation}
where $\overline{M}$ is obtained by restriction to $\Gamma_G \subset
Q\Gamma_G$. Moreover, the map \eqref{adj.1} is an isomorphism, and
it is compatible with the maps $f^*$ of \eqref{f.st.st}, so that it
gives an isomorphism
$$
M \cong Q(\lambda_N)^*Q(\phi^N)^*M = Q(\lambda_N)^*\Infl^NM.
$$
In general, the map \eqref{adj.2} does not commute with the maps
$f^*$, but it clearly does so if $M$ lies in $\M_N(G,R) \subset
\M(G,R)$, and in this case, it is an isomorphism. Since by
definition, $\Infl^N = Q(\phi^N)^*$ sends $\M(W,R)$ into
$\M_N(G,R)$, we conclude that $\Infl^N$ and $Q\lambda_N^*$ are
indeed inverse equivalences between $\M(W,R)$ and
$\M_N(G,R)$. Now, the cokernel of the map $\sum f_*$ in
\eqref{phi.exp} is clearly functorial in $S$ and trivial for $S$
with empty $S^N$, thus defines a functor $\Phi:\M(G,R) \to
\M_N(G,R)$ equipped with a map $\Id \to \Phi$. This map becomes an
isomorphism after restricting to $\M_N(G,R)$, thus gives an
adjunction between $\Phi$ and the embedding $\M_N(G,R) \subset
\M(G,R)$. This identifies $\Phi$ and $\Phi^N$.
\endproof

\subsection{Products.}\label{ma.prod.subs}

Assume now given two rings $R$, $R'$. Then one can define a natural
tensor product functor
\begin{equation}\label{ma.bo}
\M(G,R) \times \M(G,R') \to \M(G,R \otimes R').
\end{equation}
To do it, use Lindner's description. Cartesian product of $G$-sets
commutes with fibered products, thus extends to a functor
\begin{equation}\label{m.prod}
m:Q\Gamma_G \times Q\Gamma_G \to Q\Gamma_G.
\end{equation}
The product $M \circ M'$ of two Mackey functors $M \in \M(G,R)$, $M'
\in \M(G,R')$
is then defined by
\begin{equation}\label{ma.produ}
M \circ M' = \Add(m_!(M \boxtimes M')),
\end{equation}
where $M \boxtimes M'$ in the right-hand side is given by $(M
\boxtimes M')(S \times S') = M(S) \otimes M'(S')$. This is clearly
associative and commutative in the obvious sense. Moreover, for any
subgroup $H \subset G$, the functors $Q\rho^H:Q\Gamma_G \to
Q\Gamma_H$, $Q\phi^H:Q\Gamma_G \to Q\Gamma_W$ commute with $m$;
therefore the fixed points functors $\Psi^H$, $\Phi^H$ are tensor
functors.

For any three Mackey functors $M_1 \in \M(G,R_1)$, $M_2 \in
\M(G,R_2)$, $M \in \M(G,R_1 \otimes R_2)$, a map $M_1 \circ M_2 \to
M$ by adjunction corresponds to a map
\begin{equation}\label{adj.mack.prod}
\wt{\mu}:M_1 \boxtimes M_2 \to m^*M.
\end{equation}
To see such maps more explicitly, one can first restrict to the
subcategory $\Gamma_G^o \subset Q\Gamma_G$. Then the functor $m$ is
left-adjoint to the diagonal embedding $\delta:\Gamma_G^o \to
\Gamma_G^o \times \Gamma_G^o$, and again by adjunction,
\eqref{adj.mack.prod} induces a map
$$
\mu:M_1 \otimes M_2 = \delta^*(M_1 \boxtimes M_2) \to M
$$
in $\Fun(\Gamma_G^o,R_1 \otimes R_2)$. In other words, we have a map
\begin{equation}\label{mu.h}
\mu:M_1^H \otimes M_2^H \to M^H
\end{equation}
for any cofinite subgroup $H \subset G$, and for any two cofinite
subgroups $H_1,H_2 \subset G$ and a $G$-map $f:G/H_1 \to G/H_2$, we
have
\begin{equation}\label{mack.prod.1}
\mu(f^*(a_1) \otimes f^*(a_2)) =f^*\mu(a_1 \otimes a_2),
\ \ \quad
a_1 \in M_1^{H_2}, a_2 \in M_2^{H_2}.
\end{equation}

\begin{lemma}
A collection of maps \eqref{mu.h} satisfying \eqref{mack.prod.1}
define a map \eqref{adj.mack.prod} if and only if we also have
\begin{equation}\label{mack.prod.2}
\begin{aligned}
f_*(\mu(a_1 \otimes f^*(a_2)) &= \mu(f_*(a_1) \otimes a_2), \quad
a_1 \in M_1^{H_1}, a_2 \in M_2^{H_2},\\
f_*(\mu(f^*(a_1) \otimes a_2) &= \mu(a_1 \otimes f_*(a_2)), \quad
a_1 \in M_1^{H_2}, a_2 \in M_2^{H_1},
\end{aligned}
\end{equation}
for any two cofinite subgroups $H_1,H_2 \subset G$ and a map
$f:G/H_1 \to G/H_2$.
\end{lemma}

\proof{} By adjunction, the maps \eqref{mu.h} define a map
\eqref{adj.mack.prod} in the category $\Fun(\Gamma_G^o \times
\Gamma_G^o,R)$; explicitly, it is given by
\begin{equation}\label{mu.eq.1}
\wt{\mu}(b_1 \otimes b_2) = \mu(p^*_1(b_1) \otimes p^*_2(b_2))
\end{equation}
for any cofinite $H_1,H_2 \subset G$ and any elements $b_1 \in
M_1(G/H_1)$, $b_2 \in M_2(G/H_2)$, where $p_1:(G/H_1) \times (G/H_2)
\to G/H_1$, $p_2:(G/H_1) \times (G/H_2) \to G/H_2$ are the natural
projections. This map $\wt{\mu}$ is a map in $\Fun(Q\Gamma_G \times
Q\Gamma_G,R)$ if and only if for any cofinite $H_1' \subset G$, a
map $g:G/H_1' \to G_H$, and an element $b_1 \in M_1(G/H_1')$, we
have
\begin{equation}\label{mu.eq.2}
\wt{\mu}(g_*b_1 \otimes b_2) = (g \times \id)_*(\wt{\mu}(b_1 \otimes
b_2)),
\end{equation}
and similarly for the product in the opposite order. Combining
\eqref{mu.eq.1} and \eqref{mu.eq.2} gives
\begin{equation}\label{mu.eq.3}
\mu(p_1^*g_*(b_1) \otimes p_2^*(b_2)) = (g \times
\id)_*\mu(p_1^*(b_1) \otimes p_2^{'*}(b_2)),
\end{equation}
where $p_2':(G/H_1') \times (G/H_2) \to G/H_2$ is the natural
projection. Since $p_1^* \circ g_* = (g \times \id)_* \circ p_1^*$
and $p_2' = p_2 \circ (g \times \id)$, \eqref{mu.eq.3} is exactly
the first equation in \eqref{mack.prod.2} with $f = g \times \id$,
$a_1 = p_1^*b_1$, $a_2 = p_2^*b_2$. Thus \eqref{mu.eq.1} implies
that $\wt{\mu}$ is a map in $\Fun(Q\Gamma_G \times
Q\Gamma_G,R)$. Conversely, $\mu$ is expressed in terms of $\wt{\mu}$
by
$$
\mu(a_1 \otimes a_2) = \delta^*\wt{\mu}(a_1 \otimes a_2), \qquad
a_1,a_2 \in M^H,
$$
where $\delta:G/H \to (G/H \times G/H)$ is the diagonal embedding,
and since $\delta^* \circ (g \times \id)_* = g_* \circ \delta^*$,
the equation \eqref{mu.eq.2} together with the opposite equation
imply \eqref{adj.mack.prod}.
\endproof

In particular, if our base ring $R$ is commutative, we have a
natural product map $R \otimes R \to R$, so that the product
\eqref{ma.bo} induces a symmetric tensor product on the category
$\M(G,R)$.  A {\em Green functor} is a ring object in the category
$\M(G,R)$. Explicitly, a Green functor structure on a Mackey functor
$M$ is given by a ring structure on every $M^H$, $H \subset G$
cofinite; the maps $f^*$ are ring maps, and the maps $f_*$ satisfy
the projection formula \eqref{mack.prod.2}. For example, the
Burnside Mackey functor $\A$ of Example~\ref{burn.exa} is a Green
functor (in fact, $\A$ is the unit object for the product on
$\M(G,\Z)$, and more generally, the unit object for the product
\eqref{ma.bo}).

\newpage

\section{Mackey profunctors.}\label{mackey.sec}

\subsection{Definitions.}\label{pro.def.subs}

In the classical theory of Mackey functors, one always assumes that
the group $G$ is finite (or compact, in the topological versions of
the theory). We do not do so since there is no formal need: it is
only the $G$-sets that have to be finite (otherwise $\B^G(-,-)$ of
\eqref{B.G} would be identically $0$). However, for a infinite group
$G$, there is the following useful alternative version of the
category $\M(G,R)$.

\begin{defn}\label{adm.G.set}
For any group $G$, a $G$-set $S$ is {\em admissible} if
\begin{enumerate}
\item for any element $s \in S$, its stabilizer $G_s \subset G$ us
  cofinite, and
\item for any cofinite subgroup $H \subset S$, its fixed point set
  $S^H \subset S$ is finite.
\end{enumerate}
\end{defn}

Denote the category of admissible $G$-sets by $\wGamma_G$, and note
that it is a small category (with the number of isomorphism classes
of objects depending on the cardinality of $G$). A subset of an
admissible $G$-set and the product of two admissible $G$-set is
admissible, so that we can form the quotient category $Q\wGamma_G$
and consider the functor category $\Fun(Q\wGamma_G,R)$ for any ring
$R$. As in the case of finite sets, we have natural embeddings
$\wGamma_G \subset Q\wGamma_G$, $\wGamma_G^o \subset Q\wGamma_G$.

\begin{defn}\label{add.pro}
  A functor $M$ from $\wGamma_G^o$ to an additive category with
  arbitrary products is {\em additive} if for any admissible $G$-set
  $S \in \wGamma_G$, with the quotient map $q:S \to S/G$, the
  natural map
\begin{equation}\label{add.pro.eq}
M(S) \to \prod_{s \in S/G}M(q^{-1}(s))
\end{equation}
induced by the decomposition $S = \coprod_{s \in S/G}q^{-1}(s)$ is an
isomorphism. An {\em $R$-valued $G$-Mackey profunctor} $M$ is a
functor from $Q\wGamma_G$ to the category of $R$-modules whose
restriction to $\wGamma_G^o$ is additive.
\end{defn}

We will denote the category of $R$-valued $G$-Mackey profunctors by
$$
\wM(G,R) \subset \Fun(Q\wGamma_G,R).
$$
We note right away that since the abelian category of $R$-modules
satisfies $AB4^*$, additivity is preserved under kernels and
cokernels; thus $\wM(G,R)$ is also abelian and satifies $AB4^*$.

By definition, a finite $G$-set is admissible, and additivity in the
sense of Definition~\ref{add.pro} obviously implies additivity in
the sense of Definition~\ref{add.def}, so that we have an embedding
$\Gamma_G \subset \wGamma_G$ and a restriction functor
\begin{equation}\label{M.wM.eq}
\wM(G,R) \to \M(G,R).
\end{equation}
If $G$ is finite, the restriction functor is an equivalence. For
some infinite groups -- for example, for the group $G = \Z_p$ of
$p$-adic integers -- the restriction functor, while not an
equivalence, is at least fully faithful. In general, it is not even
that.

\begin{exa}\label{wB.G.exa}
Let $\wBG$ be an additive category whose objects are admissible
$G$-sets, and whose $\Hom$-sets are given by \eqref{B.G} with
$\Gamma_G$ replaced by $\wGamma_G$. Then for any $S \in \wGamma_G$,
the representable functor $\wBG(S,-)$ is a functor from
$Q\wGamma_G$ to abelian groups, and it is obviously additive in the
sense of Definition~\ref{add.pro}. Thus it gives a $\Z$-valued
$G$-Mackey profunctor, and we obtain an embedding
$$
\wBG \subset \wM(G,\Z).
$$
This embedding is fully faithful.
\end{exa}

In particular, Example~\ref{burn.exa} generalizes to Mackey
profunctors, so that we have the {\em completed Burnside ring}
$\wAG = \wBG(\ppt,\ppt)$ and the corresponding $G$-Mackey
profunctor $\wA \in \wM(G,\Z)$. As in the case of usual Mackey
functor, $\wAG$ maps into the center of the category $\wM(G,R)$ for
any ring $R$. In fact, the full subcategory in $\wBG$ spanned by
finite $G$-sets can be alternatively described as follows. For any
cofinite subgroup $H \subset G$, classes of $G$-sets $S \in
\Gamma_G$ with empty $S^H$ form a two-side ideal $\B^G_H \subset
\B_G$. Then for any $S_1,S_2 \in \Gamma_G$, we have
$$
\wBG(S_1,S_2) =
\lim_{\overset{H}{\gets}}\B^G(S_1,S_2)/\B^G_H(S_1,S_2),
$$
where the limit is taken over all cofinite subgroups $H \subset
G$. Explicitly, we have
\begin{equation}\label{wB.G}
\wBG(S_1,S_2) = \prod_{H \subset G}\Z[(S_1 \times S_2)^H/W_H],
\end{equation}
where the product is taken over all conjugacy classes of cofinite
subgroups $H \subset G$, and for any such $H$, $W_H = N_H/H$ is the
quotient of the normalizer $N_H \subset G$ by the subgroup $H
\subset N_H \subset $. In particular, the completed Burnside ring is
indeed a completion of the usual Burnside ring -- we have
\begin{equation}\label{wA.G.eq}
\wAG = \prod_{H \subset G}\Z,
\end{equation}
where the product is over all conjugacy classes of cofinite
subgroups $H \subset G$. For the usual Burnside ring $\A^G$, the
formula is the same, but the product is replaced with the sum.

\subsection{Normal systems.}\label{norm.subs}

For any cofinite subgroup $H \subset G$ and any admissible $G$-set
$S \in \wGamma_G$, the fixed point set $S^H$ is by definition
finite, so that as in the usual Mackey functors case, we have a
natural left-exact functor $\phi^H:\wGamma_G \to \Gamma_W$ and
the corresponding adjoint pair of functors
\begin{equation}\label{phi.infl.eq}
\begin{aligned}
\Phi^H &= Q\phi^H_!:\Fun(\wGamma_G,R) \to \Fun(\Gamma_W,R),\\
\Infl^H &= Q\phi^{H*}:\Fun(\Gamma_W,R) \to \Fun(\wGamma_G,R).
\end{aligned}
\end{equation}
For the same reasons as in the usual case, these functors
preserve additivity, thus induce functors between categories of
$G$-Mackey profunctors. For any $M \in \wM(G,R)$ and any pair $N
\subset N' \subset G$ of cofinite normal subgroup, we have a
natural isomorphism
\begin{equation}\label{m_n}
\Phi^{N'}M \cong \Phi^{N'/N}\Phi^NM.
\end{equation}
Thus every Mackey profunctor $M \in \wM(G,R)$ induces a collection
of Mackey functors $\Phi^NM$ related by the isomorphisms
\eqref{m_n}.

To axiomatize the situation, denote by $\Nn(G)$ the partially
ordered set of cofinite normal subgroups $N \subset G$, ordered by
inclusion, treat it as a small category in the usual way, and let
$$
\hGamma_G \subset \Gamma_G \times \Nn(G)^o
$$
be the full subcategory spanned by pairs $\langle S,N \rangle$ such
that $N$ acts trivially on $S$. Thus explicitly, morphisms from
$\langle S,N \rangle$ to $\langle S',N' \rangle$ are pairs of a
$G$-equivariant map $f:S \to S'$ and an inclusion $N' \subset N$.
We have a natural forgetful functor
\begin{equation}\label{to.N.eq}
\nu:\hGamma_G \to \Nn(G)^o, \qquad \langle S,N \rangle \mapsto N.
\end{equation}
This functor is a fibration, with fiber over $N \in \Nn(G)$
equivalent to $\Gamma_W$, $W = G/N$, and with transition functors
$\phi^{N'/N}$, $N \subset N' \subset G$.

Now let $Q\hGamma_G$ be the category with the same objects as
$\hGamma_G$, with morphisms from $\langle S,N \rangle$ to $\langle
S',N' \rangle$ given by isomorphism classes of diagrams
$$
\begin{CD}
\langle S,N \rangle @<<< \langle \wt{S},N' \rangle @>>> \langle S',N'
\rangle
\end{CD}
$$
in $\hGamma_G$, and with compositions given by pullbacks (we note
that the necessary pullbacks in $\hGamma_G$ do exist). Then the
fibration \eqref{to.N.eq} defines a cofibration
$$
Q\hGamma_G \to \Nn(G), \qquad \langle S,N \rangle \mapsto N,
$$
with fibers $Q\Gamma_W$ and transition functors
$Q(\phi^{N'/N})$. Explicitly, an object $M \in \Fun(Q\hGamma_G,R)$
is given by a collection of objects $M_N \in \Fun(Q\Gamma_W,R)$, $N
\in \Nn(G)$, $W = G/N$, and transition morphisms
\begin{equation}\label{norm.trans}
\Phi^{N'/N}M_N \to M_{N'}
\end{equation}
for any $N,N' \in \Nn(G)$, $N \subset N'$.

\begin{defn}\label{ada.def}
For any group $G$, a {\em normal system} of $G$-Mackey functors is
an object $M \in \Fun(Q\hGamma_G,R)$ such that for any $N \in
\Nn(G)$, $W = G/N$, the object $M_N \in \Fun(\Gamma_W,R)$ is additive
in the sense of Definition~\ref{add.def}, and for any $N,N' \in
\Nn(G)$, $N \subset N'$, the transition morphism \eqref{norm.trans}
is an isomorphism.
\end{defn}

Normal systems obviously form an additive category with arbitrary
sums; we denote it by $\N(G,R)$. Moreover, we have a natural functor
\begin{equation}\label{phi.eq}
\phi:\Nn(G)^o \times \wGamma_G \to \hGamma_G
\end{equation}
given by $\phi(N \times S) = \langle S^N,N \rangle$, and it defines
a functor $Q(\phi):\Nn(G) \times Q\wGamma_G \to Q\hGamma_G$. If we
denote by $p:\Nn(G) \times Q\Gamma_G \to Q\Gamma_G$ the natural
projection, then for any Mackey profunctor $M \in \wM(G,R)$,
$\phi_!p^*M \in \Fun(Q\hGamma_G,R)$ is a normal system in the sense
of Definition~\ref{ada.def}, so that we have a natural functor
$$
\Phi = Q(\phi)_! \circ p^*:\wM(G,R) \to \N(G,R).
$$
with a right-adjoint functor
$$
\Infl = p_* \circ Q(\phi)^*:\N(G,R) \to \wM(G,R).
$$
Explicitly, for any $M \in \wM(G,R)$, we have $\Phi(M)_N \cong
\Phi^NM$, and the transition maps \eqref{norm.trans} are the
isomorphisms \eqref{m_n}. Conversely, for any normal system
$M = \{M_N\} \in \N(G,R)$, we have
\begin{equation}\label{infl}
\Infl M = \lim_{\overset{N}{\gets}}\Infl^NM_N,
\end{equation}
where the limit is over $N \in \Nn(G)$ and with respect to the
transition maps adjoint to \eqref{norm.trans}. We note that by
Lemma~\ref{infl.le}, all the transition maps in the inverse system
$\Infl^NM_N$ are surjective. For any $M \in \wM(G,R)$, the
adjunction map $M \to \Infl(\Phi(M))$ is the natural map
\begin{equation}\label{can.filt.mack}
M \to \Infl(\Phi(M)) \cong \lim_{\overset{N}{\gets}}\Infl^N\Phi^NM,
\end{equation}
induced by the adjunction maps \eqref{infl.adj}.

\begin{prop}\label{can.pro.prop}
\begin{enumerate}
\item For any normal cofinite subgroup $N \subset G$, the inflation
  functor $\Infl^N$ is fully faithful, and the inflation functor
  $\Infl$ of \eqref{infl} is also fully faithful.
\item For any $R$-valued $G$-Mackey profunctor $M \in \wM(G,R)$, the
  natural map \eqref{can.filt.mack} is surjective.
\item If $\Phi(M) = 0$ for some $M \in \wM(G,R)$, then $M=0$.
\end{enumerate}
\end{prop}

\proof{} For any normal cofinite $N \subset G$ with the quotient $W
= G/N$, let $\wM_N(G,R) \subset \wM(G,R)$ be the full subcategory
spanned by Mackey profunctors $M$ such that $M(S) = 0$ whenever $S^N
= \emptyset$. Then as in the proof of Lemma~\ref{infl.le},
$\wM_N(G,R)$ is equivalent to $\M(W,R)$, with an equivalence given
by restriction to $\Gamma_W \subset \wGamma_G$, and the inverse
equivalence given by $\Infl^N$. In particular, $\Infl^N$ is fully
faithful. Moreover, again as in the proof of Lemma~\ref{infl.le},
for any $M \in \wM(G,R)$ we have exact sequences \eqref{phi.exp},
where $S$ and $S'$ are admissible $G$-sets.

To prove that $\Infl$ is fully faithful, we have to prove that the
adjunction map $\Phi \circ \Infl \to \Id$ is an isomorphism. In
other words, for any normal system $M = \{M_N\} \in \N(G,R)$ with
$\wh{M} = \Infl(M)$, and any normal cofinite subgroup $\overline{N}
\subset G$ with quotient $W = G/\overline{N}$, we have to show that
the map $\Phi^{\overline{N}}(\wh{M}) \to M_{\overline{N}}$ is an
isomorphism. By \eqref{phi.exp}, this is equivalent to proving that
for any $S \in \Gamma_W \subset \wGamma_G$, the kernel of the
natural surjective map
\begin{equation}\label{wh.M.eq}
\wh{M}(S) \to M_{\overline{N}}(S)
\end{equation}
is the union of the images of the maps $f_*$, with $f$ as in
\eqref{phi.exp}. Indeed, by definition, we have
$$
\wh{M}(S) = \lim_{\overset{N}{\gets}}M_N(S),
$$
where the limit is taken over all cofinite normal subgroups $N
\subset G$, and by \eqref{phi.exp}, for any $N' \subset N \subset
G$, the kernel of the map $M_{N'}(S) \to M_N(S) \cong
\Phi^{N/N'}(M_{N'})(S)$ consists of elements of the form
$f^N_*(m_N)$, $m_N \in M_{N'}(S_N)$, $f:S_N \to S$, where $S_N \in
\wGamma_G$ is an admissible $G$-set with no elements fixed under
$N$. Then by induction on normal cofinite subgroups $N \subset
\overline{N}$, any element $\wt{m}$ in the kernel of the map
\eqref{wh.M.eq} can be represented as a series
\begin{equation}\label{ser.eq}
\wt{m} = \sum_{N \subset \overline{N}}f^N_*(\wt{m}_N)
\end{equation}
for some $\wt{m}_N \in \wh{M}(S_N)$, $f^N:S_N \to S$, $S_N^N =
\emptyset$. But then the union
\begin{equation}\label{S.copr}
S' = \coprod_{N}S_N
\end{equation}
is an admissible $G$-set, with a natural map $f:S' \to S$ and no
elements fixed under $\overline{N}$, and $\wt{m}$ is the image of
the element
$$
\left(\prod_N \wt{m}_N\right) \in \wh{M}(S')
$$
under the map $f_*$.

The argument for \thetag{ii} is exactly the same: for any admissible
$G$-set $S$, every element $\wt{m}$ in the target of the map
\eqref{can.filt.mack} evaluated at $S$ can be represented as a
series \eqref{ser.eq}, with $\overline{N} = G$, and if we take $S'$
as in \eqref{S.copr}, with the map $f:S' \to S$, then $\wt{m}$ is
the image of the element
$$
f_*\left(\prod_N \wt{m}_N\right) \in M(S')
$$
under the canonical map \eqref{can.filt.mack}.

Finally, for \thetag{iii}, assume given $M \in \wM(G,R)$ with
trivial $\Phi(M)$. Then by \eqref{phi.exp}, for any $S \in
\wGamma_G$, any element $m \in M(S)$, and any normal cofinite
subgroup $N \subset G$, there exists a map $f:S_N \to S$ in
$\wGamma_G$ and an element $m_N \in M(S_N)$ such that $m = f_*(m_N)$,
while $S^N_N$ is empty. By induction, we can choose a decreasing
sequence of cofinite normal subgroups $N_i \subset G$, $i \geq 1$,
admissible $G$-sets $S_i \in \wGamma_G$ and elements $m_i \in
M(S_i)$ such that $N_{i+1} \subset N_i$, the intersection of all the
subgroups $N_i$ is empty, $S_i^{N_i}$ is empty, and $f^i_*(m_i) =
m_{i-1}$ for some maps $f^i:S_i \to S_{i-1}$, $i \geq 1$, where we
let $S_0=S$ and $m_0=m$. Note that since $\cap N_i = \emptyset$,
$S_i^N$ is empty for any fixed cofinite normal subgroup $N \subset
G$ and almost all $i$. Therefore the union
$$
S' = \coprod_{i \geq 1}S_i
$$
is an admissible $G$-set, and so is $S'_0 = S' \copr S$. We have the
natural projection $p:S'_0 \to S$ and two maps $i,f:S' \to S'_0$ --
the natural embedding $i$, and the disjoint union $f$ of all the
maps $f^i$. Moreover, we have $p \circ f = p \circ i$. It remains to
notice that if we let
$$
\wt{m} = \prod_{i \geq 0}m_i \in \prod_{i \geq 0}M(S_i) = M(S'_0),
$$
then we have $f_*(\wt{m}) = \wt{m}$ and
$$
\begin{aligned}
m &= p_*(\wt{m}) - p_*(i_*(\wt{m})) = p_*(f_*(\wt{m})) - (p \circ
i)_*(\wt{m}) =\\
&= (p \circ f)_*(\wt{m}) - (p \circ i)_*(\wt{m}) = 0.
\end{aligned}
$$
Since $S \in \wGamma_G$ was arbitrary, this means that $M = 0$.
\endproof

\subsection{Separated profunctors.}\label{pro.fp.subs}

We note that Proposition~\ref{can.pro.prop}~\thetag{iii} does not
imply that the canonical map \eqref{can.filt.mack} is always an
isomorphism -- {\em a priori}, the functor $\Phi$ is only
right-exact, so that the kernel of the map \eqref{can.filt.mack} can
have a non-trivial image under this functor. We do not know whether
this can happen, and for which groups. Therefore we introduce the
following.

\begin{defn}\label{sepa.def}
An $R$-valued $G$-Mackey profunctor $M$ is {\em separated} if the
natural surjective map \eqref{can.filt.mack} is an isomorphism.
\end{defn}

We denote by $\wM_s(G,R) \subset \wM(G,R)$ the full subcategory
spanned by separated $G$-Mackey profunctors. By
Proposition~\ref{can.pro.prop}~\thetag{i}, the inflation functor
$\Infl$ induces an equivalence between $\wM_s(G,R)$ and the category
$\N(G,R)$ of normal systems in the sense of
Definition~\ref{ada.def}. The subcategory $\wM_s(G,R)$ is closed
under subobjects (indeed, in the definition, it suffices to require
that the map \eqref{can.filt.mack} is injective, and this property
is inherited by subobjects). We do not know whether $\wM_s(G,R)
\subset \N(G,R)$ is closed under quotients or extensions.

\begin{remark}
Note that while $\wM_s(G,R) \cong \N(G,R)$ has infinite sums, they
do {\em not} commute with evaluation -- for any $S \in \wGamma_G$
and $M_i \in \wM(G,R)$, $i \geq 0$, the map
$$
\bigoplus_i M_i(S) \to \left( \bigoplus_i M_i\right)(S)
$$
is not an isomorphism (the right-hand side is a certain completion of
the left-hand side).
\end{remark}

\begin{defn}
The {\em canonical filtration} on a separated $G$-Mackey profunctor
$M \in \wM_s(G,R)$ is the descreasing filtration $F^NM$ indexed by
normal cofinite subgroups $N \subset G$ such that $F^NM \subset M$
is the kernel of the adjunction map
$$
M \to \Infl^N\Phi^NM.
$$
\end{defn}

By Proposition~\ref{can.pro.prop}, a separated Mackey profunctor is
automatically complete with respect to the canonical filtration.

\begin{lemma}\label{Add.pro.corr}
The embedding $\wM_s(G,R) \subset \Fun(Q\wGamma_G,R)$ admits a
left-ad\-joint additivization functor $\Add:\Fun(Q\wGamma_G,R) \to
\wM^s(G,R)$.
\end{lemma}

\proof{} For any two normal subgroups $N' \subset N \subset G$, the
fixed points functor $\Phi^{N/N'}$ preserves additivity, thus
commutes with the addivization functor $\Add$. Therefore for any
object $E \in \Fun(Q\hGamma_G,R)$, we can apply addivization
fiberwise with respect to the projection $Q\hGamma_G \to \Nn(G)$ and
obtain an object $\Add(E) \in \Fun(Q\hGamma_G,R)$ such that
$\Add(E)_N$ is additive for any $N \subset G$. Moreover, if the
transition morphisms \eqref{norm.trans} for $E$ were isomorphisms,
then the transition morphisms for $\Add(E)$ are also isomorphisms,
so that $\Add(E)$ is a normal system. Now to define the
additivization functor $\wM(G,R) \to \wM_s(G,R) \cong \N(G,R)$, it
suffices to set
$$
\Add(E) = \Infl(\Add(\Phi(E)))
$$
for any $E \in \Fun(Q\wGamma_G,R)$. We then have a natural morphism
$\alpha_E:E \to \Add(E)$, $\Add(E)$ is additive, and
$\Add(\alpha_E)$ is an isomorphism for any $E$, so that $\Add$ is
indeed left-adjoint to the embedding $\wM_s(G,R) \subset \wM(G,R)$.
\endproof

We can now extend the material of Subsection~\ref{ma.fp.subs} and
Subsection~\ref{ma.prod.subs} to Mackey profunctors. First, we note
that as in the Subsection~\ref{ma.prod.subs}, we have the adjoint
pair of functors
\begin{equation}\label{rho.gamma.eq}
\rho^H:\wGamma_G \to \wGamma_H, \qquad \gamma_H:\wGamma_H \to
\wGamma_H,
\end{equation}
for any cofinite subgroup $H \subset G$, and the functor
$$
Q\rho^H_! \cong Q\gamma_H^*:\Fun(\wGamma_G,R) \to \Fun(\wGamma_H,R)
$$
obviously preserves additivity, thus induces a categorical fixed
points functor
\begin{equation}\label{psi.pro.eq}
\Psi^H:\wM(G,R) \to \wM(H,R).
\end{equation}
As in \eqref{psi.wt}, the functor $\Psi^H$ can be refined to a
functor $\wt{\Psi}^H:\wM(G,R) \to \wM(H,R[Z_H])$. We also have the
isomorphism \eqref{phi.psi.iso}. Therefore in particular, the
functors $\Psi^H$ and $\wt{\Psi}^H$ send separated Mackey
profunctors to separated ones.

Moreover, even if a subgroup $H \subset G$ is not cofinite, we still
have a well-defined functor $\phi^H:\wGamma_G \to \wGamma_W$ of
\eqref{phi.H}, where $W = N_H/H$, and $N_H \subset G$ is the
normalizer of $H \subset G$. Therefore we have the inflation functor
$$
\Infl^H = Q(\phi^H)^*:\wM(W,R) \to \wM(G,R).
$$
We also have the inflation functor $\Infl^H:\N(W,R) \to \N(G,R)$,
and it commutes with the functor $\Infl$ of
Proposition~\ref{can.pro.prop}, so that $\Infl^H$ sends separated
Mackey profunctors to separated ones. On the subcategories of
separated Mackey profunctors, $\Infl^H$ has a left-adjoint geometric
fixed points functor
$$
\Phi^H = \Add \circ S(\phi^H)_!:\wM_s(G,R) \to \wM_s(W,R),
$$
where $\Add$ is the additivization functor provided by
Lemma~\ref{Add.pro.corr}. The isomorphism \eqref{phi.psi.iso} also
holds, as long as $H \subset G$ is cofinite.

Finally, we note that for any rings $R_1$, $R_2$,
\eqref{ma.produ} with the addivization functor $\Add$ provided by
Lemma~\ref{Add.pro.corr} defines an associative commutative product
\begin{equation}\label{pro.ma.pro}
\wM_s(G,R_1) \times \wM_s(G,R_1) \to \wM_s(G,R_1 \otimes R_2).
\end{equation}
The fixed point functors $\Phi^H$, $\Psi^H$ are tensor functors with
respect to this product, and the completed Burnside Mackey functor
$\wA \in \wM_s(G,\Z)$ is the unit object. If $R_1=R_2=R$ is a
commutative ring, we can compose the product with the natural
functor induced by the product map $R \otimes R \to R$, and this
turns $\wM_s(G,R)$ into a symmetric monoidal unital category.

\section{Generalities on the $S$-construction.}\label{S.sec}

To proceed further, we need to recall and generalize slightly a
certain construction introduced in \cite[Section 4]{mackey}. It
bears some similarity to Waldhausen's $S$-construction used in the
definition of algebraic $K$-theory, so we denote it by the same
letter. In fact, one can introduce an even more general construction
that would include both Waldhausen's construction and the
construction of \cite[Section 4]{mackey} as particular
cases. However, we will not do it here; we restrict our attention to
what is really needed for the theory of Mackey functors.

\subsection{The constructions.}

We start with the following general definition.

\begin{defn}\label{I.def}
A class of morphisms $I$ in a category $\C$ is {\em admissible} if
it is closed under compositions, contains all identity maps, and for
any two maps $f:c_1 \to c$, $i:c_2 \to c$ in $\C$ with $i$ in
$I$, there exists a cartesian square
\begin{equation}\label{sq}
\begin{CD}
c_{12} @>{f'}>> c_2\\
@V{i'}VV @VV{i}V\\
c_1 @>{f}>> c
\end{CD}
\end{equation}
in $\C$ with $i'$ in $I$. Given an admissible class $I$, we denote
by $\C_I \subset \C$ the subcategory with the same objects, and
those morphisms between them that lie in the class $I$.
\end{defn}

\begin{exa}\label{S.Q.exa}
If a category $\C$ has all fibered products, then the class of all
maps is admissible. We will denote this class by $\Id$.
\end{exa}

\begin{exa}\label{S.I.exa}
In Example~\ref{S.Q.exa}, we can also let $I$ consist of all
monomorphisms in $\C$. We denote this class by $\Inj$.
\end{exa}

\begin{defn}\label{Q.def}
  Assume given a small category $\C$ with an admissible class of
  maps $I$.  For any two objects $c,c' \in \C$, the {\em category
    $\Q_I\C(c,c')$} is the category of diagrams
\begin{equation}\label{dom}
\begin{CD}
c @<{i}<< \wt{c} @>{f}>> c'
\end{CD}
\end{equation}
in $\C$, with $i$ in $I$, and isomorphisms between such
diagrams. The {\em category $Q_I\C$} is the category with the same
objects as $\C$, with morphisms from $c$ to $c'$ given by
isomorphism classes of objects of the category $\Q_I\C(c,c')$, and
with composition given by pullbacks.
\end{defn}

To simplify notation, in the situation of Example~\ref{S.Q.exa}, we
will denote $\Q\C(-,-) = \Q_{\Id}\C(-,-)$ and $Q\C = Q_{\Id}\C$.

\begin{exa}
  Consider the situation of Example~\ref{S.Q.exa} with
  $\C=\Gamma_G$, the category of finite $G$-sets for some group
  $G$. Then $Q\C = Q\Gamma_G$ is the category considered in
  Section~\ref{mack.class.sec}.
\end{exa}

Note that the categories $\Q^I\C(c,c')$ of Definition~\ref{Q.def}
fit together to form a $2$-category $\Q_I\C$; the category $Q_I\C$
is only a truncation of this $2$-category. In general, this
truncation loses some information, and the resulting category may
not be the right one to consider. In particular, this is the case
for the category $Q\Gamma_G$ -- diagrams \eqref{domik} can have
non-trivial automorphisms, and when we pass to $Q\Gamma_G$, we
completely forget them. Thus when considering derived Mackey
functors, it would be more natural to work not with functors from
$Q\Gamma_G$ to complexes of modules over a ring $R$, but with
functors from $\Q\Gamma_G$ to such complexes. Of course, to do this,
one has to define functors from a $2$-category to complexes in the
correct way. This has been done in \cite{mackey}, in several
different ways that are all proved to be equivalent. The resulting
category $\DM(G,R)$ of ``derived Mackey functors'' of \cite{mackey}
indeed behaves better than the derived category $\D(\M(G,R))$ of the
abelian category $\M(G,R)$.

\medskip

In the present paper, we will adopt the least technically cumbersome
of the approaches of \cite{mackey}, namely, that of \cite[Subsection
4.2]{mackey}. Let us recall the construction and its main properties
(for an heuristic description and motivation, see \cite[Subsection
4.1]{mackey}).

\medskip

The construction uses the simplicial technology overviewed in
Subsection~\ref{simpl.subs}. Recall that any object $[n] \in \Delta$
in the category $\Delta$ of finite non-empty ordinals can be treated
as a small category: a functor $f:[n] \to \C$ to some category $\C$
is the same thing as a diagram
\begin{equation}\label{n-ex}
\begin{CD}
c_1 @>{f_1}>> \dots @>{f_{n-1}}>> c_n
\end{CD}
\end{equation}
in $\C$. A morphism $\alpha:f \to f'$ between two such functors is
given by a commutative diagram
\begin{equation}\label{n-ex-2}
\begin{CD}
c_1 @>{f_1}>> \dots @>{f_{n-1}}>> c_n\\
@V{\alpha_1}VV @. @VV{\alpha_n}V\\
c'_1 @>{f'_1}>> \dots @>{f'_{n-1}}>> c'_n,
\end{CD}
\end{equation}
where the bottom row represents $f'$.

\begin{defn}\label{SC.def}
  For any small category $\C$ equipped with an admissible class of
  morphisms $I$, the {\em category $S_I\C$} is given by the
  following:
\begin{enumerate}
\item objects are pairs $\langle [n],f \rangle$ of $[n] \in \Delta$
  and a functor $f:[n] \to \C$,
\item morphisms from $\langle [n],f \rangle$ to $\langle [n'],f'
  \rangle$ are given by a pair of a morphism $\phi:[n] \to [n']$ and
  a morphism $\alpha:f' \circ \phi \to f$ such that for any $i,j \in
  [n]$, the commutative square
\begin{equation}\label{n-ex-3}
\begin{CD}
c'_{\phi(i)} @>>> c'_{\phi(j)}\\
@V{\alpha_i}VV @VV{\alpha_j}V\\
c_i @>>> c_j
\end{CD}
\end{equation}
induced by \eqref{n-ex-2} is cartesian, and the maps $\alpha_i$,
$\alpha_j$ are in $I$.
\end{enumerate}
\end{defn}

As in Definition~\ref{Q.def}, if $I = \Id$ is the class of all maps,
we will drop it from notation and denote $S\C = S_{\Id}\C$.

By definition, the category $S_I\C$ is equipped with a forgetful
functor $\pi:S_I\C \to \Delta$ sending $\langle [n],f \rangle$ to
$[n]$. This functor is a fibration. Its fiber $(S_I\C)_{[1]}$ over
$[1] \in \Delta$ is naturally identified with the category $\C_I^o$
opposite to the category $\C_I$. For $n \geq 2$, the fiber
$(S_I\C)_{[n]}$ is opposite to the category of diagrams
\eqref{n-ex}, with morphisms given by diagrams \eqref{n-ex-2}, with
$\alpha_i$ in $I$ and cartesian squares. Note that by
Definition~\ref{I.def}, it suffices to require that $\alpha_n$ is in
$I$.

Note that for any two objects $c,c' \in \C = (S_I\C)_{[1]}$, the
category $Q_I\C(c,c')$ of Definition~\ref{Q.def} can be described in
terms of the category $S_I\C$ in the following way. Consider the
diagram
\begin{equation}\label{st.eq}
\begin{CD}
[1] @>{s}>> [2] @<{t}<< [1]
\end{CD}
\end{equation}
in $\Delta$, where $s,t:[1] \to [2]$ send the unique element in $[1]
\in \Delta$ to the initial resp. terminal object of the ordinal $[2]
\in \Delta$. Then $Q_I\C(c,c')$ is equivalent to the category of
diagrams
$$
\begin{CD}
c @>{\wt{s}}>> c_{[2]} @<{\wt{t}}<< c'
\end{CD}
$$
in $S_I\C$ that are sent to \eqref{st.eq} by the projection
$\pi:S_I\C \to \Delta$, and such that $\wt{t}$ is cartesian with
respect to $\pi$. The equivalence sends an object in $\Q_I(c,c')$
represented by a diagram \eqref{dom} to the diagram
\begin{equation}\label{Q.S.eq}
\begin{CD}
c @>>> [\wt{c} \to c'] @<<< c',
\end{CD}
\end{equation}
where $[\wt{c} \to c']$ is an object in $(S_I\C)_{[2]}$ represented
by the arrow $\wt{c} \to c'$.

\begin{defn}\label{sp.defn}
  A map $f:[n] \to [n']$ in the category $\Delta$ is {\em special}
  if it sends the terminal element in $[n]$ to the terminal element
  in $[n']$. For small category $\C$ with an admissible class of
  maps $I$, a morphism $\langle f,\alpha \rangle:\langle [n],f
  \rangle \to \langle [n'],f' \rangle$ is {\em special} if $f =
  \pi(\langle f,\alpha \rangle)$ is special and $\alpha_n:c'_{n'}
  \to c_n$ is an isomorphism (that is $\langle f,\alpha \rangle$ is
  cartesian with respect to the fibration $\pi:S_I\C \to \Delta$).
\end{defn}

\begin{defn}\label{DS.def}
  For any ring $R$, an object $M \in \D(S_I\C,R)$ is {\em special}
  if it can be represented by such a complex $M_\idot$ in
  $\Fun(S_I\C,R)$ that $M_\idot(f)$ is a quasiisomorphism for any
  special map $f$ in $S\C$.  The full subcategory in $\D(S_I\C,R)$
  spanned by special objects is denoted by $\DS_I(\C,R)$.
\end{defn}

\begin{remark}
The notion of a special map in Definition~\ref{sp.defn} is different
from the one used in \cite[Subsection 4.2]{mackey}. However, the
resulting category $\DS_I(\C,R)$ is clearly the same.
\end{remark}

For any $M \in \DS_I(\C,R)$ and any object $c \in \C$, we can
evaluate $M$ at $c$ by considering $c$ as an object of $\C^o =
(S_I\C)_{[1]} \subset S_I\C$; this gives an object $M(c) \in
\D(R)$. For any two objects $c,c' \in \C$, an object in
$\Q_I\C(c,c')$ represented by a digram \eqref{dom} defines a natural
map from $M(c)$ to $M(c')$ -- namely, the composition map
\begin{equation}\label{M.dia}
\begin{CD}
M(c) @>{M(\wt{s})}>> M([\wt{c} \to c']) @>{M(\wt{t})^{-1}}>> M(c')
\end{CD}
\end{equation}
induced by the corresponding diagram \eqref{Q.S.eq}. If $f$
resp. $i$ in \eqref{dom} is the identity map, then we will denote
the map \eqref{M.dia} by $i^*$ resp. $f_*$; in general,
\eqref{M.dia} is the composition $f_* \circ i^*$. This is compatible
with the composition of diagrams, so that every map $M \in
\DS_I(\C,R)$ induces a functor from $Q_I\C$ to the derived category
$\D(R)$.

This construction can be refined in the following way. For any
object $c \in S_I\C$ represented by a diagram \eqref{n-ex}, let
$q(c) = c_n$, and for any morphism $f = \langle f,\alpha \rangle:c
\to c'$ in $S_I\C$, $c \in (S_I\C)_{[n]}$, $c' \in (S_I\C)_{[n']}$,
let $q(f) \in \Q_I\C(c_n,c'_{n'})$ be the diagram
\begin{equation}\label{q.f.eq}
\begin{CD}
c_n @<{\alpha}<< c'_{f(n)} @>>> c'_{n'},
\end{CD}
\end{equation}
where the map on the right is the composition of the natural maps in
the diagram \eqref{n-ex}. Then sending $c$ to $q(c)$ and $f$ to the
isomorphism class of $q(f)$ defines a functor
\begin{equation}\label{q.SQ}
q:S_I\C \to Q_I\C
\end{equation}
such that $q(f)$ is invertible for any special morphism $f$ in
$S_I\C$. Therefore for any ring $R$ and any $M \in \D(Q_I\C,R)$, the
pullback $q^*M \in \D(S_I\Q,R)$ is special, so that we have a natural
pullback functor
\begin{equation}\label{q.st}
q^*:\D(Q_I\C,R) \to \DS_I(\C,R).
\end{equation}

\begin{lemma}\label{M.DM}
  The truncation functors with respect to the standard
  $t$-struc\-ture on $\D(S_I\C,R)$ send special objects to special
  objects, thus induce a $t$-structure on $\DS_I(\C,R) \subset
  \D(S_I\C,R)$ given by
$$
\DS^{\leq i}_I(\C,R) = \DS_I(\C,R) \cap \D^{\leq i}(S_I\C,R) \subset
\DS_I(\C,R)
$$
for any integer $i$. The functor \eqref{q.st} gives an equivalence
$$
\Fun(Q_I\C,R) \cong \Fun(S_I\C,R) \cap \DS_I(\C,R) \subset \D(S_I\C,R)
$$
between the heart of this $t$-structure and the category
$\Fun(Q_I\C,R)$.
\end{lemma}

\proof{} Clear. \endproof

\begin{defn}\label{adm.mor.def}
  Assume given two small categories $\C$, $\C'$ with admissible
  classes of maps $I$, $I'$. A {\em morphism} $\phi:\langle \C,I
  \rangle \to \langle \C',I' \rangle$ is a functor $\phi:\C \to \C'$
  that sends morphisms in $I$ to morphisms in $I'$ and cartesian
  squares \eqref{sq} in $\C$ to cartesian squares in $\C'$. For any
  two such morphisms $\phi,\phi':\C \to \C'$, a map $\alpha:\phi \to
  \phi'$ is {\em $Q$-compatible} if for any object $c \in \C$,
  $\alpha:\phi(c) \to \phi'(c)$ is in $I'$, and for any morphism
  $f:c \to c'$ in $\C$, the commutative square
$$
\begin{CD}
\phi(c) @>{\alpha}>> \phi'(c)\\
@V{\phi(f)}VV @VV{\phi'(f)}V\\
\phi(c') @>{\alpha}>> \phi'(c')
\end{CD}
$$
is cartesian.
\end{defn}

\begin{exa}\label{c.adj}
Assume given a morphism $f:c \to c'$ in a small category $\C$ with
fibered products. Then the functors $f_!$, $f^*$ of \eqref{f.l},
\eqref{f.r} both give morphisms between $\langle \C/c,\Id \rangle$
and $\langle \C/c',\Id \rangle$, and both adjunction maps $\Id \to
f^* \circ f_!$, $f_1 \circ f^*\to \Id$ are $Q$-compatible.
\end{exa}

\begin{lemma}\label{adm.pb.le}
\begin{enumerate}
\item A morphism $\phi:\langle \C,I \rangle \to \langle \C',I'
  \rangle$ in the sense of Definition~\ref{adm.mor.def} induces a
  functor $S(\phi):S_I\C \to S_{I'}\C'$. This functor commutes with
  projections to $\Delta$ and sends special maps to special maps,
  thus induces a pullback functor
$$
S(\phi)^*:\DS_{I'}(\C',R) \to \DS_I(\C,R)
$$
for any ring $R$. Moreover, $\phi$ induces a functor $Q(\phi):Q_I\C
\to Q_{I'}\C'$, and we have $q \circ S(\phi) \cong Q(\phi) \circ q$.
\item For any two such morphisms $\phi,\phi':\langle \C,I \rangle
  \to \langle \C',I' \rangle$, a $Q$-compatible map $\alpha:\phi \to
  \phi'$ induces maps of functors $S(\alpha):S(\phi) \to S(\phi')$,
  $Q(\alpha):Q(\phi) \to Q(\phi')$.
\end{enumerate}
\end{lemma}

\proof{} Clear. \endproof

\subsection{Basic properties.}

On the level of full derived categories, the functor $q^*$ of
\eqref{q.st} is usually {\em not} an equivalence -- $\DS_I(\C,R)$ is
an example of a triangulated category with a $t$-structure that is
not equivalent to the derived category of its heart. From our point
of view, it is $\DS_I(\C,R)$ that is the right category to consider,
and it is the correct replacement of the derived category
$\D(Q_I\C,R)$.

To study it, we need one result essentially proved in \cite[Section
  4]{mackey}. Denote by $\wt{S}_I\C$ the category of diagrams
\begin{equation}\label{codomik}
\begin{CD}
c_1 @<{s_1}<< c_0 @>{s_2}>> c_2
\end{CD}
\end{equation}
in $S_I\C$ with special maps $s_1$, $s_2$, and let
$\pi_1,\pi_2:\wt{S_I\C} \to S_I\C$ be the projections sending such a
diagram to $c_1$ resp. $c_2$. Both projections are
cofibrations, and we have the natural endofunctor
\begin{equation}\label{sp.eq}
L^\hdot\pi_{2!}\pi_1^*:\D(S_I\C,R) \to \D(S_I\C,R)
\end{equation}
of the derived category $\D(S_I\C,R)$. Moreover, cofibrations $\pi_1$, $\pi_2$ admit a common
section $\lambda:S_I\C \to \wt{S}_I\C$ sending $c$ to a diagram
\eqref{codomik} with $c_1=c_2=c_0=c$, $s_1=s_2=\id$, and the
adjunction map $L^\hdot\lambda_!  \circ \lambda^* \to \Id$ induces a
natural map
\begin{equation}\label{sp.adj}
a:\Id \cong L^\hdot\pi_{2!} \circ L^\hdot\lambda_! \circ \lambda^*
\circ \pi_1^* \to L^\hdot\pi_{2!}\pi_1^*.
\end{equation}

\begin{lemma}\label{sp.le}
Assume given a small category $\C$ and an admissible class of
morphisms $I$. Then the functor \eqref{sp.eq} takes values in
$\DS_I(\C,R) \subset \D(S_I\C,R)$ and induces a functor
$$
\Sp:\D(S_I\C,R) \to \DS_I(\C,R)
$$
left-adjoint to the embedding $\DS_I(\C,R) \hookrightarrow
\D(S_I\C,R)$, with the adjunction induced by the natural map
\eqref{sp.adj}.
\end{lemma}

\proof{} Say that a map $f:[n] \to [n']$ is {\em co-special} if it
identifies the ordinal $[n]$ with an initial segment of the
ordinal $[n']$, and say that a map $f$ in $S_I\C$ is {\em
  co-special} if so is its image $\pi(f)$ with respect to the
canonical projection $\pi:S_I\C \to \Delta$. Then by the same
argument as in \cite[Lemma 4.8]{mackey}, the classes of special and
co-special maps form a complementary pair on $S_I\C$ in the sense of
\cite[Definition 4.3]{mackey}, and the result immediately follows
from \cite[Lemma 4.6]{mackey}.
\endproof

\begin{corr}\label{sp.corr}
For any morphism $\phi:\langle \C,I \rangle \to \langle \C',I'
\rangle$ in the sense of Definition~\ref{adm.mor.def}, and any ring
$R$, the pullback functor $S(\phi)^*$ of Lemma~\ref{adm.pb.le}
admits a left-adjoint functor
$$
S(\phi)_!:\DS_I(\C,R) \to \DS_{I'}(\C',R).
$$
\end{corr}

\proof{} It suffices to set $S(\phi)_! = \Sp \circ L^\hdot
S(\phi)_!$, where the functor $\Sp$ is provided by
Lemma~\ref{sp.le}, and $L^\hdot S(\phi)_!:\D(S_I\C,R) \to
\D(S_{I'}\C',R)$ is the derived left Kan extension functor with
respect to $S(\phi):S_I\C \to S_{I'}\C'$.
\endproof

In some special cases, $S(\phi)_!$ is easy to compute. In
particular, assume given two categories $\C$, $\C'$ with fibered
products, equip them with admissible classes $\Id$ of all morphisms,
and assume given a pair of adjoint functors $\phi:\C' \to \C$,
$\psi:\C \to \C'$ such that both send fibered products to fibered
products, and both adjunction maps $\Id \to \psi \circ \phi$, $\phi
\circ \psi \to \Id$ are $Q$-compatible in the sense of
Definition~\ref{adm.mor.def} (such a situation occurs, for instance,
in Example~\ref{c.adj}). Then by Lemma~\ref{adm.pb.le}, the
adjunction maps induce natural maps
$$
S(\psi) \circ S(\phi) \to \Id, \qquad \Id \to S(\phi) \circ S(\psi),
$$
and these maps give an adjunction between $S(\phi)$ and $S(\psi)$,
so that we have a natural isomorphism
\begin{equation}\label{phi.psi.eq}
S(\psi)_! \cong S(\phi)^*:\DS(\C,R) \to \DS(\C',R).
\end{equation}
The specialization functor $\Sp$ does not enter into the picture at
all.

In the general case, to control $S(\phi)_!$, it is useful to have a
reasonably explicit description of objects $\Sp(F)$, $F \in
\D(S_I\C,R)$. Such a description is in fact provided in \cite[Lemma
  4.6]{mackey}. We will need it in one particular case -- namely,
for objects of the form $\Sp(M_c) \in \DS_I(\C,R)$, where $M$ is an
$R$-module, $c \in \C^o \cong (S_I\C)_{[1]} \subset S_I\C$ is an
object of $\C$ considered as an object of $S_I\C$, and $M_c \in
\Fun(S_I\C,R)$ is the functor \eqref{m.c.def}. To obtain such a
description, note that the categories $\Q_I\C(c,-)$ fit together
into a cofibration
\begin{equation}\label{rho.c.eq}
\rho_c:\Q_I^c\C \to S_I\C
\end{equation}
whose fiber over $c' \in S_I\C$ is the category $\Q_I\C(c,q(c'))$,
and whose transition functor corresponding to a morphism $f$ is
given by the composition with the diagram $q(f)$ of \eqref{q.f.eq}.

\begin{lemma}\label{sp.mc}
For any $c \in \C$ and any $R$-module $M$, we have a natural
identification
\begin{equation}\label{sp.exp.eq}
\Sp(M_c) \cong L^\hdot\rho_{c!}M,
\end{equation}
where $\rho_c$ is the cofibration \eqref{rho.c.eq}, and $M \in
\Fun(\Q_I^c\C,R)$ is the constant functor with value $M$.
\end{lemma}

We note that by definition, for any special map $f$ in $S_I\C$, the
transition functor of the cofibration $\rho_c$ is an equivalence of
categories, so that by base change, the right-hand side of
\eqref{sp.exp.eq} is a special object in $\D(S_I\C,R)$. Explicitly,
for any $c' \in S_I\C$, we have a natural identification
\begin{equation}\label{sp.m.c.eq}
\Sp(M_c)(c') \cong H_\idot(\Q_I\C(q(c),q(c')),M),
\end{equation}
where the right-hand side denotes the homology of
$\Q^I\C(q(c),q(c'))$ with coefficients in the constant functor with
value $M$. We also note that by adjunction, we have a natural map
$M_c \to q^*M_{q(c)}$, and since $q^*M_{q(c)}$ is special, it
factors through a map
\begin{equation}\label{sp.q}
\Sp(M_c) \to q^*M_{q(c)}.
\end{equation}
In terms of the identification \eqref{sp.m.c.eq}, this map becomes
the map
$$
H_\idot(\Q_I\C(q(c),q(c')),M) \to H_0(\Q_I\C(q(c),q(c')),M) \cong
M[Q_I\C(q(c),q(c'))]
$$
corresponding to the identfication of $Q_I(q(c),q(c'))$ with the set
of isomorphism classes of objects in the groupoid
$\Q_I(q(c),q(c'))$. In particular, \eqref{sp.q} identifies
$q^*M_{q(c)}$ with the truncation at $0$ of the object $\Sp(M_c)$
with respect to the standard $t$-structure on $\DS_I(\C,R)$.

\proof[Proof of Lemma~\ref{sp.mc}.] As shown in the proof of
\cite[Lemma 4.6]{mackey}, to compute the specialization functor
\eqref{sp.eq}, we can replace $\wt{S}_I\C$ with its subcategory
$\overline{S}_I\C$ spanned by diagrams \eqref{codomik} with $c_0$
lying in $(S_I\C)_{[1]} \subset S_I\C$. To apply this to $M_c$, use
its explicit description \eqref{m.c.eq}. Then by base change, we
obtain an isomorphism
$$
\Sp(M_c) \cong L^\hdot\wt{\rho}_{c!}M,
$$
where $\wt{\rho}_c:\wt{\Q}_I^c \to S_I\C$ is a cofibration whose
fiber over $c' \in S_I\C$ is the category of diagrams
$$
\begin{CD}
c @>{p}>> \wt{c} @<{i}<< q(c')
\end{CD}
$$
in $S_I\C$ with special $i$ and co-special $p$ (the meaning of
co-special is the same as in the proof of Lemma~\ref{sp.le}). But
the category $\Q_I^c$ is actually embedded into $\wt{\Q}_I^c$, with
the embedding sending a diagram \eqref{dom} to the corresponding
diagram \eqref{Q.S.eq}. It remains to notice that the embedding
$\iota:\Q_I^c \to \wt{\Q}_I^c$ admits a right-adjoint functor, so
that $L^\hdot\iota_!M \cong M$, and $L^\hdot\wt{\rho}_{c!}M \cong
L^\hdot\rho_{c!}M$.
\endproof

Lemma~\ref{sp.mc} has one immediate corollary. Say that a pair
$\langle \C,I \rangle$ of a small category $\C$ and an admissible
class of morphisms $I$ is {\em discrete} if all the categories
$\Q_I(c,c')$ of Definition~\ref{Q.def} are essentially discrete
categories (that is, groupoids with trivial automorphism groups).

\begin{corr}\label{S.discr.corr}
  Assume given a small category $\C$ and a class of morphisms $I$,
  and assume that the pair $\langle \C,I \rangle$ is discrete. Then
  the functor $q^*$ of \eqref{q.st} is an equivalence of categories.
\end{corr}

\proof{} Under the assumptions of the Corollary, the cofibration
$\rho_c$ of \eqref{rho.c.eq} is actually discrete for any $c \in
\C$. Therefore by base change, the functors $\rho_{c!}$ is exact,
and we have
$$
\Sp(M_c) \cong L^\hdot\rho_{c!}M \cong \rho_{c!}M
$$
for any $R$-module $M$. Moreover, it is easy to see that
$$
\rho_{c!}M \cong q^*M_{\wt{c}},
$$
where $\wt{c}$ denotes $c$ considered as an object in
$Q_I\C$. But by adjunction, we have
$$
L^\hdot q_!\Sp(M_c) \cong M_{\wt{c}},
$$
so that the adjunction map $L^\hdot q_!q^*E \to E$ is an isomorphism
for an object $E \in \D(Q_I\C,R)$ of the form $E = M_{\wt{c}}$, $c
\in \C$. Since such objects generate the category $\D(Q_I\C,R)$, we
have $L^\hdot q_!q^* \cong \Id$, so that $q^*$ is fully
faithful. Since objects of the form $\Sp(M_c)$, $c \in \C$, generate
the category $\DS_I(\C,R)$, the functor \eqref{q.st} is also
essentially surjective.
\endproof

\subsection{A base change lemma.}\label{bc.subs}

We will also need one general result that in some cases greatly
simplifies the computation of the functors $S(\phi)_!$ provided by
Corollary~\ref{sp.corr}.

The setup is the following. Assume given a small category $\C$ with
fibered products and a full subcategory $\C' \subset \C$ closed
under fibered products. Moreover, assume that the embedding $\C'
\subset \C$ admits a right-adjoint functor
$$
\phi:\C \to \C'.
$$
Then by adjunction, $\phi$ automatically sends
fibered products to fibered products, thus gives a morphism
$\phi:\langle \C,\Id \rangle \to \langle \C',\Id \rangle$ in the
sense of Definition~\ref{adm.mor.def}.

Moreover, assume also given an admissible class of maps $I$ in
$\C$. Since $\C' \subset \C$ is closed under taking fibered
products, $I' = I \cap \C'$ is then an admissible class of maps in
$\C'$, with $\C'_{I'} = \C' \cap \C_I \subset \C_I$. Assume further
that $\phi$ sends morphisms in $I$ to morphisms in $I'$, so that it
induces a functor $\phi:\C_I \to \C'_{I'}$ and a morphism
$\phi:\langle \C,I \rangle \to \langle \C',I' \rangle$. In addition,
assume that $\langle \C,I \rangle$, and therefore also $\langle
\C',I' \rangle$, is discrete in the sense of
Corollary~\ref{S.discr.corr}, so that we have natural equivalences
of categories
$$
\DS_I(\C,R) \cong \D(Q_I\C,R), \qquad \DS_{I'}(\C',R) \cong
\D(Q_{I'}\C',R).
$$

\begin{remark}
  If $I = \Inj$ is the class of all monomorphisms in $\C$, as in
  Example~\ref{S.I.exa}, then all these assumption on $I$ are
  automatically satisfied.
\end{remark}

Then we have tautological embedding morphisms $i:\langle \C,I
\rangle \to \langle \C,\Id \rangle$, $i:\langle \C',I' \rangle \to
\langle \C',\Id \rangle$, and $i \circ \phi \cong \phi \circ i$, so
that we have an isomorphism
$$
S(\phi)^* \circ S(i)^* \cong S(i)^* \circ S(\phi)^*
$$
of the corresponding pullback functors of Lemma~\ref{adm.pb.le}. By
adjunction, it induces a base change map
\begin{equation}\label{phi.pb.eq}
S(\phi)_! \circ S(i)^* \to S(i)^* \circ S(\phi)_!,
\end{equation}
where the adjoint functors $S(\phi)_!$ are those of
Corollary~\ref{sp.corr}.

\begin{prop}\label{bc.prop}
  Under the assumptions above, assume further that for any $c \in
  \C$, the adjunction map $\phi(c) \to c$ is in the class $I$. Then
  the base change map \eqref{phi.pb.eq} is an isomorphism.
\end{prop}

Before we prove this, we need a technical lemma. Fix an object $c
\in \C$, and let $(\C/c)_I$ be the category of objects $c' \in \C$
equipped with a map $f:c' \to c$, with morphisms from $f_1:c'_1 \to
c$ to $f_2:c'_2 \to c$ given by morphisms $i:c'_1 \to c_2'$ in the
class $I$ such that $f_2 \circ i = f_1$. Then we have a natural
forgetful functor $(\C/c)_I \to \C_I$. Composing the opposite
functor with the obvious embedding $\C_I^o \to Q_I\C$, we obtain a
functor
$$
j_c:(\C/c)_I^o \to Q_I\C.
$$

\begin{lemma}\label{jc.le}
For any $c \in \C$ and any $R$-module $M$, we have a natural
isomorphism
$$
S(i)^*\Sp(M_c) \cong L^\hdot j_{c!}M \in \DS_I(\C,R) \cong \D(Q_I\C,R),
$$
where in the right-hand side, $M \in \Fun((\C/c)^{Io},R)$ stands for
the constant functor with value $M$.
\end{lemma}

\proof{} To compute $L^\hdot j_{c!}M$, we can use the decomposition
\eqref{facto.co} and the isomorphism \eqref{facto.t}. Spelling out
the definition, we check that the comma-category $j_c\backslash Q_I\C$ is the
category of diagrams
\begin{equation}\label{Q.jc}
\begin{CD}
c @<{f}<< c_1 @<{i}<< c_2 @>{p}>> c_3
\end{CD}
\end{equation}
in $\C$, with $i$ in the class $I$, and morphisms between such
diagrams are given by commutative diagrams
$$
\begin{CD}
c @<{f}<< c_1 @<{i}<< c_2 @>{p}>> c_3\\
@| @A{\alpha_1}AA @VV{\alpha_2}V @VV{\alpha_3}V\\
c @<{f'}<< c'_1 @<{i'}<< c'_2 @>{p}>> c'_3
\end{CD}
$$
in $\C$ such that $\alpha_1$ is in $I$ and $\alpha_2$ is an
isomorphism. The cofibration $t:j_c\backslash Q_I\C \to Q_I\C$ sends a diagram
\eqref{Q.jc} to $c_3 \in Q^I\C$.

Inside the category $j_c\backslash Q_I\C$, we have a full
subcategory $\overline{j_c\backslash Q_I\C}$ spanned by diagrams
\eqref{Q.jc} with invertible $i$, and the embedding
$\overline{j_c\backslash Q_I\C} \subset j_c\backslash Q_I\C$ has a
right-adjoint $\rho:j_c\backslash Q_I\C \to \overline{j_c\backslash
  Q_I\C}$. Therefore we have canonical identifications
$$
L^\hdot j_{c!}M \cong L^\hdot t_! M \cong L^\hdot \overline{t}_!
\rho^*M \cong L^\hdot \overline{t}_!M,
$$
where $\overline{t}$ is the restriction of the projection $t$ to
$\overline{j_c\backslash Q_I\C}$. It remains to notice that we have
a cartesian square
$$
\begin{CD}
\Q^I_c\C @>>> \overline{j_c\backslash Q_I\C}\\
@V{\rho_c}VV @VV{\overline{t}}V\\
S^I\C @>{q}>> Q^I\C
\end{CD}
$$
of categories and functors, and apply Lemma~\ref{sp.mc} and the base
change isomorphism \eqref{bc.iso}.
\endproof

\proof[Proof of Proposition~\ref{bc.prop}.] We have to prove that
for any $E \in \D(S\C,R)$, the map
\begin{equation}\label{bc-toprove.eq}
S(\phi)_!S(i)^*\Sp(E) \to S(i)^*S(\phi)_!\Sp(M)
\end{equation}
induced by \eqref{phi.pb.eq} is an isomorphism. Since the derived
category $\D(S\C,R)$ is generated by objects of the form $M_c$, $c
\in \C$, $M$ an $R$-module, it suffices to prove this for $E =
M_c$. By adjunction, we have $S(\phi)_!\Sp(M_c) \cong
\Sp(M_{\phi(c)})$, and under the identification of
Lemma~\ref{jc.le}, \eqref{bc-toprove.eq} becomes the natural map
$$
L^\hdot\phi_!L^\hdot j_{c!}M \to L^\hdot j_{\phi(c)!}M.
$$
Moreover, the functor $\phi$ induces a natural functor
$\phi_c:(\C/c)_I^o \to (\C'/\phi(c))_{I'}^o$ such that $\phi \circ
j_c \cong j_{\phi(c)} \circ \phi_c$, so that it suffices to prove
that $L^\hdot\phi_{c!}M$ is the constant functor with value $M$.

But the assumptions of the proposition imply that if we denote by
$\iota:\C' \to \C$ the embedding functor, then the adjunction maps
\begin{equation}\label{adj.ci}
\Id \to \iota \circ \phi, \qquad \phi \circ \iota \to \Id
\end{equation}
both lie pointwise in the class $I$. Therefore they also induce an
adjunction between the embedding $\C'_{I'} \subset \C_I$ and the
functor $\phi:\C_I \to \C'_{I'}$. When we pass to the opposite
categories, the embedding $\C^{'o}_{I'} \subset \C^o_I$ becomes
right-adjoint to $\phi:\C^o_I \to \C^{'o}_{I'}$.

It remains to notice that the embedding $\iota:\C'_{I'} \to \C_I^o$
extends to a functor $\iota_c:(\C'/\phi(c))_{I'}^o \to (\C/c)_I^o$
sending $f:c' \to \phi(c)$ to its composition with the adjunction
map $\phi(c) \to c$. Moreover, the adjunction maps \eqref{adj.ci}
induce maps of functors
$$
\Id \to \iota_c \circ \phi_c, \qquad \phi_c \circ \iota_c \to \Id,
$$
and these maps give an adjunction between $\phi_c$ and
$\iota_c$. Thus $L^\hdot\phi_{c!}M \cong \iota_c^*M \cong M$ is
indeed the constant functor.
\endproof

\section{Additivization.}\label{add.sec}

We will also need one result concerning a more special
situation -- that of a small category $\C$ that admits finite
coproducts.

\subsection{The setup.} We start with a formal definition and the
statement.

\begin{defn}\label{copr.def}
A pair $\langle \C,I \rangle$ of a small category $\C$ and an
admissible class of morphisms $I$ in $\C$ {\em has finite
  coproducts} if the following holds.
\begin{enumerate}
\item The category $\C$ has finite coproducts and an initial object
  $\emptyset \in \C$.
\item For any $c \in \C$, the sole morphism $\emptyset \to c$ is a
  monomorphism and lies in the class $I$.
\item The coproduct of two morphisms in the class $I$ is in the
  class $I$, and the coproduct of two cartesian squares \eqref{sq}
  in $\C$ is cartesian.
\end{enumerate}
\end{defn}

\begin{defn}\label{add.S.def}
  Assume that a small category $\C$ has finite coproducts, and an
  admissible class of morphisms $I$ in $\C$ is compatible with
  coproducts in the sense of Definition~\ref{copr.def}. For any ring
  $R$, an object $E \in \DS_I(\C,R)$ is {\em additive} if for any
  two objects $c_1,c_2 \in \C$ with the natural maps $i_1:c_1 \to
  c_1 \copr c_@$, $i_2:c_2 \to c_1 \copr c_2$, the natural map
\begin{equation}\label{add.S.eq}
\begin{CD}
M(c_1 \copr c_2) @>{i_1^* \oplus i_2^*}>> M(c_1) \oplus M(c_2)
\end{CD}
\end{equation}
is an isomorphism. The full subcategory in $\DS_I(\C,R)$ spanned by
additive objects is denoted by $\DS^{add}_I(\C,R) \subset
\DS_I(\C,R)$.
\end{defn}

\begin{lemma}\label{c.c.le}
  In the situation of Definition~\ref{add.S.def}, an object $M \in
  \DS_I(\C,R)$ is additive if and only if the map \eqref{add.S.eq}
  is an isomorphism whenever $c_1 \cong c_2$.
\end{lemma}

\proof{} Definition~\ref{copr.def}~\thetag{ii} and \thetag{iii}
immediately imply that for any two objects $c_1,c_2 \in \C$ with
coproduct $c = c_1 \copr c_2$, the natural maps $i_1:c_1 \to c$,
$i_2:c_2 \to c$ are in the class $I$. Moreover,
Definition~\ref{copr.def}~\thetag{ii} implies that we have
$\emptyset \times_{c'} \emptyset = \emptyset$ for any $c' \in \C$,
and then Definition~\ref{copr.def}~\thetag{iii} implies that $c_1
\times_c c_1 \cong c_1$ and $c_2 \times_c c_2 \cong c_2$ (in other
words, $i_1$ and $i_2$ are monomorphisms). Therefore for any $M \in
\DS_I(\C,R)$, the endomorphisms $i_1^* \circ i_{1*}:M(c_1) \to
M(c_1)$, $i_2^* \circ i_{2*}:M(c_2) \to M(c_2)$ are equal to the
identity maps. Moreover, if we consider the map
$$
i = i_1  \copr i_2:c = c_1 \copr c_2 \to c \copr c,
$$
then $i^* \circ i_*:M(c) \to M(c)$ is also the identity map. Thus the
compositions 
$$
\begin{aligned}
p_1 &= i_{1*} \circ i_1^*:M(c) \to M(c),\\
p_2 &= i_{2*} \circ i_2^*:M(c) \to M(c),\\
p &= i_* \circ i^*:M(c \copr c) \to M(c \copr c),
\end{aligned}
$$
are idempotent endomorphisms, with images $M(c_1)$, $M(c_2)$,
$M(c)$. But the natural map $M(c \copr c) \to M(c) \oplus M(c)$
intertwines $p$ with $p_1 \oplus p_2$. Therefore if this map is an
isomorphism, then so is the map \eqref{add.S.eq}.
\endproof

The additivity condition is obviously preserved under truncation
functors with respect to the $t$-structure of Lemma~\ref{M.DM} on
$\DS_I(\C,R)$ (and those are by definition the truncation functors
with respect to the standard $t$-structure on
$\D(S_I\C,R)$). Therefore the subcategory $\DS^{add}_I(\C,R) \subset
\DS_I(\C,R)$ inherits a natural $t$-structure.

\begin{prop}\label{add.prop}
  In the situation of Definition~\ref{add.S.def}, the embedding
$$
\DS^{add}_I(\C,R) \subset \DS_I(\C,R)
$$
admits a left-adjoint additivization functor
$$
\Add:\DS_I(\C,R) \to \DS_I^{add}(\C,R)
$$
that is right-exact with the respect to the natural $t$-structures.
\end{prop}

One example of a category with coproducts is the category $\Gamma_G$
of finite $G$-sets for some group $G$, and the class of all maps is
compatible with coproducts in the sense of
Definition~\ref{copr.def}. In this case, Proposition~\ref{add.prop}
has been essentially proved in \cite{mackey}. However, the proof is
pretty roundabout, and the details are hard to trace. Therefore we
reprove Proposition~\ref{add.prop} from scratch using the
additivization technique of \cite[Subsection 3.2]{cartier} (which is
essentially due to T. Pirashvili). The proof is split into a
sequence of lemmas and takes up the rest of this section.

\subsection{Finite sets.}

We first consider the case $\C = \Gamma$, the category of finite
sets, with the class $\Inj$ of all monomorphisms as in
Example~\ref{S.I.exa}. These obviously satisfy all the assumptions
of Definition~\ref{I.def}. The category $Q_{\Inj}\Gamma$ is
naturally equivalent to the category $\Gamma_+$ of finite pointed
sets, that is, finite sets with a distiguished element. The
equivalence sends a pointed set $S$ with the distinguished element
$o \in S$ to the complement $\overline{S} = S \setminus \{o\}$, and
a pointed map $f:S \to S'$ goes to the diagram
\begin{equation}\label{partial.map}
\begin{CD}
\overline{S} @<<< f^{-1}(\overline{S}') @>>> \overline{S}'.
\end{CD}
\end{equation}
Moreover, the pair $\langle \Gamma,\Inj \rangle$ is discrete. Thus
Corollary~\ref{S.discr.corr} applies to $\langle \Gamma, \Inj \rangle$,
so that for any ring $R$, we have a canonical equivalence
$$
\DS_{\Inj}(\Gamma,R) \cong \D(\Gamma_+,R).
$$
For any integer $n \geq 0$, let us denote by $[n] \in \Gamma$ the
set with $n$ elements, and in keeping with our convention
\eqref{s.pl}, we denote by $[n]_+ \in \Gamma_+$ the pointed set with
$n$ non-distiguished elements. Moreover, let $j_n:\ppt \to \Gamma_+$
be the embedding onto $[n]_+$, and let $T_n = j_{n!}\Z = \Z_{[n]_+}
\in \Fun(\Gamma_+,\Z)$ be the object represented by $[n]_+ \in
\Gamma$. Then since $[0]_+$ is canonically a retract of $[1]_+$,
$T_0$ is canonically a retract of $T_1$, so that we have a canonical
direct sum decomposition
\begin{equation}\label{T.1.0}
T_1 = T_0 \oplus T
\end{equation}
for a certain object $T \in \Fun(\Gamma_+,\Z)$. Explicitly, $T$ is
given by
\begin{equation}\label{T.eq}
T(S) = \Z[S]/\Z \cdot o, \qquad S \in \Gamma_+,
\end{equation}
where $o \in S$ is the distiguished element. Now consider the functor
$$
j^T:\D(R) \to \D(\Gamma_+,R), \qquad j^T(M) = M \boxtimes T
$$
of \eqref{j.c.t}, with $\C = \ppt$ being the point category, and
$\C_1 = \Gamma_+$. Then the decomposition \eqref{T.1.0} induces an
isomorphism
\begin{equation}\label{j.T.spl.eq}
j_{1!} \cong j_{0!} \oplus j^T,
\end{equation}
so that $j^T$ has an obvious right-adjoint functor
$r^T:\D(\Gamma_+,R) \to \D(R)$ sending $E \in \D(\Gamma_+,R)$ to the
direct summand of $E([1]_+)$ complementary to its direct summand
$E([0]_+)$. Moreover, by \eqref{T.eq}, $T$ satisfies the assumptions
of Lemma~\ref{j.c.le}, so that $j^T$ has a left-adjoint functor
\begin{equation}\label{l.t.eq.g}
l^T:\D(\Gamma_+,R) \to \D(R).
\end{equation}
Note that if we let $T' = j_{1*}\Z \in \Fun(\Gamma_+,\Z)$, then the
isomorphism $T([1]_+) \cong \Z$ extends by adjunction to a map
$\eps:T \to T'$. This induces a functorial map
$$
j^T \to j^{T'} \cong j_{1*},
$$
and by adjunction, we obtain a natural map
\begin{equation}\label{eps.eq}
\eps:j_1^* \to l^T.
\end{equation}

\begin{lemma}\label{gamma.ok}
Proposition~\ref{add.prop} holds for the pair $\langle \Gamma,\Inj
\rangle$. Moreover, the functor $j^T$ induces an equivalence
$$
j^T:\D(R) \cong \DS^{add}_{\Inj}(\Gamma,R) \subset
\DS_{\Inj}(\Gamma,R) \cong \D(\Gamma_+,R),
$$
and the addivization functor
$$
\Add:\D(\Gamma_+,R) \to \D(R) = \DS^{add}_{\Inj}(\Gamma,R)
$$
coincides with the functor $l^T$ of \eqref{l.t.eq.g}.
\end{lemma}

\proof{} One observes immediately that $r^T \circ j^T = \Id$, so
that $j^T$ is a fully faithful embedding, and moreover, an object $E
\in \D(\Gamma_+,R)$ is additive if and only the adjunction map
$j^Tr^T(E) \to E$ is an isomorphism. This proves the first
claim. The second immediately follows from the definition, since
$l^T$ is left-adjoint to $j^T$.
\endproof

Note that more generally, given a small category $\C$, we have a
natural functor
$$
j_\C^T:\D(\C,R) \to \D(\C \times \Gamma_+,R)
$$
of \eqref{j.c.t}, and by Lemma~\ref{j.c.le}, it has a left-adjoint
functor $l^T_\C$. If we denote by $j_1:\C \to \C \times \Gamma_+$
the embedding sending $c \in \C$ to $c \times [1]_+$, then we have a
natural map
\begin{equation}\label{eps.eq.2}
j_1^* \to l^T_\C
\end{equation}
defined by the same procedure as the map \eqref{eps.eq}.

\subsection{The general case.}

Now consider an arbitrary small category $\C$ equip\-ped with an
admissible class of morphisms $I$. Then on one hand, we can consider
the product $\C \times \Gamma$, with the admissible class of maps
$I \times \Inj$, and obtain the corresponding category
$$
\DS_{I \times \Inj}(\C \times \Gamma,R)
$$
for any ring $R$. On the other hand, we can consider the derived
category $\D(S_I\C \times \Gamma_+,R)$. Note that any object $E \in
\D(S_I\C \times \Gamma_+,R)$ gives in particular a functor from
$S_I\C$ to the derived category $\D(\Gamma_+,R)$. Say that $E$ is
special if it sends special maps in $S_I\C$ to invertible maps, and
let
\begin{equation}\label{ds.2.eq}
\DS_I(\C,\Gamma_+,R) \subset \D(S_I\C \times \Gamma_+,R)
\end{equation}
be the full subcategory spanned by special objects. Then the functor
$q$ of \eqref{q.SQ} defines a pullback functor
\begin{equation}\label{q.st.g}
q^*:\DS_I(\C,\Gamma_+,R) \to \DS_{I \times \Inj}(\C \times
\Gamma,R),
\end{equation}
and we observe the following.

\begin{lemma}
The functor $q^*$ of \eqref{q.st.g} is an equivalence of categories.
\end{lemma}

\proof{} Same as Corollary~\ref{S.discr.corr}.
\endproof

We can now define the additivization functor required in
Proposition~\ref{add.prop}. Assume that the pair $\langle \C,I
\rangle$ has finite coproducts in the sense of
Definition~\ref{copr.def}. Then sending an object $c \times S \in \C
\times \Gamma$ to the coproduct of copies of the object $c \in \C$
numbered by elements of the finite set $S$ gives a functor $m:\C
\times \Gamma \to \C$, and the conditions of
Definition~\ref{copr.def} insure that $m$ is a morphism
\begin{equation}\label{m.eq}
m:\langle \C \times \Gamma,I \times \Inj\rangle \to \langle \C,I
\rangle
\end{equation}
in the sense of Definition~\ref{adm.mor.def}. Therefore we have a
pullback functor
$$
S(m)^*:\DS_I(\C,R) \to \DS_{I \times \Inj}(\C \times \Gamma,R)
\cong \DS_I(\C,\Gamma_+,R).
$$
We also have the embedding $e:\C \to \C \times \Gamma$ sending $c \in
\C$ to $c \times [1]$, and $e$ is a morphism
$$
e:\langle \C,I \rangle \to \langle \C \times \Gamma,I \times \Inj
\rangle,
$$
so that we have a pullback functor
$$
S(e)^*:\DS_{I \times \Inj}(\C \times \Gamma,R) \to \DS_I(\C,R).
$$
Since $m \circ e \cong \Id$, the composition $S(e)^* \circ S(m)^*$
is the identity functor.

On the other hand, the functor $j^T_{S_I\C}:\D(S_I\C,R) \to \D(S_I\C
\times \Gamma_+,R)$ obviously sends $\DS_I(\C,R)$ into
$\DS_I(\C,\Gamma_+,R)$, and by Lemma~\ref{j.c.le}, the adjoint
functor $l^T_{S_I\C}$ also sends special objects into special
objects, thus induces a functor
$$
L^T:\DS_I(\C,\Gamma_+,R) \to \DS_I(\C,R).
$$
We define an endofunctor $\Add$ of the category $\DS_I(\C,R)$ by
\begin{equation}\label{add.eq}
\Add = L^T \circ S(m)^*.
\end{equation}
We note that this functor is right-exact with respect to the
standard $t$-structures. Moreover, the embedding $j_1$ of
\eqref{eps.eq.2} coincides with the composition
$$
q \circ S(e):S_I\C \to \C \times \Gamma_+,
$$
so that $j_1^* \circ S(m)^* \cong S(e)^* \circ S(m)^*$ is the
identity functor, and the canonical map \eqref{eps.eq.2} induces a
natural map
\begin{equation}\label{eps.eq.3}
\eps:\Id \cong j_1^* \circ S(m)^* \to \Add.
\end{equation}
Then Proposition~\ref{add.prop} immediately follows from the
following fact.

\begin{lemma}
For any $E \in \DS_I(\C,R)$, $\Add(E) \in \DS_I(\C,R)$ is additive
in the sense of Definition~\ref{add.S.def}, and if $E$ were already
additive, then the natural map $\eps:E \to \Add(E)$ of
\eqref{eps.eq.3} is an isomorphism.
\end{lemma}

\proof{} For any object $c \in \C$, restricting the morphism $m$ of
\eqref{m.eq} to $\{c\} \times \Gamma \subset \C \times \Gamma$ gives
a morphism
$$
m_c:\langle \Gamma,\Inj \rangle \to \langle \C,I \rangle,
$$
and by Lemma~\ref{c.c.le}, $E \in \DS_I(\C,R)$ is additive if and
only if so are all the restrictions $S(m_c)^*E$. By
Lemma~\ref{j.c.le}, the additivization functor $\Add$ commutes with
the restriction functors $S(m_c)^*$. We conclude that both claims of
the lemma can be checked after applying the functors $S(m_c)^*$, $c
\in \C$. In other words, we may assume right away that $\langle \C,I
\rangle \cong \langle \Gamma,\Inj \rangle$. Then we can identify
$\DS_I(\C,R) \cong \D(\Gamma_+,R)$, the product functor $m$ becomes
the smash product functor
$$
m:\Gamma_+ \times \Gamma_+ \to \Gamma_+
$$
sending $[n]_+ \times [n']_+$ to their smash product $[n]_+ \wedge
[n']_+ \cong [nn']_+$ of \eqref{sm.eq}, and by Lemma~\ref{gamma.ok},
it suffices to prove that we have
$$
L^T \circ m^* \cong j^T \circ l^T,
$$
where $L^T:\D(\Gamma_+ \times \Gamma_+,R) \to \D(\Gamma_+,R)$ is the
functor $l^T$ applied along the first variable. By adjunction, this
is equivalent to checking that $m_* \circ j^T \cong j^T \circ r^T$,
or in other words, that
\begin{equation}\label{t.p.m}
m_*(T \boxtimes E) \cong T \otimes r^T(E)
\end{equation}
for any $E \in \Fun(\Gamma_+,R)$, where $r^T$ is the right-adjoint
functor to $j^T$. To check this, it is convenient to use the
equivalence
$$
P:\Fun(\Gamma_+,R) \cong \Fun(\Gamma_-,R),
$$
where $\Gamma_-$ is the category of finite sets and surjective maps
between them (with the same notation as here, this statement is
proved in \cite[Lemma 3.10]{cartier}, and wrongly asserted to belong
to folklore in \cite[Remark 3.11]{cartier} -- in fact it is also
explicitly due to Pirashvili \cite{pira.dk}). The equivalence sends
$T$ to the object $t \in \Fun(\Gamma_-,\Z)$ given by $t([1]) = \Z$
and $t([n]) = 0$ for $n \neq 1$. The functor $r^T$ becomes
evaluation at $[1]$, and $j^T$ becomes the embedding sending $M$ to
$M \otimes t$. And crucially, as shown in \cite[Lemma
  3.18]{cartier}, the functor $m_*$ becomes the restriction to the
diagonal $\Gamma_- \to \Gamma_- \times \Gamma_-$, so that
\eqref{t.p.m} reads as
$$
t \otimes P(E) \cong t \otimes (P(E)([1])).
$$
Since $t([n]) = 0$ for $n \neq 1$, this is obvious.
\endproof

\section{Derived Mackey profunctors.}\label{der.sec}

We can now develop the derived version of the theory of
Section~\ref{mackey.sec}. We assume given a group $G$, and we
consider the category $\Gamma_G$ of finite $G$-sets. It has fibered
products, so that Example~\ref{S.Q.exa} and Example~\ref{S.I.exa}
both apply, and produce two classes of admissible maps in $\Gamma_G$
-- all maps, and injective maps.

\subsection{Pointed sets.}

We first consider the second class -- that is, we equip $\Gamma_G$
with the admissible class $\Inj$ of all injective maps. The category
$Q_{\Inj}\Gamma_G$ of Definition~\ref{Q.def} is then naturally
equivalent to the category $\Gamma_{G+}$ of {\em finite pointed
  $G$-sets} -- that is, finite $G$-sets equipped with a
distinguished $G$-fixed element. The equivalence is the same as in
the case of the category of finite sets: it sends a pointed $G$-set
$S$ with the distinguished element $o \in S$ to the complement
$\overline{S} = S \setminus \{o\}$, and a pointed map $f:S \to S'$
goes to the diagram \eqref{partial.map}. We have a natural embedding
$\Gamma_G \to \Gamma_{G+}$ sending a finite $G$-set $S$ to the set
$S_+$ of \eqref{s.pl}, with $G$ acting triviality on the added
distiguished point $o$. The category $\Gamma_{G+}$ has finite
coproducts given by \eqref{vee.eq}, and every $G$-set $S$ decomposes
as
\begin{equation}\label{s.vee}
S = [G/H_1]_+ \vee \dots \vee [G/H_n]_+ = \left( \coprod_{1 \leq i
    \leq n}[G/H_i]\right)_+,
\end{equation}
where $H_i \subset G$, $1 \leq i \leq n$ are cofinite subgroups, and
$[G/H_i]$ are the corresponding $G$-orbits. For any $S_1 \in
\Gamma_{G+}$, we have a natural identification
\begin{equation}\label{s.prod}
\Gamma_{G+}(S,S_1) = \prod_i\Gamma_{G+}([G/H_i]_+,S_1) =
\prod_iS_1^{H_i},
\end{equation}
where $S_1^{H_i} \subset S_1$ is the subset of $H_i$-fixed
points.The pair $\langle \Gamma_G,\Inj \rangle$ is discrete in the
sense of Corollary~\ref{S.discr.corr}, so that we have a canonical
equivalence
\begin{equation}\label{gamma.pl}
\DS_{\Inj}(\Gamma_G,R) \cong \D(\Gamma_{G+},R)
\end{equation}
for any ring $R$. 

Analogously, let $\wGamma_G$ be the category of $G$-sets admissible
in the sense of Definition~\ref{adm.G.set}, with the admissible
class $\Inj$ of all monomorphisms. Then the same equivalence
identifies $Q_{\Inj}\wGamma_G$ with the category $\wGamma_{G+}$ of
admissible pointed $G$-sets. The category $\wGamma_{G+}$ also has
coproducts, and any $S \in \wGamma_{G+}$ has a decomposition
\eqref{s.vee} although possibly with an infinite number of terms. We
have a natural embedding $\wGamma_G \subset \wGamma_{G+}$, $S
\mapsto S_+$. The pair $\langle \wGamma_G,\Inj \rangle$ is discrete
in the sense of Corollary~\ref{S.discr.corr}, so that we have an
equivalence
\begin{equation}\label{wgamma.pl}
\DS_{\Inj}(\wGamma_G,R) \cong \D(\wGamma_{G+},R)
\end{equation}
for any ring $R$.

Now fix a cofinite subgroup $H \subset G$, with the normalizer $N_H
\subset G$ and the quotient $W = N_H/H$. Then the fixed points
functor $\phi^H$ of \eqref{phi.H} tautologically extends to a
functor
$$
\phi^H:\wGamma_{G+} \to \Gamma_{W+}
$$
sending $S$ to $S^N$. This functor preserves smash products, that
is, we have
\begin{equation}\label{wedge.H}
(S_1 \wedge S_2)^H = S_1^H \wedge S_2^H.
\end{equation}
for any $S_2,S_2 \in \wGamma_{G+}$. Moreover, assume that $H=N$ is
normal, so that $N_H = G$ and $W$ is the quotient $G/N$. Then
$\phi^N$ has a left-adjoint $\lambda^N:\Gamma_{W+} \to \wGamma_{G+}$
sending a $W$-set $S$ to the same same on which $G$ acts via the
quotient map $G \to W$. Thus the pullback functor $\lambda^{N*}$ is
right-adjoint to the pullback functor $\phi^{N*}$, and since $\phi^N
\circ \lambda^N$ is obviously the identity functor, $\phi^{N*}$ is a
full embedding.

We will also need an explicit description of the left-adjoint functor
\begin{equation}\label{phi.n.pl}
L^\hdot\phi^N_!:\D(\wGamma_{G+},R) \to \D(\Gamma_{W+},R).
\end{equation}
To obtain it, we use simplicial combinatorics of
Subsection~\ref{simpl.subs} and Subsection~\ref{contr.subs}. So, we
consider $n$-simplicial pointed admissible $G$-sets -- that is,
functors $X:(\Delta^o)^n \to \wGamma_{G+}$ from the $n$-fold
self-product $(\Delta^o)^n$ to the category $\wGamma_{G+}$.

\begin{defn}\label{homo.point.def}
  The {\em homology complex $C_\idot(X,E)$} of an $n$-simplicial
  pointed admissible $G$-set $X$ with coefficients in an
  object $E \in \D(\wGamma_{G+},R)$ is the complex
$$
C_\idot(X,E) = C_\idot((\Delta^o)^n,X^*E),
$$
where $X^*:\D(\wGamma_{G+},R) \to \D((\Delta^o)^n,R)$ is the pullback
functor associated to $X:(\Delta^o)^n \to \Gamma_{G+}$.
\end{defn}

\begin{lemma}\label{contr.ex.le}
Assume given a pointed $n$-simplicial admissible $G$-set $X$ such
that for any cofinite subgroup $H \subset G$, the fixed point set
$X^H$ is contractible in the sense of
Definition~\ref{contr.def}. Then for any $E \in \D(\wGamma_{G+},R)$,
the map
$$
E([0]_+) \to C_\idot(X,E)
$$
induced by the distinguished point embedding $[0]_+ \to X$ is a
quasiisomorphism.
\end{lemma}

\proof{} Instead of an arbitrary $E$, it suffices to consider the
generators $M_S$ of the category $\D(\wGamma_{G+},R)$ given by
\eqref{m.c.eq}, $M$ an $R$-module, $S \in \wGamma_{G+}$ an
admissible pointed $G$-set. Fix $M$ and $S$, and consider the
$n$-simplicial set $\wGamma_{G+}(S,X)$ of\eqref{C.cc}. Then by
\eqref{m.c.eq}, we have
$$
C_\idot(X,M_S) \cong C_\idot(\wGamma_{G+}(S,X),\Z) \otimes M,
$$
the homology of the $n$-simplicial set $\wGamma_{G+}(S,X)$ with
coefficients in $M$. But $S$ has the decomposition \eqref{s.vee},
and \eqref{s.prod} gives an isomorphism
$$
\wGamma_{G+}(S,X) \cong \prod_i X^{H_i}.
$$
Since all the terms $X^{H_i}$ in the product are by assumption
contractible, the product itself is contractible, and we are done by
Lemma~\ref{contr.le}.
\endproof

\begin{defn}\label{ada.N.def}
  For any normal cofinite subgroup $N \subset G$, an $n$-simplicial
  pointed $G$-set $X \in (\Delta^o)^n\Gamma_{G+}$ is {\em
    $N$-adapted} if
\begin{enumerate}
\item the fixed points subset $X^N$ is isomorphic to the constant
  $n$-simplicial set pointed $[1]_+$ with the trivial $G$-action, and
\item for any cofinite subgroup $H \subset G$ not containing $N$,
  the fixed points subset $X^H \subset X$ is contractible in the
  sense of Definition~\ref{contr.def}.
\end{enumerate}
\end{defn}

\begin{lemma}\label{ada.ex.le}
  For any cofinite normal subgroup $N \subset G$, there exists an
  admissible $2$-simplicial pointed $G$-set $S$ that is $N$-adapted
  in the sense of Definition~\ref{ada.N.def}.
\end{lemma}

\proof{} Let
$$
S = \coprod_{N \subset G}[G/N]
$$
be the disjoint union of all $G$-orbits $[G/N]$, $N \subset G$ a
cofinite normal subgroup. Then $S$ is an admissible $G$-set. Fix a
cofinite normal subgroup $N \subset G$, and let
$$
S_N = S \setminus S^N.
$$
Then $(S_N)^N$ is empty, and $(S_N)^H$ is not empty for any cofinite
subgroup $H \subset G$ not containing $N$. Therefore the $2$-simplicial
set $C(ES)$ obtained by combining Example~\ref{E.exa} and
Example~\ref{cone.exa} satisfies all the assumptions.
\endproof

\begin{lemma}\label{ada.le}
  For any $n$-simplicial pointed admissible $G$-set $X$ adapted to
  $N$, any $S \in \wGamma_{G+}$, and any $E \in \D(\wGamma_{G+},R)$,
  the natural map
\begin{equation}\label{s.n.x}
C_\idot(S^N \wedge X,E) \to C_\idot(S \wedge X,E)
\end{equation}
induced by the embedding $S^N \to S$ is a quasiisomorphism.
\end{lemma}

\proof{} As in the proof of Lemma~\ref{contr.ex.le}, we may assume
that $E=M_{S'}$, $M$ an $R$-module, $S' \in \Gamma_{G+}$ an
admissible pointed $G$-set. Fix $M$ and $S'$, and for any
$n$-simplicial pointed admissible $G$-set $Y$, consider the
$n$-simplicial set $\wGamma_{G+}(S',Y)$ of\eqref{C.cc}. Then by
\eqref{m.c.eq}, we have
\begin{equation}\label{c.prod}
C_\idot(Y,M_{S'}) \cong C_\idot(\wGamma_{G+}(S',Y),\Z) \otimes M.
\end{equation}
The pointed set $S'$ has the decomposition \eqref{s.vee}, and
\eqref{s.prod} gives an isomorphism
$$
\wGamma_{G+}(S',Y) \cong \prod_i Y^{H_i}.
$$
Since $S$ is admissible, all but a finite number $H_1,\dots,H_n$ of
subgroups $H_i \subset G$ do not contain $N \subset G$. Denote by
$\overline{Y}$ the product of the terms corresponding to these
subgroups, so that we have
\begin{equation}\label{s.pro.d}
\wGamma_{G+}(S',Y) \cong \overline{Y} \times \prod_{i=1}^n Y^{H_i}.
\end{equation}
Then by the K\"unneth formula, \eqref{c.prod} gives a
quasiisomorphism
\begin{equation}\label{times.c}
C_\idot(Y,M_S) \cong M \otimes C_\idot(\overline{Y},\Z) \otimes
\bigotimes_{i=1}^n C_\idot(Y^{H_i},\Z).
\end{equation}
Now, by Definition~\ref{ada.N.def}~\thetag{ii}, for every cofinite
$H \subset G$ not containing $N$, the fixed point set $X^H$ is
contractible. Thus if we take $Y = S \wedge X$, then for any such
$H$,
$$
Y^H = (S \wedge X)^H = S^H \wedge X^H
$$
is also contractible, and then $\overline{Y}$, being the product of
contractible pointed $n$-simplicial sets, is itself
contractible. Then by Lemma~\ref{contr.ex.le}, \eqref{times.c}
reduces to a quasiisomorphism
$$
C_\idot(S \wedge X,M_S) \cong M \otimes \bigotimes_{i=1}^n C_\idot(S^{H_i}
\wedge X^{H_i},\Z),
$$
where the product is over those subgroups $H_i \subset G$ in
\eqref{s.vee} that contain $N$. But for such subgroups, $S^{H_i}
\subset S$ lies in $S^N \subset S$. Therefore if we do the same
reduction for $Y = S^N \wedge X$, then the result is exactly the
same.
\endproof

Now fix an $N$-adapted $2$-simplicial admissible pointed $G$-set $X$
-- for example, the one provided by Lemma~\ref{ada.ex.le} -- and
consider the functor
$$
m_X:(\Delta^o)^2 \times \wGamma_{G+} \to \wGamma_{G+}
$$
given by
$$
m_X([n_1] \times [n_2] \times S) = X([n_1] \times [n_2]) \wedge S.
$$
Let $p:(\Delta^o)^2 \times \wGamma_{G+} \to \wGamma_{G+}$ be the
natural projection, and for any $E \in \D(\wGamma_{G+},R)$, let
$$
\Av_X(E) = L^\hdot p_!m_X^*E.
$$
Then $\Av_X$ is an endofunctor of the category $\D(\wGamma_{G+},R)$,
a sort of ``averaging'' over the pointed $2$-simplicial $G$-set
$X$. By base change, for any $S \in \wGamma_{G+}$ and any $E \in
\D(\wGamma_{G+},R)$, we have
$$
\Av_X(E)(S) \cong C_\idot(S \wedge X,E).
$$
The natural embedding $[1]_+ \cong X^N \to X$ induces a functorial
map $S \to S \wedge X$, so that we have a functorial map
\begin{equation}\label{id.ad}
\Id \to \Av_X.
\end{equation}
Moreover, define a functor $\Phi^N:\D(\wGamma_{G+},R) \to
\D(\Gamma_{W+},R)$ by
$$
\Phi^N = \lambda^{N*} \circ \Av_X.
$$
Then for any $S \in \wGamma_{W+}$ and any $E \in
\D(\wGamma_{G+},R)$, we have
\begin{equation}\label{Phi.N.exp}
\Phi^N(E)(S) = C_\idot(\lambda^N(S) \wedge X,E),
\end{equation}
and for any $S \in \wGamma_{G+}$ and any $E \in \D(\wGamma_{G+},R)$,
we have
$$
\phi^{N*}\Phi^N(E)(S) = C_\idot(S^N \wedge X,E).
$$
The map $S^N \to S$ induces a map $\phi^{N*}\Phi^N \to \Av_X$, and
by Lemma~\ref{ada.le}, this map is a quasiisomorphism. Thus
\eqref{id.ad} induces by adjunction a map
\begin{equation}\label{phi.phi}
L^\hdot\phi^N_! \to \Phi^N.
\end{equation}

\begin{prop}\label{ada.prop}
The map \eqref{phi.phi} is an isomorphism of functors.
\end{prop}

\proof{} As in Lemma~\ref{contr.ex.le}, we may assume $E = M_S$, $M$ an
$R$-module, $S \in  \wGamma_{G+}$ an admissible $G$-set with
decomposition \eqref{s.vee}. Then by adjunction,
$$
L^\hdot\phi^N_!M_S = M_{\phi^N(S)} = M_{S^N},
$$
so that for any $\wt{S} \in \Gamma_{W+}$, we have
$$
L^\hdot\phi^N_!M_S(\wt{S}) = M \otimes
\Z\left[\Gamma_{W+}(S^N,\wt{S})\right] =
M \otimes \bigotimes_i\Z[\wt{S}^{H_i}],
$$
where the product is taken over all subgroups $H_i \subset G$ in
\eqref{s.vee} that contain $N$. On the other hand, the value
$\Phi^N(M_S)(\wt{S})$ of the functor $\Phi^N$ can be computed by
\eqref{Phi.N.exp}. Then as in the proof of Lemma~\ref{ada.le}, we
have
$$
\Phi^N(M_S)(\wt{S}) \cong \M \otimes \bigotimes_i
C_\idot(\wt{S}^{H_i} \wedge X^{H_i},\Z),
$$
where again, the product is taken over all subgroups $H_i \subset G$
in \eqref{s.vee} that contain $N$. To finish the proof, it remains
to notice that for each of these subgroups, we have
$$
\wt{S}^{H_i} \cong \wt{S}^{H_i} \wedge [1]_+ \cong \wt{S}^{H_i} \wedge X^{H_i}
$$
by Definition~\ref{ada.N.def}~\thetag{i}.
\endproof

\subsection{Profunctors.}

Now equip the category $\Gamma_G$ with the admissible class of all
maps, fix a ring $R$, and consider the category $\DS(\Gamma_G,R)$ of
Definition~\ref{DS.def}. We have an embedding $\Gamma_G^o \cong
(S\Gamma_G)_{[1]} \subset S\Gamma_G$ and the corresponding
restriction functor $\DS(\Gamma_G,R) \to \D(\Gamma_G^o,R)$. Note
that every object $E$ in $\D(\Gamma_G^o,R)$ defines in particular a
functor from $\Gamma_G^o$ to the additive category $\D(R)$.

\begin{defn}
An {\em $R$-valued derived $G$-Mackey functor} is an object $E \in
\DS(\Gamma_G,R)$ such that the corresponding functor $\Gamma_G^o \to
\D(R)$ is additive in the sense of Definition~\ref{add.def}. The
full subcategory in $\DS(\Gamma_G,R)$ spanned by derived Mackey
functors is denoted by $\DM(G,R) \subset \DS(\Gamma_G,R)$.
\end{defn}

We note that the pair $\langle \Gamma_G,\Id \rangle$ has finite
coproducts in the sense of Definition~\ref{copr.def}, and the
additivity conditions of Definition~\ref{add.def} and
Definition~\ref{add.S.def} are identically the same. Therefore we
have
$$
\DM(G,R) = \DS^{add}(\Gamma_G,R) \subset \DS(\Gamma_G,R),
$$
and Proposition~\ref{add.prop} insures that the embedding $\DM(G,R)
\subset \DS(\Gamma_G,R)$ admits a left-adjoint additivization
functor
\begin{equation}\label{add.der}
\Add:\DS(\Gamma_G,R) \to \DM(G,R).
\end{equation}
The triangulated subcategory $\DM(G,R) \subset \DS(\Gamma_G,R)$
inherits a natural $t$-structure, and the additivization functor
$\Add$ is right-exact with respect the natural $t$-structures. The
induced functor on the hearts $t$-structures is the additivization
functor \eqref{Add}.

Analogously, consider the category $\wGamma_G$ of $G$-sets
admissible in the sense of Definition~\ref{adm.G.set}. It also has
fibered products, so that for every ring $R$, we have the category
$\DS(\wGamma_G,R)$ and the restriction functor
$$
\DS(\wGamma_G,R) \to \D(\wGamma_G^o,R).
$$
Every object $E \in \D(\wGamma_G^o,R)$ in turns gives a functor
$\wGamma_G^o \to \D(R)$.

\begin{defn}
  An {\em $R$-valued derived $G$-Mackey profunctor} is an object $E
  \in \DS(\wGamma_G,R)$ such that the corresponding functor
  $\wGamma_G^o \to \D(R)$ is additive in the sense of
  Definition~\ref{add.pro}. The subcategory of derived Mackey
  profunctors is denoted by $\wDM(G,R) \subset \DS(\wGamma_G,R)$.
\end{defn}

The additivity condition of Definition~\ref{add.pro} is preserved
under truncation with respect to the natural $t$-structure on
$\DS(\wGamma_G,R)$, so that $\wDM(G,R) \subset \DS(\wGamma_G,R)$
inherits a natural $t$-structure. We denote by
$$
\wDM^-(G,R) = \bigcup_i \wDM^{\leq i}(G,R)
$$
the full subcategory of objects bounded from above with respect to
this $t$-structure. The pullback functor $q^*$ of \eqref{q.st}
commutes with restrictions to $\wGamma_G^o$, so that we have a
natural functor
$$
q^*:\D(\wM(G,R)) \to \wDM(G,R),
$$
where $\D(\wM(G,R))$ is the derived category of the abelian category
$\wM(G,R)$ of $R$-valued $G$-Mackey profunctors of
Definition~\ref{add.pro}. By Lemma~\ref{M.DM}, this induces a fully
faithful embedding
\begin{equation}\label{M.DM.eq}
\wM(G,R) \subset \wDM^{\leq 0}(G,R) \subset \wDM(G,R)
\end{equation}
that identifies the category $\wM(G,R)$ with the heart of the
natural $t$-structure on $\wDM(G,R)$. We denote by
\begin{equation}\label{tau.eq}
\tau:\wDM^{\leq 0}(G,R) \to \wM(G,R)
\end{equation}
the associated truncation functor.

As in the non-derived case, if $G$ is finite, all admissible
$G$-sets are also finite, so that $\Gamma_G = \wGamma_G$ and
$\wDM(G,R) = \DM(G,R)$. In the general case, we still have an
embedding $\Gamma_G \to \wGamma_G$ and the restriction functor
$\DS(\wGamma_G,R) \to \DS(\Gamma_G,R)$. This clearly sends objects
additive in the sense of Definition~\ref{add.pro} to objects
additive in the sense of Definition~\ref{add.def}, so that we have a
natural forgetful functor
\begin{equation}\label{d.M.wM.eq}
\wDM(G,R) \to \DM(G,R).
\end{equation}
This is compatible with the natural $t$-structures and induces the
restriction functor \eqref{M.wM.eq} on their hearts.

The identity functor $\Gamma_G \to \Gamma_G$ defines a morphism
\begin{equation}\label{i.i.eq}
i:\langle \Gamma_G,\Inj \rangle \to \langle \Gamma_G,\Id \rangle
\end{equation}
in the sense of Definition~\ref{adm.mor.def}, so that we have a
restriction functor
\begin{equation}\label{s.i.eq}
S(i)^*:\DM(G,R) \to \DS_{\Inj}(\Gamma_G,R) \cong \D(\Gamma_{G+},R).
\end{equation}
Analogously, we have a morphism $i:\langle \wGamma_G,\Inj \rangle \to
\langle \wGamma_G,\Id \rangle$ and a restriction functor
\begin{equation}\label{s.i.w.eq}
S(i)^*:\wDM(G,R) \to \DS_{\Inj}(\wGamma_G,R) \cong \D(\wGamma_{G+},R).
\end{equation}

\begin{defn}\label{homo.DM.def}
  The {\em homology complex $C_\idot(X,E)$} of an $n$-simplicial
  pointed admissible $G$-set $X$ with coefficients in a derived
  Mackey profunctor $E \in \wDM(G,R)$ is the complex
$$
C_\idot(X,E) = C_\idot(X,S(i)^*E),
$$
where $S(i)^*$ is the restriction functor \eqref{s.i.w.eq}, and
$C_\idot(X,S(i)^*E)$ is the homology complex of
Definition~\ref{homo.point.def}.
\end{defn}

While the homology complex of Definition~\ref{homo.point.def} is
obviously covariant with respect to $X$, the homology complex of
Definition~\ref{homo.DM.def} is also contravariant: we have a
natural map
\begin{equation}\label{homo.f}
f^*:C_\idot(X',E) \to C_\idot(X,E)
\end{equation}
for any $E \in \wDM(G,R)$ and any map $f:X \to X'$ of $n$-simplicial
pointed admissible $G$-sets. It also enjoys the following continuity
properties that we will need later.

\begin{lemma}\label{prod.le}
\begin{enumerate}
\item Assume given an $n$-simplicial pointed admissible $G$-set $X$
  and an inverse system $\{E_i\}$ of derived Mackey profunctors $E_i
  \in \wDM^{\leq 0}(G,R)$. Then we have a natural identification
$$
C_\idot(X,\Tel(E_i)) \cong \Tel(C_\idot(X,E_i)),
$$
where $\Tel$ is the telescope of \eqref{Tel}.
\item Assume given a derived Mackey profunctor $E \in \wDM^{\leq
  0}(G,R)$ and an inverse system of subsets $X_{i+1} \subset X_i
  \subset X$ of an $n$-simplicial pointed admissible
  $G$-set $X$ with intersection $\overline{X} = \cap_iX_i$. Then we
  have a natural identification
$$
C_\idot(\overline{X},E) \cong \Tel(C_\idot(X_i,E)).
$$
\end{enumerate}
\end{lemma}

\proof{} For \thetag{i}, note that we have $\Tel(X^*E_i) \cong
X^*\Tel(E_i)$, and then note that for an inverse system $M_i$ of
objects in $\D^{\leq n}(\Delta^o,R)$, with some fixed integer $n$
indepedent of $i$, we have
$$
C_\idot(\Delta^o,\Tel(M_i)) \cong \Tel(C_\idot(\Delta^o,M_i)).
$$
For \thetag{ii}, observe that we have $\Tel(X_i^*E) \cong
\overline{X}^*E$ by \eqref{add.pro.eq}, and use the same argument.
\endproof

Assume now given a subgroup $H \subset G$, with normalizer $N_H
\subset G$ and the quotient $W = N_H/H$. Then the functor
$\phi^H$ of \eqref{phi.H} defines a morphism
$$
\phi^H:\langle \wGamma_G,\Id \rangle \to \langle \wGamma_W,\Id
\rangle
$$
in the sense of Definition~\ref{adm.mor.def}, so that we have a
pullback functor
\begin{equation}\label{S.phi.H.st}
S(\phi^H)^*:\DS(\wGamma_W,R) \to \DS(\wGamma_G,R).
\end{equation}
This functor clearly preserves additivity in the sense of
Definition~\ref{add.pro}, thus sends a derived Mackey profunctor to
a derived Mackey profunctor. We also have the left-adjoint functor
\begin{equation}\label{S.phi.H.sh}
S(\phi^H)_!:\DS(\wGamma_G,R) \to \DS(\wGamma_W,R)
\end{equation}
provided by Corollary~\ref{sp.corr}.

\begin{defn}\label{infl.def}
The {\em inflation functor}
$$
\Infl^H:\wDM(W,R) \to \wDM(G,R)
$$
is the functor induced by the pullback functor \eqref{S.phi.H.st}.
\end{defn}

For any normal subgroup $N \subset G$ with the quotient $W = G/N$,
say that an object $E \in \DS(\wGamma_G,R)$ is {\em supported at $N$} if
for any admissible $G$-set $S$, the natural map $E(S) \to E(S^N)$ is
a quasiisomorphism. Then by definition, $\Infl^N$ factors
through the full subcategory
\begin{equation}\label{dm.h.eq}
\wDM_N(G,R) \subset \wDM(G,R)
\end{equation}
spanned by objects $E \in \wDM(G,R)$ supported at $N$. Note that by
\eqref{add.pro.eq}, a derived Mackey profunctor $E \in \wDM(G,R)$ is
supported at $N$ if and only if $M([G/H]) = 0$ for any cofinite $H
\subset G$ not contained in $N \subset G$.

\begin{lemma}\label{infl.der}
For any normal subgroup $N \subset G$, the inflation functor induces
an equivalence
$$
\Infl^N:\wDM(W,R) \cong \wDM_N(G,R).
$$
\end{lemma}

\proof{} More generally, we will prove that $S(\phi^N)^*$ induces an
equivalence between $\DS(\wGamma_W,R)$ and the full subcategory in
$\DS(\wGamma_G,R)$ spanned by objects supported at $N$. This is
equivalent to proving that the adjunction map
$$
S(\phi^N)_!S(\phi^N)*E \to E
$$
is an isomorphism for any $E \in \DS(\wGamma_W,R)$, and the
adjunction map
$$
E \to S(\phi^N)^*S(\phi^N)_!E
$$
is an isomorphism for any $E \in \DS(\wGamma_G,R)$ supported at
$N$. But by Proposition~\ref{bc.prop}, we have
\begin{equation}\label{bc.S.eq}
S(i)^* \circ S(\phi^N)_! \cong L^\hdot\phi^N_! \circ S(i)^*,
\end{equation}
where $S(i)^*$ are the canonical restriction functors
\eqref{s.i.eq}, \eqref{s.i.w.eq}, and $L^\hdot\phi^N_!$ is the
functor \eqref{phi.n.pl}. Since we obviously also have $S(i)^* \circ
S(\phi^N)^* \cong \phi^{N*} \circ S(i)^*$, it suffices to prove both
claims after restriction to $\wGamma_{G+}$. But then, it is
equivalent to proving that $\phi^{N*}$ induces an equivalence between
$\D(\wGamma_{W+},R)$ and the full subcategory in
$\D(\wGamma_{G+},R)$ spanned by objects supported at $N$.  Since
$\phi^N$ has a fully faithful left adjoint $\lambda^N$, this is
obvious.
\endproof

Now assume that the subgroup $H \subset G$ is cofinite. Then we have
the adjoint pair of functors $\rho^H$, $\gamma^H$ of
\eqref{rho.gamma.eq}, and we note that this is in fact a situation
of Example~\ref{c.adj} with $\C = \wGamma_G$, $c = [G/H]$, $c' =
[G/G]$, and $f:c \to c'$ the natural projection. Indeed, for any a
$G$-set $S$ equipped with a map $f:S \to [G/H]$, the preimage
$f^{-1}(e) \subset S$ of the image $e \in [G/H]$ of the unity
element $e \in G$ is an $H$-set, and sending $f:S \to [G/H]$ to
$f^{-1}(e)$ gives an equivalence of categories
$$
\wGamma_G/[G/H] \cong \wGamma_H.
$$
Under this equivalence, $f^*$ of Example~\ref{c.adj} becomes the
restriction functor $\rho^H$, and $f_!$ becomes its left-adjoint
functor $\gamma^H$. Then \eqref{phi.psi.eq} becomes the isomorphism
\begin{equation}\label{gamma.rho.H}
S(\gamma^H)^* \cong S(\rho^H)_!,
\end{equation}
and this functor preserves additivity in the sense of
Definition~\ref{add.pro}.

\begin{defn}\label{cat.fp.def}
The {\em categorical fixed points functor}
$$
\Psi^H:\wDM(G,R) \to \wDM(H,R)
$$
is the functor induced by the functor \eqref{gamma.rho.H}.
\end{defn}

As in the underived situation, the centralizer $Z_H \subset G$ of
any cofinite group $H \subset G$ acts on the functor $\rho^H$ and
its adjoint $\gamma^H$, thus on $\Psi^H = S(\gamma^H)^*$, so that it
can be promoted to a functor
\begin{equation}\label{psi.wt.der}
\wt{\Psi}^H:\M(G,R) \to \M(H,R[Z_H]),
\end{equation}
where $R[Z_H]$ is the group algebra of the centralizer.

Now we note that for any cofinite subgroup $H \subset G$ with
normalizer $N_H \subset G$ and quotient $W_H = N_H/N$, we have
$\phi^H \cong \phi^H \circ \rho^{N_H}$, so that
\begin{equation}\label{W.G.eq}
S(\phi^H)_! \cong S(\phi^H)_! \circ S(\rho^{N_H})_! \cong
S(\phi^H)_! \circ S(\gamma^{N_H})^*.
\end{equation}

\begin{lemma}\label{phi.add}
  For any cofinite normal subgroup $N \subset G$ with the quotient
  $W = G/N$, the functor $S(\phi^N)_!$ of \eqref{S.phi.H.sh} sends
$\wDM(G,R) \subset \DS(\wGamma_G,R)$ to $\DM(W,R) \subset
\DS(\Gamma_W,R)$.
\end{lemma}

\proof{} By \eqref{W.G.eq}, we may replace $G$ with $N_H$ -- in
other words, we may assume right away that $H = N \subset G$ is
normal, and $W = G/N$. Then as in the proof of Lemma~\ref{infl.der},
Proposition~\ref{bc.prop} provides a canonical isomorphism
\eqref{bc.S.eq}, and by Proposition~\ref{ada.prop}, for any $E \in
\DS(\wGamma_G,R)$, the values $S(\phi^N)_!(E)(S)$, $S \in \Gamma_W$
can be computed by \eqref{Phi.N.exp}, where the right-hand side is
the homology complex of Definition~\ref{homo.DM.def}. Then the claim
becomes obvious -- we have
$$
(S \copr S')_+ \wedge X \cong (S_+ \wedge X) \vee (S'_+ \wedge X)
$$
for any two $G$-sets $S_1$, $S_2$, .
\endproof

\begin{defn}\label{geo.fp.def}
The {\em geometric fixed points functor}
$$
\Phi^H:\wDM(G,R) \to \DM(W,R)
$$
is the functor induced by the functor \eqref{S.phi.H.sh}.
\end{defn}

By definition, the functor $\Phi^H$ is left-adjoint to the inflation
functor $\Infl^H$ of Definition~\ref{infl.def}. For any subgroup $H'
\subset H \subset G$, the isomorphism \eqref{Phi.Psi.iso} induces an
isomorphism
\begin{equation}\label{phi.psi.iso.der}
\Phi^{H'} \circ \Psi^H \cong \Psi^{W^H_{H'}} \circ \Phi^{H'},
\end{equation}
a derived version of the isomorphism \eqref{phi.psi.iso}.

\begin{lemma}\label{der.noder.le}
The functors $\Infl^H$ of Definition~\ref{infl.def}, $\Psi^H$ of
Definition~\ref{cat.fp.def} are exact with respect to the standard
$t$-structures, and induce the functors $\Infl^H$, $\Psi^H$ of
\eqref{phi.infl.eq} resp. \eqref{psi.pro.eq} on their hearts. The
functor $\Phi^H$ of Definition~\ref{geo.fp.def} is right-exact with
respect to the standard $t$-structures, and induces the functor
$\Phi^H$ of \eqref{phi.infl.eq} on their hearts (that is, we have
$\tau \circ \Phi^H \cong \Phi^H \circ \tau$, where $\tau$ is the
truncation functor \eqref{tau.eq}).
\end{lemma}

\proof{} The first claim is obvious, the second immediately follows
by adjunction.
\endproof

\section{Mackey functors and representations.}\label{repr.sec}

If the group $G$ is finite, then any admissible $G$-set $S$ is
finite, so that derived Mackey profunctors are the same as derived
Mackey functors studied in \cite{mackey}. However, most of the
proofs in \cite{mackey} are rather involved. Throughout this
section, we assume that $G$ is a finite group, and we reprove some
of the results from \cite{mackey} using the technology we have
developed, especially Proposition~\ref{ada.prop}.

\subsection{Fixed points as representations.}

Fix a finite group $G$. For any subgroup $H \subset G$, denote
\begin{equation}\label{bphi.H}
\bPhi^H = \Phi^H \circ \Psi^H:\DM(G,R) \to \D(R),
\end{equation}
and let
\begin{equation}\label{bphi.dot}
\bPhi_\idot = \bigoplus_{H \subset G}\bPhi^H:\DM(G,R) \to \prod_{H
  \subset G}\D(R)
\end{equation}
be the sum of the functors $\bPhi^H$ over all conjugacy classes of
subgroups $H \subset G$.

\begin{lemma}\label{phi.dot.le}
\begin{enumerate}
\item If the image $\bPhi_\idot(M)$ of a derived Mackey functor $M
  \in \DM(G,R)$ under the functor \eqref{bphi.dot} lies in
  $\prod_H\D^{\leq n}(R)$ for some integer $n$, then $M$ lies in
  $\DM^{\leq n}(G,R) \subset \DM(G,R)$.
\item The functor $\bPhi_\idot$ of \eqref{bphi.dot} is conservative
  (that is, if $\Phi_\idot(f)$ is invertible, then $f$ is
  invertible).
\end{enumerate}
\end{lemma}

\proof{} For \thetag{i}, let $M'=\tau_{> n}M$ be the truncation of
$M$ with respect to the standard $t$-structure on $\DM(G,R)$. We
need to prove that $M' = 0$. By \eqref{phi.psi.iso.der} and
induction on cardinality of $G$, we may assume that $\Psi^HM'=0$ for
any proper subgroup $H \subset G$. Then $M'$ is supported at $G$,
and by Lemma~\ref{infl.der}, $M' \cong \Infl^G(\Phi^G(M'))$ and
$\Phi^G(M') \cong M'([G/G])$. In particular, $\Phi^G(M')$ is
non-trivial only in homological degrees $> n$. Since the functor
$\Phi^G$ is left-exact with respect to the standard $t$-structure,
the natural map $\Phi^GM \to \Phi^GM'$ is an isomorphism in these
homological degrees, and since by assumption $\Phi^GM = \bPhi^GM$
lies in $\D^{\leq n}(R)$, we have $\Phi^GM' = 0$ and $M' \cong
\Infl^G(\Phi^G(M')) = 0$.

For \thetag{ii}, let $M$ be the cone of the map $f$. Then
$\bPhi^\hdot(M) = 0$, so that $M$ must lie in $\DM^{\leq n}(G,R)$
for any integer $n$. Therefore $M=0$, and $f$ is invertible.
\endproof

By \eqref{phi.psi.iso.der}, we have $\bPhi^H \cong \Psi^{\{e\}}
\circ \Phi^H$, where $\{e\} \subset W_H$ is the trivial subgroup in
the quotient $W_H=N_H/H$, $N_H \subset G$ the normalizer of $H$. We
can promote $\Psi^{\{e\}}$ to the functor $\wt{\Psi}^{\{e\}}$ of
\eqref{psi.wt.der}, and this promotes $\bPhi^H$ to a functor
\begin{equation}\label{wphi.H}
\wPhi^H = \wt{\Psi}^{\{e\}} \circ \Phi^H:\DM(G,R) \to \D(R[W_H]).
\end{equation}

\begin{lemma}\label{R.H.le}
For any $H \subset G$, the functor $\wPhi^H$ of \eqref{wphi.H}
admits a fully faithful right-adjoint functor
$$
R_H:\D(R[H]) \to \DM(H,R).
$$
\end{lemma}

\proof{} The functor $\bPhi^H = S(\phi^H \circ \rho^H)_!$ has an
obvious right-adjoint $\overline{R}_H = S(\phi^H \circ
\rho^H)^*$. Explicitly, if we consider the object $T_H \in \M(G,Z)$
given by
$$
T_H(S) = \Z[S^H], \qquad S \in \Gamma_G,
$$
then we have $\overline{R}_H(M) = M \otimes T$ for any $M \in \D(R)$. To promote
it to a right-adjoint functor $R_H$, observe that $W=W_H$ acts on $T$,
so that it can be considered as an object $T \in \D(\ppt_W \times
S\Gamma_G,\Z)$, where $\ppt_W$ is the groupoid with one object with
automorphism group $W$. Then $R_H$ is given by
\begin{equation}\label{r.g.eq}
R_H(M) = R^\hdot\pi_{2*}(\pi_1^*M \otimes T),
\end{equation}
where $\pi_1:\ppt_W \times S\Gamma_G \to \ppt_W$, $\pi_2:\ppt_W
\times S\Gamma_G \to S\Gamma_G$ are the natural
projections. Explicitly, for any $S \in \Gamma_G$, we have
\begin{equation}\label{r.exp}
R_H(M)(S) \cong C^\hdot(W,M \otimes T(S)) = C^\hdot(W,M[S^H]).
\end{equation}
This implies that for any $S \in \Gamma_H$, we have
\begin{equation}\label{psi.r}
(\Psi^H \circ R_H)(M)(S) \cong R_H(M)(\gamma_H(S)) \cong
C^\hdot(W,M[\gamma_H(S)^H]).
\end{equation}
However, points in the orbit $G/H$ fixed under $H$ correspond to
elements $g \in G$ that normalize $H$, so that we have a natural
identification $[G/H]^H \cong [N_H/H]^H \cong W$. Then by the
definition of the functor $\gamma^H$, we have
$$
\gamma_H(S)^H \cong S^H \times W
$$
for any $S \in \Gamma_H$. Therefore \eqref{psi.r} implies that
$\Psi^H(R_H(M)) \in \DM(H,R)$ is supported at $H$, and
$$
(\Psi^H \circ R_H)(M)([H/H]) \cong
C^\hdot(W,M[W]) \cong M,
$$
so that $\bPhi^H(R_H(M)) = \Phi^H(\Psi^H(R_H(M))) \cong M$. Thus the
adjunction map $\wPhi^H(R_H(M)) \to M$ is an isomorphism.
\endproof

\begin{lemma}\label{fs.le}
Any object $M \in \DM(G,R)$ is a finite iterated extension of
objects of the form $R_H(M_H)$, $M_H \in \D(R[W_H])$, $H \subset G$
a subgroup.
\end{lemma}

\proof{} For any $M \in \DM(G,R)$, denote by $\Supp(M)$ the set of
all conjugacy classes of subgroups $H \subset G$ such that
$M([G/H']) \neq 0$ for some subgroup $H' \subset H$. 
Take an object $M \in \DM(G,R)$,
and choose subgroup $H \in \Supp(M)$
that is minimal by inclusion -- that is, no proper subgroup $H'
\subset H$ lies in $\Supp(M)$. Let $M'$ be the cone of the
adjunction map $M \to R_H(\wPhi^H(M))$. Then by \eqref{r.exp},
$\Supp(R_H(\wPhi^H(M)))$ is contained in the set of subgroups $H'
\supset H$, and in particular, $\Supp(R_H(\wPhi^H(M)))$ lies inside
$\Supp(M)$. Then $\Supp(M')$ is also contained in
$\Supp(M)$. However, since $R_H$ is fully faithful, $H$ does not lie
in $\Supp(M')$, so that the inclusion $\Supp(M') \subset \Supp(M)$
is strict. By induction on the cardinality of $\Supp(M)$, this
finishes the proof.
\endproof

\subsection{Maximal Tate cohomology.}

By Lemma~\ref{fs.le} and Lemma~\ref{R.H.le}, the category $\DM(G,R)$
is an iterated extension of full subcategories $\D(R[W_H])$, $H
\subset G$, $W_H = N_H/H$. To describe the gluing data between these
subcategories, one needs to compute the composition functors
\begin{equation}\label{e.h.h}
E^H_{H'} = \wPhi^{H'} \circ R_H:\D(R[H]) \to \D(R[H'])
\end{equation}
for pairs of different subgroups $H,H' \subset G$. This can be done
using a certain natural generalization of Tate cohomology of finite
groups.

\begin{defn}\label{tate.def}
A finitely generated $\Z[G]$-module $M$ is {\em induced} if $M =
\Ind_H^G(M')$, where $H \subset G$ is a proper subgroup, $M'$ is a
finitely generated $\Z[H]$-module, and $\Ind^G_H:\Z[H]\amod \to
\Z[G]\amod$ is the induction functor left-adjoint to the obvious
restriction functor $\Z[G]\amod \to \Z[H]\amod$. For any finitely
generated $\Z[G]$-module $M$, the {\em maximal Tate cohomology
  groups} $\vH^\hdot(G,M)$ are given by
$$
\vH^\hdot(G,M) = \RHom^\hdot_{\D^b_f(Z[G])/\D^b_i(Z[G])}(R,M),
$$
where $\D^b_f(\Z[G])$ is the bounded derived category of the
category of finitely generated $\Z[G]$-modules, and $\D^b_i(\Z[G])
\subset \D^b_f(\Z[G])$ is the smallest Karoubi-closed full triangulated
subcategory containing all induced modules.
\end{defn}

\begin{remark}
The difference with the usual Tate cohomology is that one takes the
quotient not by the subcategory $\D^{pf}(\Z[G]) \subset
\D^b_f(\Z[G])$ of perfect complexes of $\Z[G]$-modules but by the
larger subcategory $\D^b_i(\Z[G]) \subset \D^b_f(\Z[G])$ of induced
modules.
\end{remark}

We note that since the category $\D^b(\Z[G])$ is small, taking the
quotient in Definition~\ref{tate.def} presents no problem. In
effect, one considers the category $I$ of objects $V \in
\D^b_f(\Z[G])$ equipped with a map $V \to \Z$ whose cone lies in
$\D^b_i(\Z[G])$, and one has
\begin{equation}\label{tate.lim.eq}
\vH^i(G,M) = \lim_{\overset{V \in I}{\to}}\Hom(V,M[i])
= \lim_{\overset{V \in I}{\to}}H^i(G,M \lotimes V^*), \qquad i \in \Z,
\end{equation}
where $V^* = \Hom(V,\Z) \in \D^b_f(\Z[G])$ is the object dual to
$V$. Since the index category $I$ is filtered, the limit in the
right-hand side of \eqref{tate.lim.eq} is an exact functor. One can
then use \eqref{tate.lim.eq} to define maximal Tate cohomology with
coefficients.

\begin{defn}\label{tate.big.def}
For any $M \in \D(R[G])$, the {\em maximal Tate cohomology modules}
$\vH^\hdot(G,M)$ are given by
$$
\vH^\hdot(G,M) = \lim_{\overset{V \in I}{\to}}H^\hdot(G,M \lotimes
V^*),
$$
where $I$ and $V^*$ are as in \eqref{tate.lim.eq}.
\end{defn}

One can further refine this to obtain a maximal Tate cohomology
object $\vC^\hdot(G,M) \in \D(R)$ with homology modules
$\vH^\hdot(G,M)$. To do this, it is convenient to use the following
technical gadget. Assume given a complex $V_\idot \in C_{\geq
  o}(\Z[G]\amod)$ of $\Z[G]$-modules concentrated in non-negative
homological degrees, and let $V \in \Fun(\Delta^o,\Z[G]) \cong
\Fun(\ppt_G \times \Delta^o,\Z)$ be the simplicial $\Z[G]$-module
corresponding to $V_\idot$ under the Dold-Kan equivalence
\eqref{dk.eq}. Then for any $M \in \D(R)$, we can consider the
object
\begin{equation}\label{c.p.tate}
\vC^\hdot(G,V_\idot,M) = C_\idot(\Delta^o,R^\hdot\pi_{2*}(V \otimes
\pi_1^*M)) \in \D(R),
\end{equation}
where $\pi_1:\ppt_G \times \Delta^o \to \ppt_G$, $\pi_2:\ppt_G
\times \Delta^o \to \Delta^o$ are the natural projections. This
object is functorial in $M$ and in $V_\idot$ -- for any two
complexes $V_\idot$, $V'_\idot$, any map $f:V_\idot \to V'_\idot$
induces a functorial map
\begin{equation}\label{vc.f}
\vC^\hdot(G,V_\idot,-) \to \vC^\hdot(G,V'_\idot,-).
\end{equation}

\begin{defn}\label{tate.ada.def}
A complex $P_\idot$ of finitely generated $\Z[G]$-modules is
{\em maximally adapted} if
\begin{enumerate}
\item $P_i = 0$ for $i < 0$, $P_0 = \Z$, and for every $i \geq 1$,
  $P_i$ is flat over $\Z$ and induced in the sense of
  Definition~\ref{tate.def}, and
\item for any proper subgroup $H \subset G$, the complex $P_\idot$
  is contractible as a complex of $\Z[H]$-modules.
\end{enumerate}
\end{defn}

\begin{exa}\label{ada.ada}
For any $n$-simplicial pointed $G$-set $X$ that is
termwise-ad\-mis\-sib\-le in the sense of Definition~\ref{adm.G.set} and
$G$-adapted in the sense of Definition~\ref{ada.N.def}, the complex
$P_\idot = C_\idot(X,\Z)$ with its natural $G$-action is maximally
adapted in the sense of Definition~\ref{tate.ada.def}.
\end{exa}

\begin{lemma}\label{tate-ma.le}
For any maximally adapted complex $P_\idot$, and any object $M \in
\D(R[G])$, the homology modules of the object
$\vC^\hdot(G,P_\idot,M)$ of \eqref{c.p.tate} are naturally
identified with the maximal Tate cohomology modules $\vH^\hdot(G,M)$
of Definition~\ref{tate.big.def}. For any two maximally adapted
complexes $P_\idot$, $P'_\idot$ and any map $f:P_\idot \to P'_\idot$
identical in degree $0$, this identification commutes with the map
\eqref{vc.f}.
\end{lemma}

\proof{} For any object $E \in \D(\Delta^o,R)$ represented by a
complex of simplicial $R$-modules, we can apply the Dold-Kan
equivalence termwise and obtain a bicomplex $E^\hdot_\idot$. The
sum-total complex of this bicomplex then represents the homology
object $C_\idot(\Delta^o,E)$. If the bicomplex is of finite length
in either of the two directions, then the sum-total complex is the
same as the product-total complex. Therefore for a complex $V_\idot
\in C_{\geq 0}(\Z[G]\amod)$ of finite length, we have a natural
isomorphism
$$
\vC^\hdot(G,V_\idot,M) \cong C^\hdot(G,M \otimes
C_\idot(\Delta^o,V)) \cong C^\hdot(G,M \otimes V_\idot).
$$
For a general complex $V_\idot$, we can consider its stupid
filtration $F^\hdot V_\idot$ given by
$$
F^lV_i =
\begin{cases}
V_i, &\quad i \leq l,\\
0, &\quad i > l.
\end{cases}
$$
Then $V_\idot = \displaystyle\lim_{\overset{l}{\to}}F^lV_\idot$, and
since homology commutes with direct limits, we have an
identification
\begin{equation}\label{l.lim}
\vC^\hdot(G,V_\idot,M) \cong \lim_{\overset{l}{\to}}C^\hdot(G,M
\otimes F^lV_\idot).
\end{equation}
Since the limit is filtered, its homology modules are then naturally
identified with the corresponding direct limits of the homology
groups $H^\hdot(G,M \otimes F^lV_\idot)$. Now take a maximally
adapted complex $P_\idot$, and consider the double limit
$$
\lim_{\overset{V \in I}{\to}}\lim_{\overset{l}{\to}}H^\hdot(G,M
\otimes F^lP_\idot \otimes V^*).
$$
By Definition~\ref{tate.ada.def}~\thetag{i}, we have $\vH^\hdot(G,M
\otimes P_i) = 0$ for any $i \geq 1$, so that the natural embedding
$$
\lim_{\overset{l}{\to}}H^\hdot(G,M \otimes F^lP_\idot) \to
\lim_{\overset{V \in I}{\to}}\lim_{\overset{l}{\to}}H^\hdot(G,M
\otimes F^lP_\idot \otimes V^*)
$$
is an isomorphism. By Definition~\ref{tate.ada.def}~\thetag{ii}, the
limit \eqref{l.lim} vanishes for any induced $M = \Ind^G_H(M')$, $M' \in
\D(R[H])$, $H \subset G$ a proper subgroup. By the projection
formula, for any $M \in \D(R[G])$ and $V' \in \D^b_i(\Z[G])$, the
product $M \otimes V$ is a direct summand of a finite iterated
extension of objects of this type, so that the natural embedding
$$
\lim_{\overset{V \in I}{\to}}H^\hdot(G,M \lotimes V^*) \to
\lim_{\overset{V \in I}{\to}}\lim_{\overset{l}{\to}}H^\hdot(G,M
\otimes F^lP_\idot \otimes V^*)
$$
induced by the embedding $\Z = P_0 \to P_\idot$ is also an
isomorphism. This finishes the proof.
\endproof

Now note that for any two maximally adapted complexes $P_\idot$,
$P'_\idot$, their product $P_\idot \otimes P'_\idot$ is also
maximally adapted, and the embeddings $\Z = P_0 \to P_\idot$, $\Z =
P'_0 \to P'_\idot$ induce natural maps
$$
P_\idot \to P_\idot \otimes P'_\idot, \qquad P'_\idot \to P_\idot
\otimes P'_\idot.
$$
The corresponding maps \eqref{vc.f} then provide functorial
quasiisomorphisms
$$
\vC^\hdot(G,P_\idot,M) \cong \vC^\hdot(G,P_\idot \otimes P'_\idot,M)
\cong \vC^\hdot(G,P'_\idot,M).
$$
Thus up to a quasiisomorphism, $C_\idot(G,P_\idot,-)$ does not
depend on the choice of a maximally adapted complex $P_\idot$, so
that we can drop it from notation and obtain a maximal Tate
cohomology functor
\begin{equation}\label{tate.c}
\vC^\hdot(G,-):\D(R[G]) \to \D(R).
\end{equation}
Strictly speaking, this functor is only defined up to a
non-canonical isomorphism, but we will ignore this.

Finally, we record one further refinement of maximal Tate
cohomology. For any subgroup $H \subset G$ and any $M \in \D(R[G])$,
the cohomology complex $C^\hdot(H,R)$ carries a natural action of
the group $W_H = N_H/H$, $N_H \subset G$ the normalizer of $H
\subset G$, so that $C^\hdot(H,-)$ is a functor from $\D(R[G])$ to
$\D(R[W_H])$. To obtain a similar maximal Tate cohomology functor
\begin{equation}\label{tate.r.c}
\vC^\hdot(G,-):\D(R[G]) \to \D(R[W_H]),
\end{equation}
it suffices to choose a maximally adapted complex $P_\idot$ of
$H$-modules such that the $H$-action is extended to an action of
$N_H$. Choosing an $H$-adapted pointed $n$-simplicial $N_H$-set $X$
in Example~\ref{ada.ada} does the job.

\subsection{Gluing data.}

We can now describe the gluing functors $E^H_{H'}$ of~\eqref{e.h.h}.
For any two subgroups $H,H' \subset G$, denote
\begin{equation}\label{c.h.h}
c(H',H) = \Hom_G([G/H],[G/H']) \cong [G/H']^H.
\end{equation}
The group $W_H = \Aut_G([G/H])$ acts on this set on the left, and
the group $W_H' = \Aut_G([G/H'])$ acts on the right. Denote by
\begin{equation}\label{b.c.h.h}
\overline{c}(H',H) = W_H \backslash c(H',H) / W_{H'}
\end{equation}
the double quotient. Moreover, assume that $H' \subset G$ lies in
the normalizer $N_H \subset G$, and consider the intersection
$N^H_{H'} = N_H \cap N_{H'} \subset G$. Then $H'$ lies in $N^H_{H'}$
as a normal subgroup, so that we can consider the subgroup $W^H_{H'}
= N^H_{H'}/H' \subset W_{H'} = N_{H'}/H'$, and the associated
induction functor
\begin{equation}\label{ind.w}
\iota:\D(R[W^H_{H'}]) \to \D(R[W_{H'}]).
\end{equation}
By definition, the functor $\wPhi^{H'}$ only depends on $H'$ up to a
conjugation, so that the following results gives a complete
description of the functors $E^H_{H'}$.

\begin{prop}\label{tate.le}
For any two subgroups $H,H' \subset G$, $H \neq H'$, and an object
$M \in \D(R[W_H])$, $E^H_{H'}(M) = 0$ unless $H'$ is conjugate to a
subgroup in $N_H$ containing $H \subset N_H$. If $H \subset H'
\subset N_H$, then we have a natural isomorphism
$$
E^H_{H'}(M) \cong
\iota\left(\vC^\hdot(H'/H,M)\right)\left[\overline{c}(H',H)\right],
$$
where $\vC^\hdot(-,-)$ is the refined maximal Tate cohomology
functor \eqref{tate.r.c}, $\iota$ is the induction functor
\eqref{ind.w}, and $\overline{c}(H',H)$ is the finite set
\eqref{b.c.h.h}.
\end{prop}

\proof{} Take some $M \in \D(R[W_H])$. Assume first that $H' = G$,
so that $W_{H'}$ is trivial, and \eqref{tate.r.c} reduces to
\eqref{tate.c}. Choose a $2$-simplicial finite pointed $G$-set $X$
that is $G$-adapted in the sense of Definition~\ref{ada.N.def}. Then
$R_H(M)$ is given by \eqref{r.g.eq}, and as in the proof of
Lemma~\ref{phi.add}, we have a natural identification
\begin{equation}\label{phi.x}
E^H_G(M) \cong C_\idot(X,R^\hdot\pi_{2*}(\pi_1^*M \otimes T_H))\\
\otimes T_H)).
\end{equation}
Note that by definition, we have $T_H = S(\phi^H)^*T$, where $T \in
\M(W_H,\Z)$ sends a $W_H$-set $S$ to $\Z[S]$. Therefore by base
change, \eqref{phi.x} yields an identification
\begin{equation}\label{phi.y}
E^H_G(M) \cong C_\idot(Y,R^\hdot\pi_{2*}(\pi_1^*M \otimes T)),
\end{equation}
where we denote $Y=X^H$. If $H \subset G$ is not normal, so that
$N_H \subset G$ is a proper subgroup, then for every subgroup $H_0
\subset N_H$ containing $H$, the fixed point set $Y^{H_0/H} =
X^{H_0}$ is contractible by Definition~\ref{ada.N.def}~\thetag{ii},
and $E^H_G(M)=0$ by Lemma~\ref{contr.ex.le}. If $H \subset G$ is
normal, so that $G = N_H$, then $Y$ is a $W_H$-adapted
$2$-simplicial pointed $W_H$-set. By Definition~\ref{homo.point.def}
and \eqref{shuffle.eq}, the right-hand side of \eqref{phi.y} is then
given by
\begin{equation}\label{c.c}
C_\idot(X,R^\hdot\pi_{2*}(\pi_1^*M \otimes T)) \cong
C_\idot(\Delta^o,R^\hdot\pi_2^*(\pi_1^*M \otimes \delta^*Y^*T)).
\end{equation}
But by definition, the simplicial $\Z[W_H]$-module $\delta^*Y^*T$
corresponds to the maximally adapted complex $P_\idot =
C_\idot(Y,\Z)$ of Example~\ref{ada.ada}, so that the right-hand side
of \eqref{c.c} is exactly \eqref{c.p.tate}.

In the general case, recall that by \eqref{orb.fib}, the category
$\Gamma_{H'}$ is naturally identified with the category of finite
$G$-sets $S$ equipped with a $G$-equivariant map $S \to
[G/H']$. Then every element $f:[G/H] \to [G/H']$ of the set
$c(H',H)$ of \eqref{c.h.h} defines in particular an object in
$\Gamma_{H'}$. This object is an orbit of the form $[H'/H_f]$, where
the subgroup $H_f \subset H' \subset G$ is conjugate to $H$ in $G$,
and well-defined up to conjugation in $H'$. For any $S \in
\Gamma_{H'}$, we then have a natural identification
$$
(\gamma^{H'}(S))^H \cong \coprod_{f \in c(H',H)}S^{H_f},
$$
and by definition, this gives an isomorphism
$$
\Psi^{H'}T_H \cong \bigoplus_{f \in c(H',H)}T_{H_f},
$$
where $T_{H_f} \in \M(H',\Z)$ corresponds to the subgroup $H_f
\subset H'$. The group $W_H$ acts on $C(H',H)$, and then for any $f
\in c(H',H)$, its stabilizer $W_{H_f} = \Aut_{H'}([H'/H_f]) \subset
W_H$ acts on $T_{H_f}$. Therefore by base change, we have
\begin{equation}\label{f.spl}
E^H_{H'} \cong \bigoplus_{f \in W_H\backslash c(H',H)}\wt{E}^{H_f}_{H'},
\end{equation}
where $\wt{E}^{H_f}_H$ are the functors \eqref{e.h.h} for $H_f$
considered as a subgroup in $H'$. If $H_f \subset H'$ is not normal
-- that is, $H'$ is not conjugate to a subgroup in $N_H$ -- then the
right-hand side of \eqref{f.spl} vanishes. In the case $H \subset H'
\subset N_H$, the group $W_{H'}$ acts on the quotient $W_H
\backslash c(H',H)$, and the stabilizer of any element is conjugate
to $W^H_{H'} \subset W_{H'}$. Therefore \eqref{f.spl} can be
rewritten as
$$
E^H_{H'} \cong \iota\left(\wt{E}^H_{H'}\right)\left[\overline{c}(H',H)\right],
$$
where $\iota$ is the induction functor \eqref{ind.w}, and
$\overline{c}(H',H)$ is the finite set \eqref{b.c.h.h}. Plugging in
the expression for $\wt{E}^H_{H'}$, we get the claim.
\endproof

\subsection{Autoequivalences.}

Proposition~\ref{tate.le} is essentially a reformulation of
\cite[Corollary 7.10]{mackey}. The full story in \cite{mackey} is
considerably more elaborate. Roughly speaking, one refines the
gluing data description obtained in Proposition~\ref{tate.le} to a
certain DG coalgebra, and then uses it to give a full description of
the category $\DM(G,R)$. We do not attempt to re-do it in this
paper. However, we will prove one useful corollary of
Proposition~\ref{tate.le}.

Consider the full subcategory
$$ 
\DS(\Gamma_G,\Gamma_{G+},R) \subset \D(S\Gamma_G \times
\Gamma_{G+},R)
$$
spanned by special objects, as in \eqref{ds.2.eq}. As in
\eqref{q.st.g}, we have a natural equivalence
$$
\DS(\Gamma_G,\Gamma_{G+},R) \cong \DS_{\Id \times
  \Inj}(\Gamma_G \times \Gamma_G,R).
$$
As in \eqref{m.eq}, the cartesian product functor $m:\Gamma_G \times
\Gamma_G \to \Gamma_G$ is a morphism in the sense of
Definition~\ref{adm.mor.def}, so that we have the pullback functor
$$
S(m)^*:\DS(\Gamma_G,R) \to \DS(\Gamma_G,\Gamma_{G+},R).
$$
Moreover, for any simplicial pointed finite $G$-set $X:\Delta^o \to
\Gamma_{G+}$, we have the pullback functor
$$
(\id \times X)^*:\DS(\Gamma_G,\Gamma_{G+},R) \to
\DS(\Gamma_G,\Delta^o,R),
$$
where $\DS(\Gamma_G,\Delta^o,R) \subset \D(S\Gamma_G \times
\Delta^o,R)$ is again defined as in \eqref{ds.2.eq}.

\begin{defn}\label{sm.def}
The {\em smash product} $M \wedge X$ of a derived Mackey functor
$M \in \DM(G,R)$ and a simplicial finite pointed $G$-set $X \in
\Delta^o\Gamma_{G+}$ is given by
$$
M \wedge X = L^\hdot\pi_!(\id \times X)^*S(m)^*M,
$$
where $\pi:\Delta^o \times S\Gamma_G \to S\Gamma_G$ is the
natural projection.
\end{defn} 

Let us record some of the obvious properties of the smash product of
Definition~\ref{sm.def}. First of all, if $G$ is trivial, so that
$\DM(G,R) \cong \D(R)$, then we have
\begin{equation}\label{triv.we}
M \wedge X \cong M \otimes \overline{C}_\idot(X,\Z),
\end{equation}
where $\overline{C}_\idot(X,\Z) = C_\idot(X,\Z)/\Z\cdot\{o\}$ is the
reduced chain complex of the pointed simplicial set $X$. In the
general case, by base change, we have a natural identification
$$
(M \wedge X)(S) \cong C_\idot(X \wedge S_+,M)
$$
for any $S \in \Gamma_G$, where the right-hand side is as in
Definition~\ref{homo.point.def}. This is additive in $S$, so that
for any $M \in \DM(G,R)$ and any $X \in \Delta^o\Gamma_{G+}$, the
smash product $M \wedge X$ is also a derived Mackey functor. The
K\"unneth formula provides a natural quasiisomorphism
$$
(M \wedge X_1) \wedge X_2 \cong M \wedge (X_1 \wedge X_2)
$$
for any $X_1,X_2 \in \Delta^o\Gamma_{G+}$. For any subgroup $H
\subset G$, the obvious isomorphism $\rho^H(X_1 \wedge X_2) \cong
\rho^H(X_1) \wedge \rho^H(X_2)$ induces by adjunction a projection
formula isomorphism
$$
\gamma^H(\rho^H(X) \wedge S_+) \cong X \wedge \gamma^H(S)_+,
$$
and this induces a natural functorial isomorphism
$$
\Psi^H(M \wedge X) \cong \Psi^H(M) \wedge \rho^H(X).
$$
This isomorphism respects the $W_H$-action, thus gives an
isomorphism
\begin{equation}\label{psi.we}
\wPsi^H(M \wedge X) \cong \wPsi^H(M) \wedge \rho^H(X)
\end{equation}
of extended functors of \eqref{psi.wt}. Also, we have $\phi^H(X
\wedge S_+) \cong X^H \wedge \phi^H(S)_+$, and this induces a
functorial isomorphism
$$
\Infl^H(M) \wedge X \cong \Infl^H(M \wedge X^H).
$$
By adjunction, we obtain an isomorphism
\begin{equation}\label{phi.we}
\Phi^H(M \wedge X) \cong \Phi^H(M) \wedge X^H
\end{equation}
for any $M \in \DM(G,R)$ and any $X \in \Delta^o\Gamma_{G+}$.

\begin{defn}\label{hom.ada}
A simplicial finite pointed $G$-set $X \in \Delta^o\Gamma_{G+}$ is
{\em homologically adapted} to a normal subgroup $N \subset G$ if
for any subgroup $H \subset G$ containing $N$, the reduced chain
complex $\overline{C}_\idot(X^H,\Z)$ is quasiisomorphic to $\Z$,
while for any subgroup $H \subset G$ not containing $N$, the complex
$\overline{C}_\idot(X^H,\Z)$ is acyclic.
\end{defn}

\begin{exa}
For any $n$-simplicial finite pointed $G$-set $X$ adapted to $N$ in
the sense of Definition~\ref{ada.N.def}, the diagonal $\delta^*X$ is
homologically adapted to $N$ in the sense of
Definition~\ref{hom.ada}.
\end{exa}

\begin{lemma}
Assume given a normal subgroup $N$ and a simplicial finite pointed $G$-set $X \in
\Delta^o\Gamma_{G+}$ homologically adapted to $N$ in the sense of
Definition~\ref{hom.ada}. Then for any $E \in \DM(G,R)$, we have a
natural isomorphism
$$
X \wedge E \cong \Infl^N(\Phi^N(E)).
$$
\end{lemma}

\proof{} Since $X^G$ has non-trivial homology, it is not empty, so
that we have a nontrivial $G$-equivariant map $i:[1]_+ \to X$. This
map then automatically induces a quasiisomorphism $\Z \cong
\overline{C}_\idot(X^H,\Z)$ for any $H \subset G$ containing
$N$. Consider the natural map
$$
\wt{i}:E \cong [1]_+ \wedge E \to X \wedge E
$$
induced by $i$. Then by \eqref{phi.we}, $X \wedge E$ lies inside the
subcategory $\DM_N(G,R)$, so that $\wt{i}$ factors through a map
$$
\Infl^N(\Phi^N(E)) \to X \wedge E,
$$
and again by \eqref{phi.we}, this map becomes an isomorphism after
applying the conservative functor $\bPhi_\idot$ of
Lemma~\ref{phi.dot.le}.
\endproof

\begin{defn}\label{sph.def}
A simplicial finite pointed $G$-set $X \in \Delta^o\Gamma_{G+}$ is a
{\em homological sphere} if for any subgroup $H \subset G$, the
reduced chain complex $\overline{C}_\idot(X^H,\Z)$ is
quasiisomorphic to $\Z$ placed in some degree $d_H$.
\end{defn}

\begin{exa}\label{regu.sph.exa}
For any finite set $S$, let $X(S) = I^S/\overline{I^S}$ be the
quotient of the product of copies of the simiplicial interval $I$
numbered by elements $s \in S$ by its boundary
$$
\overline{I^S} = \coprod_{s \in S}\{s\} \times I^{S \setminus \{s\}}
\copr \coprod_{s \in S}\{t\} \times I^{S \setminus \{t\}},
$$
where $s,t \in I$ are the endpoints. Then
$\overline{C}_\idot(X(S),\Z)$ is $\Z$ placed in degree $|S|$, the
cardinality of $S$. Let $G$ act on $X(G)$ via its action on
itself. Then for any subgroup $H \subset G$, we have $X(G)^H \cong
X(G/H)$, so that $X(G)$ is a homological sphere in the sense of
Definition~\ref{sph.def} (with $d_H=|G/H|$, the cardinality of the
orbit $G/H$).
\end{exa}

\begin{exa}
More generally, for any finite-dimensional representation $V$ of the
group $G$ over $\RR$, consider its one-point compactification $S_V$,
choose its finite triangulation compatible with the $G$-action, and
let $X$ be the corresponding simplicial pointed $G$-set. Then $X$ is
a homological sphere, with $d_H = \dim_{\RR}V^H$. If $V=\RR[G]$ is
the regular representation, one recovers $X(G)$ of
Example~\ref{regu.sph.exa}.
\end{exa}

\begin{prop}\label{equi.prop}
Assume given a simplicial finite pointed $G$-set $X \in
\Delta^o\Gamma_{G+}$, and assume that $X$ is a homological sphere in
the sense of Definition~\ref{sph.def}. Then the endofunctor $E_X$ of
the category $\DM(G,R)$ sending $M$ to $M \wedge X$ is an
autoequivalence.
\end{prop}

\proof{} For any subgroup $H \subset G$ and $M \in \DM(G,R)$, the
isomorphisms \eqref{psi.we}, \eqref{phi.we} and \eqref{triv.we}
provide a functorial isomorphism $\wPhi^H(M \wedge X) \cong
\wPhi^M(M)[d_H]$, where $\wPhi^H$ is the functor \eqref{wphi.H}. By
adjunction, we obtain functorial base change maps
\begin{equation}\label{r.we}
R_H(M) \wedge X \to R_H(M)[d_H], \qquad M \in \D(R[W_H]),
\end{equation}
where $R_H$ is the adjoint functor of Lemma~\ref{R.H.le}. Applying
the functors $\wPhi^{H'}$, $H' \subset G$ another subgroups, we
obtain functorial maps
\begin{equation}\label{d.d.eq}
E^H_{H'}(M)[d_{H'}] \to E^H_{H'}(M)[d_H]
\end{equation}
where we have used the identification $E^H_{H'}(M)[d_{H'}] \cong
\wPhi^{H'}(R_H(M) \wedge X)$ provided by \eqref{phi.we} and
\eqref{triv.we}.  Assume for the moment that all the maps
\eqref{d.d.eq} are isomorphisms. Then by
Lemma~\ref{phi.dot.le}~\thetag{ii}, the maps \eqref{r.we} are also
isomorphisms, so that $E_X$ preserves the full subcategories
$R_H(\D(R[H])) \subset \DM(G,R)$ of Lemma~\ref{fs.le}, and induces a
homological shift on each of these subcategories. Moreover, by
adjunction, $E_X$ is fully faithful on objects of the form $R_H(M)$,
$H \subset G$, $M \in \D(R[H])$. Then by Lemma~\ref{fs.le}, it is
fully faithful on the whole $\DM(G,R)$, and since homological shifts
are autoequivalences, it is essentially surjective, again by
Lemma~\ref{fs.le}.

It remains to prove that all the functorial maps \eqref{d.d.eq} are
isomorphisms. By Proposition~\ref{tate.le}, it suffices to consider
the case when $H$ lies in $H'$ as a normal subgroup, and in this
case, the cone of the map \eqref{d.d.eq} is a sum of terms of the
form
$$
\vC^\hdot(H'/H, M \otimes \overline{C}_\idot(X^H/X^{H'},\Z)).
$$
But since the finite simplicial pointed $W_H$-set $X^H/X^{H'}$ has
no non-trivial points fixed under $H'/H \subset W_H$, the reduced
chain homology complex $\overline{C}_\idot(X^H/X^{H'},\Z)$ is
induced in the sense of Definition~\ref{tate.def}, and the maximal
Tate cohomology in question indeed vanishes.
\endproof

\section{Derived normal systems.}\label{ns.sec}

To proceed further, we need to describe a derived version of the
theory of normal systems introduced in Subsection~\ref{norm.subs}.

\subsection{The comparison theorem.}

Consider the fibration $\nu$ of \eqref{to.N.eq}, and recall that its
fibers are identified with the categories $\Gamma_W$, $W = G/N$, $N
\in \Nn(G)$ a normal cofinite subgroup in $G$, while the transition
functors are the functors $\phi^{N'/N}$, $N \subset N' \subset
G$. Note that the $S$-construction of Definition~\ref{SC.def} can be
applied fiberwise to the fibration $\hGamma_G/\Nn(G)^o$. Namely, let
$S(\hGamma_G/\Nn(G)^o)$ be the category of pairs $\langle [n],f
\rangle$ of $[n] \in \Delta$ and a functor $f:[n] \to \hGamma_G$
such that $\nu \circ f:[n] \to \Nn(G)^o$ is the constant functor
onto an object $N \in \Nn(G)^o$, with morphisms from $\langle [n],f
\rangle$ to $\langle [n'],f' \rangle$ given by a pair of a morphism
$\phi:[n] \to [n']$ and a morphism $\alpha:f' \circ \phi \to f$ such
that for any $i,j \in [n]$, the commutative square \eqref{n-ex-3}
induced by \eqref{n-ex-2} is cartesian in $\hGamma$. Then $\nu$
induces a forgetful functor
$$
S(\nu):S(\hGamma_G/\Nn^o(G)) \to \Nn(G).
$$
This functor is a cofibration whose fibers are the categories
$S\Gamma_W$, and whose transition functors are functors
$S(\phi^{N'/N})$. Say that a morphism $f$ in the category
$S(\hGamma_G/\Nn(G)^o)$ is {\em special} if $S(\nu)(f)$ is
invertible, and $f$ is special in the sense of
Definition~\ref{DS.def} when considered as a morphism in a fiber of
the cofibration $S(\nu)$. Say that an object $E \in
\D(S\hGamma_G,R)$ is {\em special} if it can be represented by
complex $E_\idot$ such that $E_\idot(f)$ is a quasiisomorphism for
any special morphism $f$, and denote by
\begin{equation}\label{ds.rel}
\DS(\hGamma_G/\Nn(G)^o,R) \subset \D(S(\hGamma_G/\Nn(G)^o),R)
\end{equation}
the full subcategory spanned by special objects. Then every special
object $E \in \DS(\hGamma_G/\Nn(G)^o,R)$ defines a collection of
objects $E_N \in \DS(\wGamma_G,R)$, $N \in \Nn(G)$, $W = G/N$, and
transition morphisms
\begin{equation}\label{der.norm.trans}
S(\phi^{N'/N})_!E_N \to E_{N'}
\end{equation}
for any $N,N' \in \Nn(G)$, $N \subset N'$.

\begin{defn}\label{der.norm.def}
A {\em derived normal system} of $G$-Mackey profunctors is a special
object $E \in \DS(\hGamma_G/\Nn(G)^o,R)$ such that for any $N \in
\Nn(G)$, $W = G/N$, the object $E_N \in \DS(\wGamma_W,R)$ lies in
$\DM(W,R) \subset \DS(\wGamma_W,R)$, and for any $N,N' \in \Nn(G)$,
$N \subset N'$, the map
$$
\Phi^{N'/N}E_N \to E_{N'}
$$
induced by \eqref{der.norm.trans} is an isomorphism.
\end{defn}

By definition, derived normal systems form a full triangulated
subcategory in $\DS(\hGamma_G/\Nn(G)^o,R)$; we denote this category
by $\DN(G,R)$. For any integer $i$, we let
$$
\DN^{\leq i}(G,R) = \DN(G,R) \cap \DS^{\leq i}(\hGamma_G/\Nn(G)^o,R)
\subset \DS(\hGamma_G/\Nn(G)^o,R).
$$
We note that since the functors $\Phi^{N'/N}$ are only right-exact
with respect to the standard $t$-structure, these subcategories do
not automatically give a $t$-structure on $\DN(G,R)$. Nevertheless,
they are perfectly well-defined. So is the subcategory
$$
\DN^-(G,R) = \bigcup_i \DN^{\leq i}(G,R) \subset \DN(G,R),
$$
and the truncation functor $\tau$ of \eqref{tau.eq} induces a
natural functor
\begin{equation}\label{tau.N}
\tau:\DN^{\leq 0}(G,R) \to \N(G,R).
\end{equation}
Now let $\wt{S}(\hGamma_G/\Nn(G)^o)$ be the category of diagrams
\eqref{codomik} in $S(\hGamma_G/\Nn(G)^o)$ with special $s_1$,
$s_2$, and note that the projections
$\pi_1,\pi_2:\wt{S}(\hGamma_G/\Nn(G)^o) \to S(\hGamma_G/\Nn(G)^o)$
are cocartesian over $\Nn(G)$. Therefore by Lemma~\ref{bc.sga},
\eqref{sp.eq} provides a functor
$$
\Sp:\D(S(\hGamma_G/\Nn(G)^o),R) \to \DS(\hGamma_G/\Nn(G)^o,R)
$$
left-adjoint to the embedding \eqref{ds.rel} that induces the
specialization functor of Lemma~\ref{sp.le} on the fiber over every
$N \in \Nn(G)^o$. Then the functor $\phi$ of \eqref{phi.eq} induces
a functor
$$
S(\phi):S\wGamma_G \times \Nn(G) \to S(\hGamma_G/\Nn(G)^o),
$$
and $S(\phi)$ induces a pair of adjoint functors
\begin{equation}\label{s.phi.eq}
\begin{aligned}
S(\phi)^*:&\DS(\hGamma_G,R) \to \DS(\wGamma_G,\Nn(G),R),\\
S(\phi)_! = \Sp \circ L^\hdot S(\phi)_!:&\DS(\wGamma_G,\Nn(G),R) \to \DS(\hGamma_G,R),
\end{aligned}
\end{equation}
where as in \eqref{ds.2.eq}, $\DS(\wGamma_G,\Nn(G),R) \subset
\D(S\wGamma_G \times \Nn(G),R)$ is the full triangulated subcategory
spanned by objects $E$ such that the corresponding functor
$S\wGamma_G \to \D(\Nn(G),R)$ is special in the sense of
Definition~\ref{DS.def}. If we let $p:S\wGamma_G \times \Nn(G) \to
S\wGamma_G$ be the natural projection, then we obtain a pair of
adjoint functors
\begin{equation}\label{phi.infl.der}
\begin{aligned}
\Phi &= S(\phi)_! \circ p^*:\DS(\wGamma_G,R) \to
\DS(\hGamma_G/\Nn(G),R),\\
\Infl &= R^\hdot p_* \circ S(\phi)^*:\DS(\hGamma_G/\Nn(G),R) \to
\DS(\wGamma_G,R).
\end{aligned}
\end{equation}
Explicitly, for any $E \in \DS(\wGamma_G,R)$, we have $\Phi(E)_N =
S(\phi^N)_!(E)$ for any cofinite normal $N \subset G$, and the
transition morphisms \eqref{der.norm.trans} are induced by the
natural isomorphisms $\phi^{N'/N} \circ \phi^N \cong
\phi^{N'}$. Thus by Lemma~\ref{phi.add}, for any $E \in \wDM(G,R)
\subset \DS(\wGamma_G,R)$, $\Phi(E)$ is a derived normal system in
the sense of Definition~\ref{der.norm.def}, so that the functors
$\Infl$ and $\Phi$ of \eqref{phi.infl.der} induce a pair of adjoint
functors between the categories $\wDM(G,R)$ and $\DN(G,R)$.

\begin{prop}\label{can.filt.der.prop}
Assume that the group $G$ is finitely generated. Then for any
integer $i$, the functors \eqref{phi.infl.der} induce a pair of
inverse equivalences between the triangulated categories $\wDM^{\leq
  i}(G,R)$ and $\DN^{\leq i}(G,R)$, and they also induce
equivalences between $\wDM^-(G,R)$ and $\DN^-(G,R)$.
\end{prop}

\proof{} Since both functors of \eqref{phi.infl.der} are
triangulated, in particular commute with shifts, it suffices to
consider the case $i=0$. By Lemma~\ref{der.noder.le}, the functor
$\Phi$ of \eqref{phi.infl.der} sends $\wDM^{\leq 0}(G,R)$ into
$\DN^{\leq 0}(G,R)$. By definition, for any $E \in \DN^{\leq
  0}(G,R)$ with components $E_N$, $N \in \Nn(G)$, and any admissible
$G$-set $S \in \wGamma_G$, the value of the inflation $\Infl(E)$ at
$S$ is given by
\begin{equation}\label{infl.lim}
\Infl(E)(S) = \dlim_{\overset{N}{\gets}}E_N(S^N)),
\end{equation}
where $\lim^\hdot$ in the right-hand is the derived inverse limit
functor, taken over $N \in \Nn(G)$. Thus in general, even for $E \in
\DN^{\leq 0}(G,R)$, $\Infl(E)$ does not have to lie in $\wDM^{\leq
  0}(G,R)$. However, since the group $G$ is finitely generated, for
any integer $i \geq 1$, there exists at most a finite number of
cofinite normal subgroups $N \subset G$ with $|G/N| \leq i$, so that
their intersection $N_i \subset G$ is also a normal cofinite
subgroup. The subset $\Nn(G)$ formed by the subgroups $N_i$ is
cofinal, so that we have
\begin{equation}\label{dlim.eq}
\dlim_{\overset{N}{\gets}}E_N(S^N)) \cong
\dlim_{\overset{i}{\gets}}E_{N_i}(S^{N_i}).
\end{equation}
By \eqref{tel.eq}, the limit in the right-hand side is isomorphic to
the telescope $\Tel(E_{N_i}(S^{N_i}))$ of the inverse system
$E_{N_i}(S^{N_i})$. Therefore $\Infl$ at least sends $\DN^{\leq
  0}(G,R)$ into $\wDM^{\leq 1}(G,R)$. Moreover, for any $E \in
\DN^{\leq 0}(G,R)$, the truncation at $1$ of $\Infl(E)(S)$ is given
by
\begin{equation}\label{R.lim.eq}
R^1\lim_{\overset{i}{\gets}}E_i(S^{N_i}),
\end{equation}
where we denote $E_i = \tau(E)_{N_i}$, and $\tau(E) \in \N(G,R)$ is
obtained by the truncation functor \eqref{tau.N}.  The projective
system $E_i(S^{N_i})$, $i \geq 1$ is a part of the natural
projective system \eqref{infl} for the normal system $\tau(E)$. In
particular, its transition maps are surjective. Therefore the
group \eqref{R.lim.eq} vanishes, so that for a finitely generated
group $G$, we have
$$
\Infl\left(\DN^{\leq 0}(G,R)\right) \subset \wDM^{\leq 0}(G,R) \subset
\wDM(G,R).
$$
Moreover, as in the proof of Lemma~\ref{phi.add}, we have
\begin{equation}\label{phi.S.N}
\Phi^N(E)(S) \cong C_\idot(S_+ \wedge X_N,E)
\end{equation}
for any $N \in \Nn(G)$, $W = G/N$, $S \in \Gamma_W \subset
\wGamma_G$, $E \in \wDM(G,R)$, $X_N$ an $N$-adapted pointed
admissible $n$-simplicial $G$-set in the sense of
Definition~\ref{ada.N.def}. By Lemma~\ref{prod.le}~\thetag{i}, for
any $E \in \DN^{\leq 0}(G,R)$, we have
$$
C_\idot(S_+ \wedge X_N,\Tel(\Infl^{N_i}E_{N_i})) \cong
\Tel(C_\idot(S_+ \wedge X_N,\Infl^{N_i}E_{N_i})),
$$
so that this imples that
$$
\Phi^N(\Infl(E)) \cong
\Tel(\Phi^N(\Infl^{N_i}(E_{N_i}))),
$$
and as soon as $N_i \subset G$ becomes contained in $N \subset G$,
the inverse system in the right-hand side becomes the constant
system with value $E_N$. Thus the adjunction map
$$
\Phi(\Infl(E)) \to E
$$
is an isomorphism for any $E \in \DN^{\leq 0}(G,R)$.

To finish the proof, we have to prove that the adjunction map $E \to
\Infl(\Phi(E))$ is an isomorphism for any $E \in \wDM^{\leq
  0}(G,R)$. Choose adapted $2$-simplicial sets $X_N = C(ES_N)$ for
all cofinite normal subgroups $N \subset G$ as in
Lemma~\ref{ada.ex.le}, so that effectively, we have
$$
X_{N'} \subset X_N
$$
whenever $N' \subset N$, and we have $[1]_+ \cong \bigcap_NX_N$.
Then \eqref{phi.S.N} gives an isomorphism
$$
\Infl^N(\Phi^N(E))(S) \cong C_\idot(S_+ \wedge X_N,E)
$$
for any $S \in \wGamma_G$ and $E \in \wDM(G,R)$, and by
\eqref{infl.lim} and \eqref{dlim.eq}, this gives an isomorphism
$$
\Infl(\Phi(E))(S) = \dlim_{\overset{i}{\gets}}C_\idot(S_+ \wedge
X_i,E)),
$$
where we denote $X_i = X_{N_i}$. Then by
Lemma~\ref{prod.le}~\thetag{ii}, for $E \in \wDM^{\leq 0}(G,R)$ we
have
$$
\dlim_{\overset{i}{\gets}}C_\idot(S \wedge X_i,E)) \cong
C_\idot(\cap_i(S_+ \wedge X_i),E) = C_\idot(S_+ \wedge \cap_iX_i,E),
$$
and since $\cap_iX_i$ is the constant $2$-simplicial set $[1]_+$,
the right-hand side is naturally identified with $E(S)$.
\endproof

\subsection{Corollaries.}

Comparing Proposition~\ref{can.filt.der.prop} and
Proposition~\ref{can.pro.prop}, we see that the situation in the
derived case is better: every derived Mackey profunctor bounded from
above is separated in the derived sense. In particular, this is true
for a Mackey profunctor $E$ considered as a derived Mackey
profunctor via the embedding \eqref{M.DM.eq}, even if $M$ is not
separated in the sense of Definition~\ref{sepa.def}. The reason for
the discrepancy is that the fixed points functor $\Phi$ is only
exact on the right, while the inverse limit functor is only exact on
the left. Thus even for $E \in \wM(G,R) \subset \wDM(G,R)$,
$\Infl(\Phi(E))$ can contain contributions coming from the derived
functor $R^1\lim_{\gets}$ that do not appear in
\eqref{can.filt.mack}. Thus we could expect the derived theory to
allow us to improve results in the underived setting. Here is a
first example.

\begin{lemma}\label{add.der.pro.corr}
For any finitely generated group $G$, the natural embedding
$\wDM^-(G,R) \subset \DS^-(\wGamma_G,R)$ admits a left-adjoint
additivization functor
$$
\Add:\DS^-(\wGamma_G,R) \to \wDM^-(G,R).
$$
This functor is right-exact with respect to the standard
$t$-structures and induces the functor
\begin{equation}\label{add.under}
\Add = \tau \circ \Add \circ q^*:\Fun(Q\wGamma_G,R) \to \wM(G,R)
\end{equation}
left-adjoint to the embedding $\wM(G,R) \subset \Fun(Q\wGamma_G,R)$.
\end{lemma}

\proof{} If $G$ is finite, so that $\wDM(G,R) = \DM(G,R)$, take the
functor \eqref{add.der} provided by Proposition~\ref{add.prop}. In
the general case, applying $\Add$ fiberwise to the fibration
$S\hGamma_G \to \Nn(G)$ gives a functor
$\Add:\DS^-(S\hGamma_G/\Nn(G),R) \to \DS^-(S\hGamma_G/\Nn(G),R)$,
and as in Lemma~\ref{Add.pro.corr}, the functor
$$
\Add = \Infl \circ \Add \circ \Phi:\DS^-(\wGamma_G,R) \to \wDM^-(G,R)
$$
does the job. By construction, it is right-exact with respect to the
standard $t$-structures, so that by adjuncion, \eqref{add.under}
indeed defines a functor adjoint to the embedding $\wM(G,R) \subset
\Fun(Q\wGamma_G,R)$.
\endproof

As a consequence of this, for finitely generated groups, all the
material of Subsection~\ref{pro.fp.subs} that depends on
Lemma~\ref{Add.pro.corr} can be now generalized to arbitrary Mackey
profunctors. This concerns the categorical fixed points functor
$\Psi^H$ with its refinement $\wt{\Psi}^H$, the inflation functor
$\Infl^H$, the geometric fixed points functor $\Phi^H$, and the
product \eqref{pro.ma.pro}.

\medskip

We can now give the derived versions of all this material. Firstly,
we can extend the geometric fixed points functor $\Phi^H$ to the
case of a subgroup $H \subset G$ that is not cofinite -- indeed, the
fixed points functor $\phi^H:\wGamma_G \to \wGamma_H$ and the
inflation functor $\Infl^H = S(\phi^H)^*$ are perfectly well defined
in this case, and the inflation has a left-adjoint $\Phi^H$ given by
$$
\Phi^H = \Add \circ S(\phi^H)_!,
$$
where $\Add$ is the additivization functor provided by
Corollary~\ref{add.der.pro.corr}. If $H \subset G$ is normal, then
the adjunction map $\Phi^H \circ \Infl^H \to \Id$ is obviously an
isomorphism, so that $\Infl^H$ is still fully faithful, as in
Lemma~\ref{infl.der}. The isomorphism \eqref{phi.psi.iso.der} also
extends to the case when $H \subset G$ is cofinite, but $H' \subset
H \subset G$ is not.

Also, we can extend the product \eqref{ma.produ} to derived Mackey
functors. Indeed, the cartesian product $m:\wGamma_G \times
\wGamma_G \to \wGamma_G$ is a morphism in the sense of
Definition~\ref{adm.mor.def} with respect to the admissible classes
of all maps. Therefore for any rings $R_1$, $R_2$, we can define the
product $E_1 \circ E_2 \in \wDM(G,R_1 \otimes R_2)$ of two derived
Mackey profunctors $E_1 \in \wDM(G,R_1)$, $E_2 \in \wDM(G,R_2)$ as
\begin{equation}\label{der.w.pro}
E_1 \circ E_2 = \Add(S(m)_!E_1 \boxtimes E_2),
\end{equation}
where again, $\Add$ is the additivization functor of
Corollary~\ref{add.der.pro.corr}. The fixed points functors
$\Psi^H$, $\Phi^H$ are tensor functors with respect to this
product. The product is associative and commutative, with the unit
given by the {\em derived completed Burnside Mackey profunctor}
$\wA_\idot \in \wDM(G,\Z)$ expressed as
$$
\wA_\idot = \Add(S(e)_!\Z),
$$
where $e:\ppt \to \wGamma_G$ is the embedding onto the one-point
$G$-orbit $[G/G] \in \Gamma_G$. The {\em derived completed Burnside
  ring} $\wAG_\idot$ of the group $G$ is the DG ring of
endomorphisms of $\wA_\idot \in \wDM(G,\Z)$. By adjunction, it is
given by
$$
\wAG_\idot = \wA_\idot([G/G]),
$$
and since $\wA_\idot$ is the unit object for the tensor product,
$\wAG_\idot$ acts on any $E \cong \wA_\idot \circ E \in \wDM(G,R)$ --
more precisely, we have an action map
\begin{equation}\label{wAG.acts}
\wAG_\idot \otimes E \to E
\end{equation}
functorial with respect to $E$.

Another thing we can now extend to derived Mackey profunctors is the
notion of a smash product of Definition~\ref{sm.def}.

\begin{defn}\label{w.sm.def}
Assume that the group $G$ is finitely generated. The {\em smash
  product} $M \wedge X$ of a derived Mackey profunctor $M \in
\wDM^-(G,R)$ and a simplicial admissible pointed $G$-set $X \in
\Delta^o\wGamma_{G+}$ is given by
$$
M \wedge X = \Add(L^\hdot\pi_!(\id \times X)^*S(m)^*M) \in \wDM^-(G,R),
$$
where $\pi$, $(\id \times X)^*$ and $S(m)^*$ are as in
Definition~\ref{sm.def}.
\end{defn}

\begin{remark}
There is a relation between smash products and the Mackey functor
product \eqref{der.w.pro} similar to \eqref{triv.we}, but since in
this paper, we do not have enough technology to explore it properly,
we will return to it elsewhere.
\end{remark}

Proposition~\ref{equi.prop} also generalizes to Mackey profunctors
but with a certain twist. Assume given a sequence $\{d_\idot\}$ of
integers numbered by cofinite subgroups $H \subset G$, and assume
that the sequence is non-descreasing in the following sense: for any
$H \subset H' \subset G$, we have $d_H \geq d_{H'}$. For any integer
$n$, let
\begin{equation}\label{d.idot}
\wDM^{\leq n - d_\idot}(G,R) \subset \wDM^-(G,R) \subset
\wDM(G,R)
\end{equation}
be the full subcategory spanned by objects $M \in \wDM(G,R)$ such
that $M([G/H])$ lies in $\D^{\leq n - d_H}(R)$ for any cofinite
subgroup $H \subset G$.

\begin{lemma}\label{d.le}
Assume that the group $G$ is finitely generated. Then an object $M
\in \wDM^-(G,R)$ lies in the subcategory \eqref{d.idot} if an only
if for any cofinite subgroup $H \subset G$, $\bPhi^H(M) =
\Phi^H(\Psi^H(M))$ lies in $\D^{\leq n-d_H}(R) \subset \D(R)$.
\end{lemma}

\proof{} Since the sequence $\{d_\idot\}$ is non-decreasing, $M$
lies in the subcategory \eqref{d.idot} if and only if $\Psi^H(M)$
lies in $\wDM^{\leq n-d_H}(H,R)$ for any cofinite subgroup $H
\subset G$. If this happens, $\bPhi^H(M)$ lies in $\D^{\leq n -
  d_H}(R)$ since $\Phi^H$ is left-exact with respect to the standard
$t$-structures. Conversely, if $\bPhi^{H'}(M) \in \D^{\leq n -
  d_{H'}}(R)$ for any cofinite $H' \subset H$, then
$\Phi(\Psi^H(M))$ lies in $\DN^{\leq n-d_H}(H,R)$ by
Lemma~\ref{phi.dot.le}~\thetag{i}, and then $\Psi^H(M) \cong
\Infl(\Phi(\Psi^H(M)))$ lies in $\wDM^{\leq n-d_H}(H,R)$ by
Proposition~\ref{can.filt.der.prop}.
\endproof

I do not know whether for a general non-descreasing sequence
$\{d_\idot\}$, the subcategories \eqref{d.idot} for all integers $n$
define a $t$-structure on $\wDM(G,R)$. However, note the
following. The notion of a homological sphere of
Definition~\ref{sph.def} extends literally to infinite groups: a
simplicial pointed admissible $G$-set $X \in \Delta^o\wGamma_{G+}$
is a homological sphere if $\overline{C}_\idot(X^H,\Z) \cong
\Z[d_H]$ for any cofinite subgroup $H \subset G$. Moreover, if one
modifies Example~\ref{regu.sph.exa} by letting
\begin{equation}\label{regu.eq}
X(G) = \lim_{\overset{H \subset G}{\to}}X(G/H),
\end{equation}
where the limit is taken over all cofinite subgroups $H \subset G$,
then $X(G)$ is a homological sphere, with $d_H = |G/H|$.

\begin{lemma}\label{pro.equi.le}
Assume that the group $G$ is finitely generated. Assume given a
homological sphere $X\in \Delta^o\wGamma_{G+}$, with the sequence of
degrees $\{d_\idot\}$. Then for any integer $n$, the functor $E_X$,
$M \mapsto M \wedge X$ is an equivalence of categories between
$\wDM^{\leq n}(G,R)$ and $\wDM^{\leq n - d_\idot}(G,R)$.
\end{lemma}

\proof{} If the group is finite, this is
Proposition~\ref{equi.prop}. In the general case, apply
Proposition~\ref{can.filt.der.prop} and Lemma~\ref{d.le}.
\endproof

As a consequence of this, the subcategories \eqref{d.idot} do form a
$t$-structure on $\wDM(G,R)$ if the sequence $\{d_\idot\}$
corresponds to a homological sphere (for example, one can take $d_H
= |G/H|$, with the sphere \eqref{regu.eq}). However, in this case,
the subcategory
$$
\bigcup_n \wDM^{\leq n - d_\idot}(G,R) \subset \wDM^-(G,R) \subset \wDM(G,R)
$$
equivalent to $\wDM^-(G,R)$ need not coincide with the whole
$\wDM^-(G,R)$. I do not know whether the functor $E_X$ is an
autoequivalence of the whole $\wDM(G,R)$.

\medskip

Another consequence of Lemma~\ref{add.der.pro.corr} is the following
result that turns out to be very useful for constructing Mackey
profunctors.

\begin{defn}\label{l.f.s}
An $R$-valued derived Mackey functor $M \in \DM(G,R)$ is {\em
  locally finitely supported} if $M$ lies in $\DM^{\leq n}(G,R)$ for
some integer $n$, and for any integer $m$, $M([G/H])$ lies in
$\D^{\leq m}(R)$ for all but a finite number of cofinite subgroups
$H \subset G$.
\end{defn}

\begin{lemma}
Assume that the group $G$ is finitely generated. Then for any
integer $n$ and any ring $R$, the forgetful functor
$$
\wDM^{\leq n}(G,R) \to \DM^{\leq n}(G,R)
$$
induced by the forgetful functor \eqref{d.M.wM.eq} has a
left-adjoint completion functor
$$
\DM^{\leq n}(G,R) \to \wDM^{\leq n}(G,R).
$$
Moreover, if $M \in \DM(G,R)$ is locally finitely supported in the
sense of Definition~\ref{l.f.s}, then the adjunction map $M \to
\wh{M}$ from $M$ to its complection $\wh{M}$ is an isomorphism.
\end{lemma}

In other words, for finitely generated groups, locally finitely
supported derived Mackey functors canonically extend to derived
Mackey profunctors.

\proof{} The completion functor is the composition $\Add \circ S(i)_!$,
where $i:\Gamma_G \to \wGamma_G$ is the natural
embedding. Explicitly, for any $M \in \DM^{\leq n}(G,R)$, its
completion $\wh{M}$ is given by
$$
\wh{M} = \dlim_{\overset{N}{\gets}}\Infl^N\Phi^N M,
$$
where the derived inverse limit in the right-hand side is taken over
all cofinite normal subgroups $N \subset G$, and $\Phi^N$ is the
geometric fixed points functor $S(\phi^N)_!:\DM(G,R) \to \DM(W,R)$,
$W = G/N$. If $M$ is finitely supported -- that is, $M([G/H]) = 0$
for all but a finite number of cofinite subgroups $H \subset G$ --
then the inverse system stabilizes at a finite step, and we have $M
\cong S(i)^*\wh{M}$, as required. To extend the argument to a
locally finitely supported derived Mackey functor $M$, note that the
completion functor is left-exact with respect to the standard
$t$-structure, and that by Definition~\ref{l.f.s}, for any integer
$m$, the truncation $\tau_{\geq m}M$ with respect to the standard
$t$-structure is a finitely supported derived Mackey functor.
\endproof

\subsection{Computations.}

To finish the section, we will now describe explicitly the derived
versions of the profunctors $\wBG(S,-)$ of Example~\ref{wB.G.exa}
(and in particular, we will compute the derived Burnside ring
$\wAG_\idot$).

For any admissible $G$-set $S$, consider the derived Mackey
profunctor
$$
\Add(\Sp(\Z_S)) \in \wDM^{\leq 0}(G,\Z),
$$
where $\Z_S \in \Fun(S\wGamma_G,\Z)$ is the functor represented by
$S \in \wGamma_G \subset S\wGamma_G$, $\Sp$ is the specializaion
functor of Lemma~\ref{sp.le}, and $\Add$ is the additivzation
functor of Lemma~\ref{add.der.pro.corr}. For any two admissible
$G$-sets $S,S' \in \wGamma_G$, let
\begin{equation}\label{add.bg}
\wBG_\idot(S,S') = \Add(\Sp(\Z_S))(S').
\end{equation}
Note that by adjunction, we have
$$
\wBG_\idot(S,S') = \Hom^\hdot(\Add(\Sp(\Z_S)),\Add(\Sp(\Z_{S'}))),
$$
so that we have natural associative product maps
\begin{equation}\label{B.G.prod}
\wBG_\idot(S_1,S_2) \otimes \wBG_\idot(S_2,S_3) \to
\wBG_\idot(S_1,S_3).
\end{equation}
These maps can be refined so that $\wBG_\idot(-,-)$ becomes an
$A_\infty$-category, but we will not need this; we refer an
interested reader to \cite[Subsection 1.6]{mackey}. However, we note
that by adjunction, we have a natural map
\begin{equation}\label{wB.aug}
\wBG_\idot(S_1,S_2) \to \wBG(S_1,S_2)
\end{equation}
for any $S_1,S_2 \in \wGamma_G$, where the right-hand side is as in
\eqref{wB.G}, and these maps are multiplicative with respect to the
product \eqref{B.G.prod}. In particular, setting $S_1 = S_2 = \ppt$,
we obtain a natural map of DG rings
\begin{equation}\label{A.aug.eq}
\wAG_\idot \to \wAG,
\end{equation}
where $\wAG$ is the completed Burnside ring of \eqref{wA.G.eq}.

\begin{prop}\label{wAG.prop}
For any finitely generated group $G$ and any admissible $G$-sets
$S,S' \in \wGamma_G$, we have
\begin{equation}\label{wA.G.der.eq}
\wBG_\idot(S,S') \cong \prod_{H \subset G}C_\idot(W_H,\Z[(S \times
  S')^H]),
\end{equation}
where the product is over all the conjugacy classes of cofinite
subgroups $H \subset G$, $W_H = \Aut_G([G/H])$ is the quotient
$N_H/H$, as in \eqref{wB.G}, and and $C_\idot(W_H,-)$ is the
homology complex of the group $W_H$ with coefficients in the free
abelian group $\Z[(S_1 \times S_2)^H]$ spanned by the fixed point
set $(S_1 \times S_2)^H$. In particular, $\wBG_\idot(S_1,S_2)$ lies
in $\D^{\leq 0}(\Z)$, and the map \eqref{wB.aug} identifies its
truncation at $0$ with $\wBG(S_1,S_2)$ of \eqref{wB.G}.
\end{prop}

\proof{} The claim is clearly compatible with the functors $\Phi$
and $\Infl$ of Proposition~\ref{can.filt.der.prop}, so it suffices
to consider the case when the group $G$ is finite. Then the
addivization functor $\Add$ in \eqref{add.bg} is the functor of
Proposition~\ref{add.prop}. Thus to evaluate the right-hand side of
\eqref{add.bg}, we need to restrict $\Sp(\Z_S)$ via the map
$m_{S'}:\Gamma \to \Gamma_G$ and then apply the functor $L^T$. By
Lemma~\ref{sp.mc}, we have
$$
\Sp(\Z_S) = L^\hdot\rho_{S!}\Z,
$$
where $\rho_S$ is the natural cofibration \eqref{rho.c.eq} with
$c=S$. Thus if denote by $\wt{\rho}:\wt{\Q} \to \Gamma_+$ the fibration
obtained by the cartesian square
$$
\begin{CD}
\wt{\Q} @>>> \Q^S\Gamma_G\\
@V{\wt{\rho}}VV @V{\rho_S}VV\\
S_{\Inj}\Gamma @>{S(m_{S'})}>> S\Gamma_G,
\end{CD}
$$
then $\wt{\Q}$ descends to a fibration $\rho:\Q \to \Gamma_+ \cong
Q_{\Inj}\Gamma$, in the sense that there is a cartesian square
$$
\begin{CD}
\wt{\Q} @>>> \Q\\
@V{\wt{\rho}}VV @VV{\rho}V\\
S_{\Inj}\Gamma @>{q}>> \Gamma_+,
\end{CD}
$$
and we have a natural identification
$$
\wBG_\idot(S,S') \cong L^TL^\hdot\rho_!\Z.
$$
Explicitly, the category $\Q$ is the category of triples $\langle
S_1,S_2,f \rangle$ of a finite set $S \in \Gamma$, a finite $G$-set
$S_2 \in \Gamma_G$, and a map $f:S_2 \to S_1 \times S \times
S'$. Morphisms from $\langle S_1,S_2,f \rangle$ to $\langle
S'_1,S'_2,f' \rangle$ are given by pair of a map $S_1 \to S_1'$ in
$\Gamma_+ \cong Q_{\Inj}\Gamma$ represented by a diagram
\begin{equation}\label{g.pl}
\begin{CD}
S_1 @<{i}<< \wt{S}_1 @>{g_1}>> S_1'
\end{CD}
\end{equation}
in $\Gamma$ with injective $i$, and a map $g_2:S_2' \to S_2$ in
$\Gamma_G$ that fit into a commutative diagram
$$
\begin{CD}
S_2 @<{g_2}<< S'_2 @= S'_2\\
@V{f}VV @VVV @VV{f'}V\\
S_1 \times S \times S' @<{i \times \id \times \id}<< \wt{S} \times S
\times S' @>{g_1 \times \id \times \id}>> S_1' \times S \times S'
\end{CD}
$$
in $\Gamma_G$ with cartesian square on the left. The projection
$\rho$ sends $\langle S_1,S_2,f \rangle$ to $S_1 \in
\Gamma$.

Moreover, let $\R$ be the category of finite $G$-sets
$\overline{S}$ equipped with a morphism $\overline{f}:\overline{S}
\to S \times S'$, with maps from $\overline{f}:\overline{S} \to S
\times S'$ to $\overline{f}':\overline{S}' \to S \times S'$ given by
injective maps $g:\overline{S}' \to \overline{S}$ such that
$\overline{f} \circ g = \overline{f}'$. Then for any $\langle
S_1,S_2,f \rangle \in \Q$, we can compose $f$ with the natural
projection $S_1 \times S \times S' \to S \times S'$ to obtain a map
$\overline{f}:S_2 \to S \times S'$, and sending $\langle S_1,S_2,f
\rangle$ to $\overline{f}$ and sending a map $\langle i,g_1,g_2
\rangle$ to $g_2$ gives a functor $\Q \to \R$.  This functor has a
left-adjoint $\nu:\R \to \Q$ that sends a finite $G$-set
$\overline{S}$ with a map $\overline{f}:\overline{S} \to S \times
S'$ to the triple $\langle \overline{S}/G,\overline{S},q \times
\overline{f}\rangle$, where $q:\overline{S} \to \overline{S}/G$ is
the quotient map. Therefore $L^\hdot\nu_!\Z = \Z$, and we have
$$
\wBG_\idot(S,S') \cong L^T(L^\hdot\tau_!\Z),
$$
where $\tau = \rho \circ \nu:\R \to \Gamma_+$ sends $\overline{S}$
with a map $\overline{f}:\overline{S} \to S \times S'$ to the
quotient $\overline{S}/G$. By definition, $L^T \circ L^\hdot\tau_!$
is left-adjoint to the composition $\tau^* \circ j^T$ where $j^T$ is
the functor of \eqref{j.T.spl.eq}. Explicitly, this composition
sends an abelian group $M$ to $M \otimes \tau^*T \in
\Fun(\R,\Z)$. But the functor $\tau:\R \to \Gamma_+$ factors through
the subcategory $\Gamma_{\Inj}^o \subset \Gamma_+$ of finite sets
and maps between them represented by diagrams \eqref{g.pl} with
invertible $g_1$, and the functor $\tau:\R \to \Gamma_{\Inj}^o$ is a
cofibration. If we let $O = \tau^{-1}([1])$ be the fiber of this
cofibration over $[1] \in \Gamma^o_{\Inj}$, then the cartesian square
$$
\begin{CD}
O @>{\eta}>> \R\\
@VVV @VV{\tau}V\\
\ppt @>>> \Gamma_{\Inj}^o
\end{CD}
$$
induces a base change isomorphism $\tau^*(j^T(M)) = M \otimes
\tau^*T \cong L^\hdot\eta_*M$, so that by adjunction, we have
$$
L^T(L^\hdot\tau_!\Z) \cong C_\idot(O,\Z).
$$
It remains to notice that by definition, $O$ is the groupoid of all
$G$-orbits $[G/H]$ equipped with a map $[G/H] \to S \times S'$, and
giving such a map is equivalent to giving and element $s \in (S
\times S')^H$. Therefore $C_\idot(O,\Z)$ is exactly the same as the
right-hand side of \eqref{wA.G.der.eq}. This finishes the proof.
\endproof

In particular, we see that the derived Burnside ring $\wAG_\idot$ is
given by
$$
\wAG_\idot \cong \prod_{H \subset G}C_\idot(W_H,\Z),
$$
where as in \eqref{wA.G.der.eq} and \eqref{wA.G.eq}, the product is
over all conjugacy classes of cofinite subgroups $H \subset
G$. Comparing this to \eqref{wA.G.eq}, we see that in degree $0$,
the homology of $\wAG_\idot$ is exactly the underived completed
Burnside ring $\wAG$.

\section{The cyclic group case.}\label{mack.cycl.subs}

We finish the paper by showing how all the abstract machinery that
we have developed works in the particular case $G=\Z$, the infinite
cyclic group.

\subsection{Overview.}

All cofinite subgroups in $\Z$ are of the form $l\Z \subset \Z$, $l
\geq 1$. To simplify notation, denote
\begin{equation}\label{M.l.eq}
M_l = M([\Z/l\Z])
\end{equation}
for any $\Z$-Mackey functor $M \in \M(\Z,R)$. For any $l \geq 1$,
$M_l$ carries a natural action of the quotient group $\Z/l\Z$, and
we denote the generator of this group by $\sigma$. For any $l,l'
\geq 1$, the quotient map $q:\Z/ll'\Z \to \Z/l\Z$ induces
canonical maps
\begin{equation}\label{v.f.mack}
v_{l',l} = q_*:M_{ll'} \to M_{l'}, \qquad
f_{l',l} = q^*:M_{l'} \to M_{ll'}.
\end{equation}

\begin{lemma}
  The category $\M(\Z,R)$ is equivalent to the category of
  collections of $R[\Z/l\Z]$-modules $M_l$, $l \geq 1$, equipped
  with maps \eqref{v.f.mack} satisfying
\begin{equation}\label{z.mack.exp}
\begin{aligned}
v_{l',l} &= v_{l',l} \circ \sigma^l,\\
f_{l',l} &= \sigma^l \circ f_{l',l},\\
v_{l'',l} &= v_{l',l} \circ v_{l'',l'},\\
f_{l'',l} &= f_{l'',l'} \circ f_{l',l},\\
f_{l'',l} \circ v_{l',l} &= v_{l'l''/l,l''} \circ f_{l'l''/l,l'}
\quad \text{ if $l'/l$ and $l''/l$ are coprime}\\
t_{l',l} &= f_{l',l} \circ v_{l',l},
\end{aligned}
\end{equation}
where $t_{l',l} = 1 + \sigma^l + \sigma^{2l} + \dots +
\sigma^{(l'-1)l}$ is the averaging over the subgroup $\Z/l'\Z
\subset \Z/ll'\Z$.
\end{lemma}

\proof{} The first four equations represent the functoriality of the
maps $f_*$, $f^*$ of Subsection~\ref{ma.def.subs}, and the last two
equations express the double coset formula \eqref{double}.
\endproof

A $\Z$-Mackey profunctor is in particular a $\Z$-Mackey functor, so
that we will still use the notation $M_l$, and we still have the
canonical maps \eqref{v.f.mack}. For any $l \geq 1$, we denote
\begin{equation}\label{psi.l}
\begin{aligned}
\gamma_l &= \gamma^{l\Z}:\wGamma_\Z \to \wGamma_\Z,\\
\Psi^l &= \Psi^{l\Z}:\wM(\Z,R) \to \wM(\Z,R).
\end{aligned}
\end{equation}
Of course, $l\Z \subset \Z$ is canonically isomorphic to $\Z$ as an
abstract group, so that $\Psi^l$ is an endofunctor of the category
$\wM(\Z,R)$. For any $l,l' \geq 1$, we obviously have $\Psi^l \circ
\Psi^{l'} = \Psi^{ll'}$.

The completed Burnside ring $\wAZ$ is easy to compute. As a group,
it is given by
\begin{equation}\label{wAZ.eq}
\wAZ = \Z\langle \eps_1,\eps_2,\dots \rangle,
\end{equation}
the group of infinite linear combinations of generators $\eps_i$
numbered by all integers $i \geq 1$. The generator $\eps_i$
corresponds to the $\Z$-orbit $\Z/i\Z$. The product in $\wAZ$ is
given by
\begin{equation}\label{a.z.rel}
\eps_i\eps_j = \frac{ij}{\{i,j\}}\eps_{\{i,j\}},
\end{equation}
where $\{i,j\}$ stands for the least common multiple of $i$ and
$j$. In particular, we have
$$
\eps_i^2=i\eps_i.
$$
The element $\eps_1 \in \wAZ$ is the unit of the ring $\wAZ$. The
Burnside ring acts on any $M \in \wM(\Z,R)$, and the action of the
generators $\eps_i$ is given by \eqref{burn.act.eq}. In particular,
for any $l,l' \geq 1$ and any $M \in \wM(\Z,R)$, we have
\begin{equation}\label{eps.acts.eq}
\eps_{ll'} = v_{l,l'} \circ f_{l,l'}:M_l \to M_l,
\end{equation}
where $v_{l,l'}$ and $f_{l,l'}$ are the natural maps
\eqref{v.f.mack}.

\begin{remark}
  A reader might notice that the ring $\wAZ$ in fact coincides with
  the universal Witt vectors ring $\W(\Z)$. We will explore this
  coincidence elsewhere.
\end{remark}

All subgroups in $\Z$ are normal, so that for any separated $M \in
\wM_s(\Z,R)$, the terms in the canonical filtration
\eqref{can.filt.mack} are numbered by positive integers $l \geq 1$,
and we have $F^{l'\Z}M \subset F^{l\Z}M$ whenever $l'$ is divisible
by $l$. It is convenient to consider a slightly different filtration
on $M$ by letting
\begin{equation}\label{can.filt.Z}
F^lM = \cap_{n \leq l}F^{n\Z}M \subset M
\end{equation}
for any $l \geq 1$. Then $F^{l'}M \subset F^lM$ whenever $l' \geq
l$, $F^{l\Z}M \supset F^lM \supset F^{l!\Z}M$, and $M$ is
automatically complete with respect to the filtration $F^lM$.

By definition, for any group $G$, we have $\wGamma_G =
\wGamma_{\wh{G}}$, where $\wh{G}$ is the profinite completion of the
group $G$. In particular, $\M(\Z,R) = \M(\wh{\Z},R)$. Note that
$$
\wh{Z} = \prod_p \Z_p,
$$
where the product is over all primes $p$. For every prime $p$,
let
\begin{equation}\label{z.p.pr}
\Z'_p = \prod_{l \neq p}\Z_l,
\end{equation}
where the product is over all primes distinct from $p$. Then we have
$\wh{\Z} = \Z_p \times \Z'_p$. Denote by
\begin{equation}\label{phi.p}
\Phi^{(p)} = \wh{\Phi}^{\Z'_p}:\wM(\Z,R) = \wM(\wh{\Z},R) \to
\wM(\Z_p,R)
\end{equation}
the geometric fixed points functor with respect to the embedding
$\Z'_p \subset \wh{\Z}$, and let $M^{(p)} = \Phi^{(p)}(M)$ for any
$M \in \wM(\Z,R)$.

\begin{defn}\label{p.typ.def}
The $\Z_p$-Mackey profunctor $M^{(p)} \in \wM(\Z_p,R)$ is called the
{\em $p$-typical part} of the $\Z$-Mackey profunctor $M \in
\wM(\Z,R)$.
\end{defn}

Cofinite subgroups in $\Z_p$ are of the form $\Z_p \subset \Z_p$, $q
= p^{n-1}$, $n \geq 1$; as before, we simplify notation and denote
$$
M_n = M(\Z_p/q\Z_p)
$$
and $\Psi^n = \Psi^{q\Z_p}$ for any $M \in \wM(\Z_p,R)$ and any such
$q$. Since $q\Z_p$ is abstractly isomorphic to $\Z_p$, we can treat
$\Psi^n$ as an endofunctor of the category $\M(\Z_p,R)$. We further
simplify notation by setting $\Psi = \Psi^1$. We have natural
isomorphisms $\Psi^n \cong (\Psi^1)^n$, so that the notation is
consistent. If $M$ is separated, we also renumber the canonical
filtration \eqref{can.filt.mack} on $M$ by setting
\begin{equation}\label{can.filt.Zp}
F^nM = F^{q\Z_p}M
\end{equation}
for any $n \geq 1$ and $q=p^n$.

\subsection{The $p$-typical decomposition.}\label{p.typ.subs}

It turns out that under some assumptions, a $\Z$-Mackey profunctor
$M$ is completely determined by the $p$-typical parts of the
functors $\Psi^l(M)$, $l \geq 1$ not divisible by $p$. 

Namely, fix a prime $p$, and assume that the base ring $R$ is
$p$-local -- that it, all integers $i \geq 1$ coprime to $p$ are
invertible in $R$. Then the completed Burnside ring $\wAZ$ is given
by \eqref{wAZ.eq}, and in particular, $\frac{1}{i}\eps_i$ with $i$
coprime to $p$ is a well-defined idempotent in $R \otimes \wAZ$.
Any $\Z$-Mackey functor $M \in \M(\Z,R)$ comes equipped with an
action of $\wAZ \otimes R$, thus with commuting idempotent
endomorphisms $\frac{1}{i}\eps_i$. Moreover, for any $l$ prime to
$p$, $\wAZ \otimes R$ contains a well-defined idempotent element
\begin{equation}\label{eps.l.dual}
\eps_{(l)} = \frac{1}{l}\eps_l \cdot \prod_{i
  \text{ does not divide }l}\left(1 - \frac{1}{i}\eps_i\right),
\end{equation}
where the product is over $i$ prime to $p$. We have
$$
1 = \sum_{l \text{ prime to }p}\eps_{(l)},
$$
so that for any $M \in \wM(\Z,R)$, we have a canonical decomposition
\begin{equation}\label{typ.deco}
M = \prod_{l \text{ prime to } p}M_{(l)},
\end{equation}
where $M_{(l)}$ is the image of the idempotent $\eps_{(l)}$. We will
say that $M \in \wM(\Z,R)$ is {\em type $l$} if $M = M_{(l)}$, and
we will denote by $\wM_l(\Z,R)$ the full subcategory spanned by
Mackey profunctors of type $l$. We then have the natural
decomposition
\begin{equation}\label{type.deco}
\wM(\Z,R) \cong \prod_{l \text{ prime to } p}\wM_l(\Z,R).
\end{equation}
Now note that for any $l \geq 1$ not divisible by $p$,
\eqref{phi.psi.iso} with $H = \Z'_p \subset \wh{\Z}$ and $H' =
l\Z'_p \subset H$ gives an isomorphism
$$
(\Psi^lM)^{(p)} = \Phi^{(p)}\Psi^lM \cong \Psi^{\Z_p}\Phi^{l\Z'_p}M,
$$
and $\Psi^{\Z_p}$ in the right-hand side can be promoted to the
functor $\wt{\Psi}^{\Z_p}$ of \eqref{psi.wt}. In particular,
$(\Psi^lM)^{(p)}$ carries a natural action of the cyclic group
$\Z/l\Z \subset Z_{\Z_p} \subset \Z_p \times \Z/l\Z =
\wh{\Z}/l\Z'_p$, so that we can promote the functor
$\Phi^{(p)} \circ \Psi^l$ to a functor
\begin{equation}\label{phi.p.l}
\Phi^{(p)}_{(l)}:\wM(\Z,R) \to \wM(\Z_p,R[\Z/l\Z]).
\end{equation}

\begin{prop}\label{p-typ.prop}
Assume that $R$ is $p$-local. Then for any $l$ prime to $p$, the
functor $\Phi^{(p)}_{(l)}$ of \eqref{phi.p.l} induces an equivalence
\begin{equation}\label{p.l.equi}
\Phi^{(p)}_{(l)}:\wM_l(\Z,R) \cong \wM(\Z_p,R[\Z/l\Z]),
\end{equation}
where $\wM_l(\Z,R) \subset \wM(\Z,R)$ is the term of the
decompositon \eqref{type.deco} corresponding to $l$.
\end{prop}

\proof{} Since $R$ is $p$-local, for any $l \geq 1$ not divisible by
$p$, the functor
$$
R_l = R_{\{e\}}:\D(R[\Z/l\Z]) \to \DM(\Z/l\Z,R)
$$
provided by Lemma~\ref{R.H.le} is exact with respect to the standard
$t$-structure, thus induces a functor from $R[\Z/l\Z]$-modules to
$\M(\Z/l\Z,R)$ right-adjoint to the restriction functor
$\Psi^{\{e\}}$. Explicitly, we have
\begin{equation}\label{coinv}
R_l(M)_n = M_{\sigma^{l/n}} \cong M^{\sigma^{l/n}}
\end{equation}
for any $R[\Z/l\Z]$-module $M$ and divisor $n$ of the integer $l$,
where $\sigma:M \to M$ generates the action of $\Z/l\Z$. The functor
is fully faithful. Applying this functor pointwise over $Q\Z_p$, we
obtain a functor
$$
\wt{R}_l:\wM(\Z_p,R[\Z/l\Z]) \to \wM(\Z_p \times (\Z/l\Z),R)
$$
right-adjoint to the functor $\wt{\Psi}^l$, and moreover, we have
$\wt{R}_l \circ \wt{\Psi}^l \cong \Id$, so that $\wt{R}_l$ is fully
faithful. Composing it with the inflation functor $\Infl^{l\Z'_p}$,
we obtain a fully faithful embedding
$$
\overline{R}_l:\wM(\Z_p,R[\Z/l\Z]) \to \wM(\Z,R)
$$
right-adjoint to the functor $\Phi^{(p)}_{(l)}$ of
\eqref{phi.p.l}. By \eqref{coinv} and \eqref{eps.acts.eq}, for any
$E \in \wM(\Z_p,R[\Z/l\Z])$ and any $n \geq 1$ prime to $p$, the
generator $\eps_n$ of the completed Burnside ring $\wAZ$ acts on
$\overline{R}_l(E)$ by $n\id$ if $n$ divides $l$ and by $0$
otherwise. By \eqref{eps.l.dual}, this means that
$\overline{R}_l(E)$ is of type $l$, that is, the functor
$\overline{R}_l$ factors through a fully faithful embedding
$$
\overline{R}_l:\wM(\Z_p,R[\Z/l\Z]) \to \wM_l(\Z,R).
$$
This embedding is then automatically right-adjoint to the functor
\eqref{p.l.equi}, so to finish the proof, it suffices to check that
$\overline{R}_l$ is essentially surjective. By \eqref{type.deco}, we
may as well prove that the fully faithful functor
$$
\overline{R}_\idot = \prod_{l \text{ prime to }
  p}\overline{R}_l:\prod_{l \text{ prime to } p}\wM(\Z_p,R[\Z/l\Z])
\to \wM_l(\Z,R)
$$
is essentially surjective. By construction, this functor is
right-adjoint to
$$
\Phi^{(p)}_\idot = \prod\Phi^{(p)}_{(l)}:\wM(\Z,R) \to \prod_{l
  \text{ prime to } p}\wM(\Z_p,R[\Z/l\Z]),
$$
so it suffices to prove that the adjunction map $a:\Id \to
\overline{R}_\idot \circ \Phi^{(p)}_\idot$ is an
isomorphism. Moreover, since $\overline{R}_\idot$ is fully faithful,
we know that $\Phi^{(p)}_\idot(a) = \id$, so it suffices to check
that the functor $\Phi^{(p)}_\idot$ is conservative. However, we can
extend it to a functor
\begin{equation}\label{phi.p.l.der}
\Phi^{(p)}_\idot:\wDM^{\leq 0}(\Z,R) \to \prod_{l
  \text{ prime to } p}\wDM^{\leq 0}(\Z_p,R[\Z/l\Z]),
\end{equation}
with the same definition, and note that since $R$ is $p$-local, the
extended functor is exact with respect to the standard
$t$-structures. Therefore it suffices to check that the extended
functor $\Phi^{(p)}_\idot$ is conservative. This is clear:
Proposition~\ref{can.filt.der.prop} and
Lemma~\ref{phi.dot.le}~\thetag{ii} show that the functor
$$
\bPhi_\idot = \prod_{n \geq 1}\bPhi^{n\Z}:\wDM^{\leq 0}(\Z,R) \to
\prod_{n \geq 1}\D^{\leq 0}(R)
$$
is conservative, and $\bPhi_\idot$ factors through
$\Phi^{(p)}_\idot$.
\endproof

We note that \eqref{coinv} allows to describe the decomposition
\ref{typ.deco} rather explicitly. In particular, for any $M \in
\wM(\Z,R)$, and any integer $n \geq 1$ of the form $n=lp^m$, $l$
prime to $p$, we have a natural decomposition
\begin{equation}\label{typ.exp}
M_n \cong \prod_{l' \text{ prime to }
  p}\left(\Phi^{(p)}_{(ll'}(M)_m\right)_{\sigma^l},
\end{equation}
where $\sigma$ in the right-hand side is the generator of the group
$\Z/ll'\Z$, and it acts on $\Phi^{(p)}_{(ll')}(M) \in
\wM(\Z_p,R[\Z/ll'\Z])$ via the group algebra $R[\Z/ll'\Z]$.

\begin{lemma}
  Assume that $R$ is $p$-local, and assume that a $\Z$-Mackey
  profunctor $M \in \wM(\Z,R)$ is separated in the sense of
  Definition~\ref{sepa.def}. Then the canonical filtration
  \eqref{can.filt.Z} on $M$ induces the canonical filtration
  \eqref{can.filt.Zp} on $M^{(p)} = M_{(1)}$ by
$$
(F^mM^{(p)}) \cong (F^{p^m}M)_{(1)},
$$
and we have
$$
(\gr^n_FM)_{(1)} =
\begin{cases}
\gr^m_FM^{(p)}, &m=p^n,\\
0, &\text{otherwise}.
\end{cases}
$$
\end{lemma}

\proof{} Immediately follows from \eqref{typ.exp}.
\endproof

\subsection{Fixed points data.}

On the level of derived Mackey profunctors, we still have an
identification $\wDM(\Z,R) \cong \wDM(\wh{\Z},R)$, so that taking
the $p$-typical part of a derived Mackey profunctor makes sense (at
least for derived Mackey profunctors bounded from above). Moreover,
by Proposition~\ref{wAG.prop}, the elements $\eps_l \in \wAZ$, $l$
prime to $p$ lift to elements in the degree-$0$ homology of the
derived completed Burnside ring $\wAZ_\idot$, thus act on any object
$E \in \wDM(\Z,R)$ by \eqref{wAG.acts}. Thus it is reasonable to
expect that Proposition~\ref{p-typ.prop} has a derived version.

However, it turns out that if we restrict our attention to
$\wDM^-(\Z,R) \subset \wDM(\Z,R)$, we can do much better: there is
an alternative model of the category $\wDM^-(\Z,R)$ that is
reasonably simple, requires no assumptions on the ring $R$, and in
the $p$-local case, makes the $p$-typical decomposition of
Proposition~\ref{p-typ.prop} completely obvious.

To describe this model, denote by $I$ the groupoid of all finite
$\Z$-orbits and isomorphisms between them. We have a decomposition
$$
I = \coprod_{n \geq 1}\ppt_n,
$$
where $\ppt_n$ is a groupoid with one object with automoprhism group
$\Z/n\Z$. For any $M \in \D(I,R)$ and $n \geq 1$, we denote by
$$
M_n \in \D(\ppt_n,R) \cong \D(R[\Z/n\Z])
$$
its restriction to $\ppt_n \subset I$. For any prime $p$, let $I_p
\subset I$ be the full subcategory given by
\begin{equation}\label{i.p.eq}
I_p = \prod_{n\geq 1}\ppt_{np} \subset I,
\end{equation}
and let $I_\idot$ be the disjoint union of the categories $I_p$ over
all primes. Denote by
\begin{equation}\label{i.eq}
i:I_\idot \to I
\end{equation}
the functor that acts on $I_p$ by the embedding \eqref{i.p.eq}, and
let
\begin{equation}\label{pi.eq}
\pi:I_\idot \to I
\end{equation}
be the functor that acts on $I_p$ by the union of the natural
projections $\ppt_{np} \to \ppt_n$ corresponding to the quotient
maps $\Z/pn\Z \to \Z/n\Z$. For any ring $R$, any object $M \in
\Fun(I,R)$, and any integer $n \geq 1$ we denote by $M^n$ the
restriction of $M$ to $\ppt_n \subset I$, and for any prime $p$ and
an object $M \in \Fun(I_\idot,R)$, we denote by $M^{p,n}$ the
restriction of $M$ to $\ppt_{np} \subset I_p \subset I_\idot$.

\begin{defn}\label{ada.p.def}
  For any prime $p$ and integer $n \geq 1$, a $\Z[\Z/pn\Z]$-module
  $M$ is {\em $p$-adapted} if it gives a finitely generated
  projective $\Z[\Z/p\Z]$-module after restriction to $\Z/p\Z
  \subset \Z/pn\Z$. A functor $M \in \Fun(I_\idot,\Z)$ is {\em
    adapted} if for any $p$, $n$, its component $M^{p,n}$ is a
  $p$-adapted $\Z[\Z/np\Z]$-module. An {\em adapted complex} is an
  acyclic complex $P_\idot$ in $\Fun(I_\idot,\Z)$ such that $P_i =
  0$ for $i < 0$, $P_0$ is the constant functor $\Z$, and $P_i$ is
  adapted for $i > 0$.
\end{defn}

\begin{defn}
  An {\em $R$-valued fixed points datum} with respect to $P_\idot$
  is a pair $\langle M_\idot,\alpha \rangle$ of a complex $M_\idot$
  in $\Fun(I,R)$ and a map
\begin{equation}\label{alpha.eq}
\alpha:\pi^*M_\idot \to i^*M_\idot \otimes P_\idot,
\end{equation}
where $i,\pi:I_\idot \to I$ are projections \eqref{i.eq} and
\eqref{pi.eq}. A map of fixed points data from $\langle
M_\idot,\alpha \rangle$ to $\langle M'_\idot,\alpha' \rangle$ is a
map of complexes $f:M_\idot \to M_\idot'$ such that $\alpha' \circ
\pi^*f = (i^*f \otimes \id) \circ \alpha$, and the category of
$R$-valued fixed points data with respect to $P_\idot$ is denoted by
$C^\alpha(P_\idot,R)$.
\end{defn}

We say that an object $\langle M_\idot,\alpha \rangle$ in
$C^\alpha(P_\idot,R)$ is acyclic if the underlying complex $M_\idot$
is acyclic, and we say that a map in $C^\alpha(P_\idot,R)$ is a
quasiisomorphism if it is a quasiisomorphism of the underlying
complexes. Inverting quasiisomorphisms in $C^\alpha(P_\idot,R)$, one
obtains the {\em derived category of $R$-valued fixed points data}
that we denote by $\DA(P_\idot,R)$. We note that this localization
procedure presents no set-theoretical problems, by the same argument
as in \cite[Lemma 1.7]{mackey}.

We can also do the localization in two steps. First, we define the
category $H^\alpha(P_\idot,R)$ of $R$-valued fixed points data and
chain-homotopy classes of morphisms between them. This is a
triangulated category. Then we take its Verdier localization with
respect to the full subcategory spanned by acyclic complexes. This
shows that the category $\DA(P_\idot,R)$ is triangulated. The
obvious forgetful functor
\begin{equation}\label{phi.I}
\eta:\DA(P_\idot,R) \to \D(I,R), \quad \langle M_\idot,\alpha
\rangle \mapsto M_\idot
\end{equation}
is a triangulated functor. We note that it is conservative (that is,
for any morphism $f$ in $\DA(P_\idot,R)$, $\eta(f)$ is invertible if
and only if $f$ is invertible).

\begin{remark}
  Usually, another way to construct and study unbounded derived
  categories is to use model structures in the sense of
  Quillen. This does not seem to work for the category
  $C^\alpha(P_\idot,R)$. Indeed, it is a very simple example of the
  category of DG comodules over a DG coalgebra, and those are not
  known to possess reasonable model structures.
\end{remark}

\begin{lemma}\label{rho.le}
  The forgetful functor $\eta$ of \eqref{phi.I} has a right-adjoint
  functor
$$
\rho:\D(I,R) \to \DA(R).
$$
\end{lemma}

\proof{} Let $H(I,R)$ be the triangulated category of complexes in
$\Fun(I,R)$ and chain-homotopy classes of maps between them, so that
$\D(I,R)$ is the quotient of $H(I,R)$ by the subcategory of acyclic
complexes. Then the forgetful functor $\eta$ is induced by the
forgetful functor $\eta^h:H^\alpha(P_\idot,R) \to H(I,R)$, and this
has an obvious right-adjoint $\rho^h:H(I,R) \to H^\alpha(P_\idot,R)$
given by
\begin{equation}\label{rho.h}
\rho^h(M) = M \oplus \pi_*(i^*M \otimes P_\idot).
\end{equation}
Moreover, the quotient functor $q:H(I,R) \to \D(I,R)$
has a right-adjoint fully faithful functor $r:\D(I,R) \to H(I,R)$
sending $M  \in D(I,R)$ to its $h$-injective representative in
$H(I,R)$. Thus we are in the situation of Lemma~\ref{tria.le},
so that the functor
$$
\rho = q \circ \rho^h \circ r:\D(I,R) \to \DA(P_\idot,R)
$$
with $q:H^\alpha(P_\idot,R) \to \DA(P_\idot,R)$ being the
natural projection, is right-adjoint to $\eta$.
\endproof

Now denote by $R = \eta \circ \rho:\D(I,R) \to \D(I,R)$ the
composition of the forgetful functor $\eta$ of \eqref{phi.I} with
its right-adjoint $\rho$ provided by Lemma~\ref{rho.le}. Note that
by \eqref{rho.h}, the adjunction map $R \to \Id$ fits into a
functorial exact triangle
$$
\begin{CD}
\overline{R} @>>> R @>>> \Id @>>>,
\end{CD}
$$
where $\overline{R}:\D(I,R) \to \D(I,R)$ is given by
\begin{equation}\label{bar.R}
\overline{R}(M) = \pi_*(i^*r(M) \otimes P_\idot).
\end{equation}

\begin{lemma}\label{R.sq.le}
For any $M \in \D(I,R)$, we have
\begin{equation}\label{phi.tate}
  \overline{R}(M)^n \cong \bigoplus_{p \text{ a
      prime}}\vC^\hdot(\Z/p\Z,M^{np}),
\end{equation}
where $\overline{R}$ is the functor of \eqref{bar.R}, and
$\vC^\hdot(\Z/p\Z,-)$ in the right-hand side are the maximal Tate
cohomology objects \eqref{tate.c}. Moreover, we
have $\overline{R} \circ \overline{R} = 0$, and every object in
$\DA(P_\idot,R)$ is a cone of objects of the form $\rho(M)$, $M \in
\D(I,R)$.
\end{lemma}

\proof{} Note that for any $n \geq 1$, any $h$-injective complex $E$
of $R[\Z/n\Z]$-modules, and any finite-length complex $K_\idot$ of
$\Z[\Z/n\Z]$-modules that are finitely generated and flat over $\Z$,
the product $E \otimes K_\idot$ is $h$-injective. Therefore for each
term $F^lP_\idot$ of the stupid filtration on the adapted complex
$P_\idot$, the product $i^*r(M) \otimes F^lP_\idot$ is
$h$-injective. Since the functor $\pi_*$ commutes with filtered
direct limits, \eqref{bar.R} implies that
\begin{equation}\label{bar.R.lim}
\overline{R}(M) \cong \lim_{\overset{l}{\to}}R^\hdot\pi_*(i^*M
\otimes F^lP_\idot).
\end{equation}
For any $n$ and $p$, the complex $P^{p,n}_\idot$ considered
as a complex of $\Z[\Z/p\Z]$-modules is obviously maximally adapted
in the sense of Definition~\ref{tate.ada.def}, so that
\eqref{phi.tate} follows from \eqref{bar.R}. To prove that
$\overline{R} \circ \overline{R} = 0$, it then suffices to check
that
\begin{equation}\label{double.tate}
\vC^\hdot(\Z/p'\Z,\vC^\hdot(\Z/p\Z,M)) = 0
\end{equation}
for any two primes $p$, $p'$, and any $M \in \D(R[\Z/pp'\Z])$. For
any integer $n \geq 2$, the trivial representation $\Z$ of the
cyclic group $\Z/n\Z$ has a standard $2$-periodic projective
resolution $\wt{P}_\idot$. Using this resolution, one immediately
observes that the cohomology algebra $H^\hdot(\Z/n\Z,\Z)$ is the
algebra $\Z[u]$ with a single generator $u$ of degree $2$, subject
to a single relation $nu=0$. Thus for any $M \in \D(R[\Z/n\Z])$, we
have a natural map
\begin{equation}\label{u.eq}
u:M \to M[2].
\end{equation}
If $n=p$ is a prime, then the cone of the augmentation map
$\wt{P}_\idot \to \Z$ is maximally adapted in the sense of
Definition~\ref{tate.ada.def}. Using this complex to compute maximal
Tate cohomology, we conclude that for any $M \in \D(R[\Z/p\Z])$, we
have
\begin{equation}\label{t.l}
\vC^\hdot(\Z/p\Z,M) \cong
\dlim_{\overset{l}{\to}}C^\hdot(\Z/p\Z,M)[2l],
\end{equation}
where the limit is taken with respect to the map \eqref{u.eq}. This
shows that $vC^\hdot(\Z/p\Z,M)$ is $p$-local, and trivial if $p$ is
invertible in $M$. Thus \eqref{double.tate} trivially holds if $p
\neq p'$, and we may assume that $p=p'$. The natural map
$$
\pi:H^\hdot(\Z/p\Z,\Z) \to H^\hdot(\Z/p^2,\Z)
$$
induced by the quotient map $\Z/p^2\Z \to \Z/p\Z$ sends the
generator $u \in H^2(\Z/p\Z,\Z)$ to $pu \in H^2(\Z/p^2,\Z)$, so that
$\pi(u)^2=0$. Therefore for any $M \in \D(R[\Z/p\Z])$ of the form $M
= C^\hdot(\Z/p\Z,M')$, $M' \in \D(R[\Z/p^2\Z])$, the natural map
\eqref{u.eq} squares to $0$. The same is then true for arbitrary
sums of objects of this type. Then by \eqref{t.l}, we have
$$
\vC^\hdot\left(\Z/p\Z,\bigoplus_{l \geq 0}C^\hdot(\Z/p\Z,M)[2l]\right) =
0
$$
for any $M \in \D(R[\Z/p^2\Z])$. To deduce \eqref{double.tate},
evaluate its left-hand side by \eqref{t.l} and compute the limit by
the telescope construction \eqref{tel}.

Finally, for the last claim, take an object $E \in \DA(P_\idot,R)$,
and let $\overline{E}$ be the cone of the adjunction map $E \to
\rho(\eta(E))$. Then $\overline{R}(\eta(E)) = \overline{R} \circ
\overline{R}(E) = 0$, so that the adjunction map $\overline{E} \to
\rho(\eta(\overline{E}))$ becomes an isomorphism after applying
$\eta$. Since $\eta$ is conservative, $\overline{E} \cong
\rho(\eta(\overline{E}))$.
\endproof

\begin{lemma}\label{no.P.le}
The category $\DA(P_\idot,R)$ does not depend on the choice of the
adapted complex $P_\idot$.
\end{lemma}

\proof{} Any map $f:P_\idot \to P'_\idot$ between adapted complexes
$P_\idot$, $P'_\idot$ compatible with augmentations induces a
functor $f_*:C^\alpha(P_\idot,R) \to C^\alpha(P_\idot,R)$ given by
$$
f_*(\langle M_\idot,\alpha \rangle) = \langle M_\idot,(\id \otimes
f) \circ \alpha \rangle.
$$
This functor commutes with the forgetful functors $\eta$. Moreover,
by \eqref{phi.tate}, the base change map
$$
f_* \circ \rho \to \rho \circ f_*
$$
induced by the isomorphism $\eta \circ f_* \cong f_* \circ \eta$ is
also an isomorphism. Then by adjunction, $f_*$ is fully faithful on
the full triangulated subcategory in $\DA(P_\idot,R)$ spanned by
$\rho(\D(I,R))$, and by Lemma~\ref{R.sq.le}, this is the whole
category $\D(P_\idot,R)$. Moreover, its essential image contains the
full subcategory in $\DA(P'_\idot,R)$ spanned by $\rho(\D(I,R))$,
and by Lemma~\ref{R.sq.le}, this is the whole category
$\D(P'_\idot,R)$. Therefore $f_*$ is an equivalence of
categories. To finish the proof, it suffices to note that if
$\overline{P}_\idot$ is a projective resolution $P_\idot$ of the
constant functor $\Z \in \Fun(I,\Z)$, then the cone $P_\idot$ of the
augmentation map $\wt{P}_\idot \to \Z$ is obviously an adapted
complex, and for any other adapted complex $P'_\idot$, we have a map
$f:P_\idot \to P'_\idot$.
\endproof

Note that because of Lemma~\ref{no.P.le}, it is safe to drop
$P_\idot$ from notation and denote $\DA(R) = \DA(P_\idot,R)$.

Assume now given a subgroup $H \subset \wh{\Z}$, recall that the
groupoid $I$ is the groupoid of finite $\wh{\Z}$-orbits, and let
$I^H \subset I$ be the full subcategory spanned by orbits on which
$H$ acts trivially. Say that an object $M \in \D(I,R)$ is {\em
  supported at $H$} if $M(i) = 0$ for any $i \in I \setminus I^H$,
and denote by $\D_H(I,R) \subset \D(I,R)$ the full subcategory
spanned by objects supported at $H$. Then the subcategory $\D_H(I,R)
\subset \D(I,R)$ is canonically a direct summand, and the embedding
$\D_H(I,R) \subset \D(I,R)$ has an obvious left-adjoint functor
\begin{equation}\label{phi.h.I}
\phi^H:\D(I,R) \to \D_H(I,R).
\end{equation}
Explicitly, $\phi^H(M)(i) = M(i)$ if $i \in I$ lies in $I^H \subset
I$, and $0$ otherwise. Moreover, let $\DA_H(R) \subset
\DA(R)$ be the full subcategory spanned by objects $M$ such that
$\eta(M) \in \D(I,R)$ is supported at $H$. Then by
Lemma~\ref{R.sq.le}, the adjoint functor $\rho$ of
Lemma~\ref{rho.le} sends $\D_H(I,R) \subset \D(I,R)$ into $\DA_H(R)
\subset \DA(I,R)$, and the functor \eqref{phi.h.I} defines a functor
\begin{equation}\label{phi.h.a}
\phi^H:\DA(R) \to \DA_H(R)
\end{equation}
left-adjoint to the embedding $\DA_H(R) \subset \DA(R)$.

\subsection{Derived correspondence.}

Now we go back to $\wDM(\Z,R)$, the category of $R$-valued derived
$\Z$-Mackey profunctors.  We will use the same shorthand notation
\eqref{M.l.eq}, \eqref{psi.l}, \eqref{phi.p}, where $\Psi$ and
$\Phi$ now stand for fixed points functors on the derived
level. Moreover, for any $n \geq 1$, we denote
$$
\Phi^n = \Phi^{n\Z}:\wDM(\Z,R) \to \DM(\Z/n\Z,R),
$$
and we let
$$
\wPhi^n = \wt{\Psi}^{\{e\}} \circ \Phi^n:\wDM(\Z,R) \to
\D(R[\Z/n\Z]),
$$
where as in \eqref{wphi.H}, $\{e\} \subset \Z/n\Z$ is
the trivial group. Let
$$
R^n = \Infl^{n\Z} \circ R_{\{e\}}:\D(R[\Z/n\Z]) \to \wDM(\Z,R)
$$
be the right-adjoint functor to $\wPhi^n$, where $R_{\{e\}}$ is
as in Lemma~\ref{R.H.le}, and $\{e\} \subset \Z/n\Z$ is again the
trivial group. Note that for any $M \in \D(R[\Z/n\Z])$, $R^n(M)$ is
supported at $n\Z \subset \Z$.

\begin{lemma}\label{phi.tate.le}
For any $n,l \geq 1$ and $M \in \D(R[\Z/n\Z])$, we have
$$
\wPhi^l(R^n(M)) \cong
\begin{cases}
M, &\quad l=n,\\
\vC^\hdot(\Z/p\Z,M), &\quad n=pl, \quad p \text{ a prime},\\
0, &\quad \text{otherwise}.
\end{cases}
$$
\end{lemma}

\proof{} Since $R^n(M)$ is supported at $n\Z \subset \Z$,
$\wPhi^l(R^n(M))$ vanishes unless $l$ divides $n$, and in this case,
$\wPhi^l(R^n(M))$ is given by Proposition~\ref{tate.le}, with $G =
\Z/n\Z$. It remains to note that by \cite[Lemma
  7.15~\thetag{ii}]{mackey}, the maximal Tate cohomology
$\vC^\hdot(\Z/m\Z,M)$ with any coefficients vanishes unless $m$ is a
prime.
\endproof

Taken together, the functors $\wPhi^n$ define a functor
\begin{equation}\label{bar.P}
\wPhi:\wDM(\Z,R) \to \D(I,R)
\end{equation}
such that $\wPhi(M)^n = \wPhi^n(M)$, and the
functors $R^n$ together provide its right-adjoint
$$
\wt{R} = \prod_n R^n:\D(I,R) \to \wDM(\Z,R).
$$
Now for any $n \geq 1$, choose an admissible pointed $2$-simplicial
$\Z$-set $X_n$ that is adapted to $n\Z \subset \Z$ in the sense of
Definition~\ref{ada.N.def}. For any prime $p$, let $X_n^p =
X_n^{np\Z} \subset X_n$. Denote
\begin{equation}\label{p.n.p.eq}
P_\idot^{p,n} = \overline{C}_\idot(X^p_n/,\Z),
\end{equation}
and let $P_\idot$ be the complex in $\Fun(I_\idot,\Z)$ whose value
at $\ppt_{np} \subset I_p \subset I$ is given by the complex
$P_\idot^{p,n}$.

\begin{lemma}\label{P.le}
  The complex $P_\idot$ is adapted in the sense of
  Definition~\ref{ada.p.def}.
\end{lemma}

\proof{} For any $[m] \times [m'] \in \Delta^o \times \Delta^o$, if
an element $x \in X^p_n([m] \times [m'])$ is fixed under $\Z/p\Z
\subset \Z/pn\Z$, it is also fixed under $n\Z \subset \Z$, and then
it must lie in $\iota([1]_+)$ by
Definition~\ref{ada.N.def}~\thetag{i}. Therefore $\Z/p\Z \subset
\Z/np\Z$ acts freely on the quotient $X^p_n/[1]_+$, so that
$P^{p,n}_i$ is $p$-adapted for any $i \geq 1$. By construction,
$P^{n,p}_0$ is $\overline{C}_\idot([1]_+,\Z) = \Z$. To finish the
proof, we note that the pointed $n$-simplicial set
$\overline{X}^p_n$ is contractible by
Definition~\ref{ada.N.def}~\thetag{ii}, so that the right-hand side
of \eqref{p.n.p.eq} is an acyclic complex.  \endproof

We can now define our comparison functor from $\wDM(\Z,R)$ to the
category $\DA(R)$ of fixed points data. We will denote it by
\begin{equation}\label{nu.da.eq}
\nu:\wDM(\Z,R) \to \DA(R).
\end{equation}
By Lemma~\ref{no.P.le}, we are free to choose any adapted complex
$P_\idot$ to define the category of fixed points data $\DA(R)$; we
will use the one provided by Lemma~\ref{P.le}. For any $E \in
\wDM(\Z,R)$ and any integer $n \geq 1$, let
\begin{equation}\label{nu.E.eq}
\nu(E)_n = \overline{C}_\idot(X_n \wedge [\Z/n\Z]_+,E).
\end{equation}
The group $\Z/n\Z$ acts on $X_n \times [\Z/n\Z]_+$ via the second
factor, so that $\nu(E)_n$ is a complex of
$R[\Z/n\Z]$-modules. Taking together $\nu(E)_n$ for all $n \geq 1$,
we obtain a natural complex $\nu(E)$ in $\Fun(I,R)$.

To promote $\nu(E)$ to an object in $\DA(R)$, we need to construct
$(\Z/pn\Z)$-equivariant maps
$$
\alpha_{p,n}:\nu(E)_n \to P^{p,n}_\idot \otimes \nu(E)_{np}
$$
for all primes $p$ and integers $n \geq 1$. Note that the natural
inclusion $\iota:[1]_+ \to X^p_n$ defines $(\Z/pn\Z)$-equivariant
morphisms
$$
\iota_{p,n}:E_{np} = \overline{C}_\idot([Z/np\Z]_+,E) =
\overline{C}_\idot([1]_+ \wedge [Z/np\Z]_+,E) \to \nu(E)_{np},
$$
so that it suffices to construct morphisms
\begin{equation}\label{bar.a.eq}
\overline{\alpha}_{p,n}:\nu(E)_n \to P^{p,n} \otimes E_{np}
\end{equation}
and take $\alpha_{p,n} = \iota_{p,n} \circ
\overline{\alpha}_{p,n}$. By definition, the $\Z$-action on $X^p_n$
factors through $\Z/pn\Z$, and the quotient $(X_n^p \wedge
[\Z/pn\Z]_+)/(\Z/pn\Z)$ of the product $X_n^p \wedge [\Z/pn\Z]_+$ by
the diagonal $(\Z/pn\Z)$-action is of course isomorphic to
$X^p_n$. Thus the embedding $X^p_n \subset X_n$ induces a natural
map
$$
\begin{CD}
X_n^p \wedge [\Z/pn\Z]_+ @>{q}>> X^p_n @>>> X_n,
\end{CD}
$$
where $q$ is the quotient map. If we let $\Z$ act on $X^p_n \wedge
[\Z/pn\Z]_+$ via the second factor, then this map is
$\Z$-equivariant. Taking its product with the natural projection
$$
\begin{CD}
X^p_n \wedge [\Z/pn\Z]_+ @>>> [\Z/pn\Z]_+ @>>> [\Z/n\Z]_+,
\end{CD}
$$
we obtain a map
$$
a_{p,n}:X_n^p \wedge [\Z/pn\Z]_+ \to X_n \wedge [\Z/n\Z]_+.
$$
By definition, under this map, the diagonal $(\Z/pn\Z)$-action on
the left goes to the action of $\Z/pn\Z$ on the right via the second
factor, and conversely, the action via the second factor on the left
goes to the diagonal action on the right. To construct the morphism
$\overline{\alpha}_{p,n}$ of \eqref{bar.a.eq}, it remains to use the
natural identification
$$
\overline{C}_\idot(X_n^p \wedge [\Z/pn\Z]_+,E) \cong
\overline{C}_\idot(X_n^p,\Z) \otimes
\overline{C}_\idot([\Z/pn\Z]_+,E) \cong P^{p,n}_\idot
\otimes E_{pn},
$$
and take as $\overline{\alpha}_{p,n}$ the natural morphism
$a_{p,n}^*$ of \eqref{homo.f}.

We note that by construction, for any subgroup $H \subset \wh{\Z}$
with the quotient $W=\wh{\Z}/H$, the composition $\phi^H \circ \nu$
of the comparison functor $\nu$ of \eqref{nu.da.eq} with the
restriction functor $\phi^H$ of \eqref{phi.h.a} factors through the
fixed points functor $\Phi^H$, so that $\nu$ induces a natural
functor
\begin{equation}\label{nu.da.h}
\nu:\wDM(W,R) \to \DA_H(R)
\end{equation}
such that $\nu \circ \Phi^H \cong \phi^H \circ \nu$.

\begin{prop}\label{nu.prop}
\begin{enumerate}
\item If $H = n\wh{\Z} \subset \wh{\Z}$ is a cofinite subgroup, then
  the functor $\nu$ of \eqref{nu.da.h} is an equivalence of categories.
\item For any subgroup $H \subset \wh{\Z}$ with the quotient $W =
  \wh{\Z}/H$, and any integer $l$, let $\DA_H(R)^{\leq l} \subset
  \DA_H(R)$ be the full subcategory spanned by objects $M$ such that
  $\nu(M)$ lies in $\D^{\leq l}(I,R) \subset \D(I,R)$. Then the
  functor $\nu$ of \eqref{nu.da.h} induces an equivalence between
  $\wDM^{\leq l}(W,R)$ and $\DA_H(R)^{\leq l}$.
\end{enumerate}
\end{prop}

\proof{} As in the proof of Lemma~\ref{phi.add},
Proposition~\ref{ada.prop} provides an isomorphism of
functors
$$
\wPhi \cong \eta \circ \nu:\wDM(\Z,R) \to \D(I,R),
$$
where $\nu$ is the functor \eqref{nu.da.eq}. It induces a base
change map
$$
\nu \circ \wt{R} \to \rho,
$$
and by comparing Lemma~\ref{R.sq.le} and Lemma~\ref{phi.tate.le}, we
see that this map is also an isomorphism.  For any subgroup $H
\subset \wh{\Z}$, the functor $\wPhi$ induces a functor
$\wPhi:\wDM(W,R) \to \D_H(I,R)$, the functor $\wt{R}$ induces a
right-adjoint functor $\wt{R}:\D_H(I,R) \to \wDM(W,R)$, and we have
induced isomorphisms $\wPhi \cong \eta \circ \nu$, $\nu \circ \wt{R}
\cong \rho$.

Assume first that the subgroup $H$ is cofinite, that is, $H=
n\wh{\Z} \subset \wh{\Z}$ for some $n \geq 1$, and $W = \Z/n\Z$.
Denote by $\wDM(W,R)^{fr} \subset \wDM(W,R)$, $\DA_H(R)^{fr} \subset
\DA_H(R)$ the full triangulated subcategories generated by images of
the functors $\overline{R}$ resp. $\rho$. Then by adjunction, $\nu$
induces a fully faithful functor
$$
\wDM_H(\Z,R)^{fr} \to \DA_H(R)^{fr}.
$$
But Lemma~\ref{R.sq.le} implies that $\DA_H(R)^{fr} = \DA_H(R)$, so
that $\nu$ is essentially surjective, and by
Lemma~\ref{phi.dot.le}~\thetag{ii}, the functor $\wPhi:\wDM(W,R) \to
\D_H(I,R)$ is conservative, so that the same argument shows that
$\wDM(W,R)^{fr} = \wDM(W,R)$. Therefore $\nu$ is fully faithful on
the whole of $\wDM(W,R)$. This proves \thetag{i}.

For \thetag{ii}, note that that for a cofinite subgroup $H \subset
\wh{\Z}$, the statement immediately follows from
Lemma~\ref{phi.dot.le}~\thetag{i}. In the general case, it remains
to apply Proposition~\ref{can.filt.der.prop}.
\endproof

To see how Proposition~\ref{nu.prop} implies the $p$-typical
decomposition of Subsection~\ref{p.typ.subs}, note that for any $l
\geq 1$, the functor $\gamma_l$ of \eqref{psi.l} induces a functor
$\gamma_l:I \to I$, the pullback functor $\gamma_l^*$ lifts to a
natural functor
$$
\gamma_l^*:\DA(R) \to \DA(R),
$$
and we have $\gamma_l \circ \nu \cong \nu \circ \Psi^l$. Moreover,
fix a prime $p$, and take $H = \Z'_p$, so that $I^H \subset I$ is
the union of groupoids $\ppt_q \subset I$, $q = p^m$, $m \geq
0$. Then for any $l$ prime to $p$, the functor $\gamma_l$ can be
refined to a functor
$$
\wt{\gamma}_l:I^H \times \ppt_l \to I^{lH}
$$
sending $\ppt_q \times \ppt_l$ to $\ppt_{ql} \subset I^{lH}$, and
this lifts to a functor
\begin{equation}\label{wt.ga}
\wt{\gamma}_l^*:\DA_{lH}(R) \to \DA_H(R[\Z/l\Z]).
\end{equation}
We then have a natural isomorphism
\begin{equation}\label{nu.phi.p.l}
\nu \circ \Phi^{(p)}_{(l)} \cong \wt{\gamma}_l^* \circ \nu,
\end{equation}
where $\Phi^{(p)}_{(l)}$, $l$ prime to $p$ are the functors
\eqref{phi.p.l.der}.

\begin{corr}
Assume that the ring $R$ is $p$-local. Then for any integer $n$, the
functor
$$
\Phi^{(p)}_\idot:\wDM^{\leq n}(\Z,R) \to \prod_{l
  \text{ prime to }p}\wDM^{\leq n}(\Z_p,R[\Z/l\Z])
$$
of \eqref{phi.p.l.der} is an equivalence of categories.
\end{corr}

\proof{} Choose an adapted complex $P_\idot$, and let
$\overline{P}_\idot \subset P_\idot$ be the sum of its components
$P^{p,n}_\idot$ (in other words, we remove all the compoments
$P^{p',n}_\idot$ with $p' \neq p$). Then if we consider the category
$C^\alpha(\overline{P}_\idot,R)$, the only non-trivial components of
the map $\alpha$ of \eqref{alpha.eq} relate $M^n_\idot$, $n \geq 1$,
and $M^{pn}_\idot$. Therefore any object $\langle M,\alpha \rangle$
of $C^\alpha(\overline{P}_\idot,R)$ splits into a direct sum of
objects supported at
$$
\gamma_l(I^H) = \coprod_{m \geq 0}\ppt_{p^ml} \subset I
$$
over all $l \geq 1$ prime to $p$, and the pullback functors
$\wt{\gamma}_l^*$ of \eqref{wt.ga} induce an equivalence of
categories
$$
\DA(\overline{P}_\idot,R) \cong \prod_{l \text{ prime to }
  p}\DA_H(\overline{P}_\idot,R[\Z/l\Z]).
$$
By Proposition~\ref{nu.prop} and \eqref{nu.phi.p.l}, to finish the
proof, it remains to show that the natural functor
$\DA(\overline{P}_\idot,R) \to \DA(P_\idot,R) = \DA(R)$ induced by
the embedding $\overline{P}_\idot \subset P_\idot$ is an equivalence
of categories. The argument for this is exactly the same as in
Lemma~\ref{no.P.le}: one simply needs to observe that since the
ring $R$ is $p$-local, the only non-trivial Tate cohomology terms in
the right-hand side of \eqref{phi.tate} are those with our fixed prime $p$.
\endproof

{\small\noindent
Affiliations (in the precise form required for legal reasons):
\begin{enumerate}
\renewcommand{\labelenumi}{\arabic{enumi}.}
\item Steklov Mathematics Institute, Algebraic Geometry Section
  (main affiliation).
\item Laboratory of Algebraic Geometry, National Research University
Higher\\ School of Economics.
\end{enumerate}}

\noindent
{\em E-mail address\/}: {\tt kaledin@mi-ras.ru}


\begin{thebibliography}{BHM}

\bibitem[Ba]{clark} C. Barwick, {\em Spectral Mackey functors and
  equivariant algebraic K-theory (I)}, {\tt arXiv:1404.0108}.

\bibitem[BBD]{BBD} A. Beilinson, J. Bernstein, and P. Deligne, {\em
Faiscaux pervers}, Ast\'erisque {\bf 100}, Soc. Math. de France,
1983.

\bibitem[BZ]{ber} J. Bernstein and A. Zelevinski, {\em
  Representations of the group\\ GL(n,F), where F is a local
  non-Archimedean field}, Uspehi Mat. Nauk {\bf 31} (1976), 5-70.

\bibitem[Bo]{bo} M. B\"okstedt, {\em Topological Hochschild
  homology}, preprint, Bielefeld, 1985.

\bibitem[BHM]{BHM} M. B\"okstedt, W.C. Hsiang, and I. Madsen, {\em
    The cyclotomic trace and algebraic $K$-theory of spaces},
    Invent. Math.  {\bf 111} (1993), 465--539.

\bibitem[tD]{tD} T. tom Dieck, {\em Transformation groups}, De
Gruyter, Berlin-New York, 1987.

\bibitem[DTT]{TT} V. Dolgushev, D. Tamarkin, and B. Tsygan, {\em
  Noncommutative calculus and the Gauss-Manin connection}, in {\em
  Higher structures in geometry and physics}, 139-158, Progr. Math.,
  287, Birkh\"auser/Springer, New York, 2011.

\bibitem[DS]{dr-si} A. Dress and Ch. Siebeneicher, {\em The Burnside
  ring of profinite groups and the Witt vector construction},
  Adv. in Math. {\bf 70} (1988), 87-132.

\bibitem[G]{SGA} A. Grothendieck,  {\em Cat\'egories fibr\'ee et
  descente}, SGA I, Expos\'e VI, SMF, 2003.

\bibitem[H2]{hewi} L. Hesselholt, {\em Witt vectors of
  non-commutative rings and topological cyclic homology}, Acta
  Math. {\bf 178} (1997), 109--141.

\bibitem[HM]{HM} L. Hesselholt and I. Madsen, {\em On the $K$-theory
of finite algebras over Witt vectors of perfect fields}, Topology
{\bf 36} (1997), 29--101.

\bibitem[I]{ill} L. Illusie, {\em Complexe de de Rham-Witt et
  cohomologie cristalline}, Ann. Sci. \'Ecole Norm. Sup. (4) {\bf
  12} (1979), 501--661.

\bibitem[K1]{cartier} D. Kaledin, {\em Non-commutative Hodge-to-de
  Rham degeneration via the method of Deligne-Illusie}, Pure
  Appl. Math. Q. {\bf 4} (2008), 785--875.

\bibitem[K2]{ka-icm} D. Kaledin, {\em Motivic structures in
  non-commutative geometry}, Proc. ICM 2010.

\bibitem[K3]{mackey} D. Kaledin, {\em Derived Mackey functors},
  Mosc. Math. J. {\bf 11} (2011), 723--803.

\bibitem[K4]{ka-wi} D. Kaledin, {\em Witt vectors as a polynomial
  functor}, Selecta Math., {\bf 24} (2018), 359--402.

\bibitem[K5]{ka-hw} D. Kaledin, {\em Hochschild-Witt complex},
  Adv. in Math. {\bf 351} (2019), 33--95.

\bibitem[LMS]{may1} L.G. Lewis, J.P. May, and M. Steinberger, {\em
Equivariant stable homotopy theory}, with contributions by
J. E. McClure, Lecture Notes in Mathematics, {\bf 1213},
Springer-Verlag, Berlin, 1986.

\bibitem[L]{lind} H. Lindner, {\em A remark on Mackey functors},
Manuscripta Math.  {\bf 18} (1976), 273--278.

\bibitem[M]{may2} J.P. May, {\em Equivariant homotopy and cohomology
theory}, with contributions by M. Cole, G. Comezana, S. Costenoble,
A.D. Elmendorf, J.P.C. Greenlees, L.G. Lewis, Jr., R.J. Piacenza,
G. Triantafillou, and S. Waner, CBMS Regional Conference Series in
Mathematics, {\bf 91}. Published for the Conference Board of the
Mathematical Sciences, Washington, DC, by the American Mathematical
Society, Providence, RI, 1996.

\bibitem[P]{pira.dk} T. Pirashvili, {\em Dold-Kan type theorem for
  $\Gamma$-groups}, Math. Ann. {\bf 318} (2000), 277--298.

\bibitem[T]{the} J. Thevenaz and P. Webb, {\em The structure of
  Mackey functors}, Trans. AMS {\bf 347} (1995), 1865--1961.

\bibitem[V]{V} V. Voevodsky, {\em Triangulated categories of motives
  over a field}, in {\em Cycles, transfers, and motivic homology
  theories}, 188-238, Ann. of Math. Stud., {\bf 143}, Princeton
  Univ. Press, Princeton, NJ, 2000.

\end{thebibliography}
\end{document}